\documentclass[a4paper,11pt]{amsart}
\usepackage[margin=1.5cm]{geometry}
\usepackage{amssymb,amsmath,latexsym,amsthm,mathrsfs} 
\usepackage{aliascnt}
\usepackage{resmes}
\usepackage{graphicx}
\usepackage{enumerate}
\usepackage{comment}
\usepackage{aurical}

\makeatletter
\renewcommand{\l@section}{\@tocline{1}{0pt}{1pc}{20pt}{}} 
\makeatother

\usepackage[autostyle,italian=guillemets]{csquotes}
\usepackage[
backend=biber,
style=alphabetic, 
]{biblatex}
\usepackage{xcolor}
\definecolor{dark_orange}{HTML}{CC6633}
\definecolor{dark_purple}{HTML}{70004B}
\definecolor{brick}{RGB}{192, 57, 43}
\definecolor{aqua}{RGB}{22, 160, 133}

\usepackage[
hidelinks,
colorlinks=true,
urlcolor=black,
citecolor=teal,
linkcolor=dark_orange,
breaklinks=true]
{hyperref}
\usepackage[noabbrev,capitalize,nameinlink]{cleveref}

\numberwithin{equation}{section}

\newtheorem{theorem}{Theorem}[section]
\newtheorem{lemma}[theorem]{Lemma}
\newtheorem{cor}[theorem]{Corollary}
\newtheorem{prop}[theorem]{Proposition}
\theoremstyle{definition}
\newtheorem{definition}[theorem]{Definition}
\newtheorem{rem}[theorem]{Remark}
\newtheorem{example}[theorem]{Example}

\newcommand{\R}{\mathbb{R}}
\newcommand{\N}{\mathbb{N}}
  
\renewcommand{\S}{\mathbb{S}}

\newcommand{\norm}[1]{\|#1\|}
\newcommand{\intpart}[1]{\lfloor #1 \rfloor}

\newcommand{\cont}{\mathcal{C}}
\newcommand{\eqp}[1]{\stackrel{#1}{=}}
\newcommand{\degpar}{\operatorname{deg}_{\operatorname{par}}}

\newcommand{\strongeq}{\mathrel{\ooalign{\hfil$\equiv$\hfil\cr\hfil$\circ$\hfil}}}
\newcommand{\eqps}[1]{\stackrel{#1}{\strongeq}}

\newcommand{\TF}{\text{\Fontamici{TF}}}

\newcommand{\loc}{\rm{loc}}
\renewcommand{\leq}{\leqslant}
\renewcommand{\le}{\leqslant}
\renewcommand{\geq}{\geqslant}
\renewcommand{\ge}{\geqslant}
\renewcommand{\epsilon}{\varepsilon}
\allowdisplaybreaks

\title[Solutions of the space-fractional heat equation with polynomial growth]{Solutions of the space-fractional heat equation \\ with polynomial growth:\\
Liouville theorems and Schauder estimates}
\date{}

\author{\href{https://research-repository.uwa.edu.au/en/persons/serena-dipierro/}{Serena Dipierro}}
\address{
	Serena Dipierro,
	Department of Mathematics and Statistics, 
	The University of Western Australia, 
	35 Stirling Highway, Perth, WA 6009, Australia}
\email{serena.dipierro@uwa.edu.au}

\author{\href{https://www.math.columbia.edu/~savin/}{Ovidiu Savin}}
\address{
	Ovidiu Savin,
	Department of Mathematics,
	Columbia University, 
	2990 Broadway, New York, NY 10027, USA}
\email{savin@math.columbia.edu}

\author{\href{https://research-repository.uwa.edu.au/en/persons/enrico-valdinoci/}{Enrico Valdinoci}}
\address{
	Enrico Valdinoci,
	Department of Mathematics and Statistics, 
	The University of Western Australia, 
	35 Stirling Highway, Perth, WA 6009, Australia}
\email{enrico.valdinoci@uwa.edu.au}

\author{\href{https://research-repository.uwa.edu.au/en/persons/riccardo-villa-2/}{Riccardo Villa}} 
\address{
	Riccardo Villa,
	Department of Mathematics and Statistics, 
	The University of Western Australia, 
	35 Stirling Highway, Perth, WA 6009, Australia}
\email{riccardo.villa@research.uwa.edu.au}

\begin{document}
	
	\maketitle
	\begin{abstract}
		With the aim of extending to parabolic equations the framework developed by the first three authors in [Dipierro,~Savin, and~Valdinoci, \emph{Rev. Mat. Iberoam.}~(2019)], we introduce a notion of (viscosity) solutions for the space-fractional heat equation modulo~$k$. In turn, we are able to define $s$-caloric functions with polynomial growth at infinity, showing also related Schauder estimates, a Liouville-type result, and a stability property with respect to the $L^\infty_{\rm{loc}}$-convergence. 
	\end{abstract}
	
	\setcounter{tocdepth}{1}
	{
		\hypersetup{linkcolor=black}
		\tableofcontents
	}

	\section{Introduction}
	
In this paper, we extend the definition of viscosity solutions of a space-fractional heat equation to functions with polynomial growth at infinity. 
This can be seen as the parabolic counterpart of the theory put forth in~\cite{elliptic_poly}
for the fractional elliptic case.
	
	We recall that, given $s\in(0,1)$, the fractional Laplacian of a function $u:\R^n\to\R$ is defined as
	\begin{equation} \label{eq::def_sfp}
		\begin{split}
			&(-\Delta)^s u(x) := \mbox{p.v.}\, c_{n,s} \int_{\R^n} \frac{u(x)-u(y)}{|x-y|^{n+2s}}\,dy\\
			&\qquad\qquad =\lim_{\epsilon\to0}\, c_{n,s} \int_{B_\epsilon^c(x)} \frac{u(x)-u(y)}{|x-y|^{n+2s}}\,dy .
		\end{split}
	\end{equation}
	where $\mbox{p.v.}$ denotes the principal value, and $c_{n,s}>0$ is a normalizing constant depending only on~$n$ and~$s$. We also use the notation~$B_\epsilon^c(x):=\R^n\setminus B_\epsilon(x)$.
	
	For a more detailed study of the fractional Laplacian and more general non-local operators, we refer to~\cite{MR2448308, MR3168912, MR3967804}.
	
	We stress that the definition in~\eqref{eq::def_sfp} is well-posed if:
	\begin{enumerate}[i)]
		\item\label{item::frac_laplace_reg} $u$ is sufficiently regular near $x$, e.g., $u$ is $\cont^\alpha$-regular at $x$, for some $\alpha>2s$;
		\item\label{item::frac_laplace_growth} $u$ satisfies the growth condition at infinity
		\begin{equation} \label{eq::frac_laplace_growth}
			\int_{\R^n} \frac{|u(y)|}{1+|y|^{n+2s}}\,dy < +\infty.
		\end{equation}
	\end{enumerate}
	
	A defining feature of non-local operators is that the behavior of a function at infinity may deeply affect the value of the operator. This makes the growth condition in~\eqref{eq::frac_laplace_growth} a genuine source of difficulties in the regularity theory
	of non-local equations and free boundary problems, for instance when considering rescalings or blow-up sequences.
	
	\begin{example}
		Let us consider a smooth function $u:\R^n\to\R$ such that
		\begin{equation*}
			u(x) :=
			\begin{cases}
				|x|^2&\quad\mbox{in }B_1,\\
				0&\quad\mbox{in }B_2^c.
			\end{cases}
		\end{equation*}
		Clearly, $u$ satisfies both the conditions in~\eqref{item::frac_laplace_reg} and~\eqref{item::frac_laplace_growth}, and indeed
		\begin{equation*}
\left|(-\Delta)^s u(0)\right| = 
\int_{\R^n} \frac{u(y)}{|y|^{n+2s}} \,dy = \int_{B_1} \frac{dy}{|y|^{n+2s-2}}<+\infty.
		\end{equation*}
		
		However, the blow-up sequence $\{u_\lambda(x):=\lambda u(x/\lambda)\}_{\lambda>0}$ violates~\eqref{eq::frac_laplace_growth} as $\lambda\to+\infty$.
	\end{example}
	
In this spirit, several recent developments in fractional elliptic regularity theory have benefited from the notion of harmonic functions ``up to a polynomial'' (see, e.g.,~\cite{MR4038144, MR4385144, MR4395952, MR4413304, MR4808714, MR4850027, MR4884561, MR5023029, 2025arXiv251013340D}). Moreover, the regularity results established in~\cite{elliptic_poly} introduced several necessarily sharp refinements for handling the polynomial case, making them particularly well suited to delicate regularity questions, even in settings where polynomial growth is not explicitly assumed (see, e.g.,~\cite{MR3770173, MR4577340}).
	 
	The aim of this work is therefore to introduce an alternative notion of (viscosity) solution for non-local parabolic equations, in order to make sense of the fractional heat operator under a milder growth assumption.
	
	\begin{rem}
		We stress that, in the classical setting, the lack of any growth assumption allows one to consider solutions with polynomial growth (and, in fact, any growth) at infinity. For instance, the function~$u(t,x):=2nt+|x|^2$ satisfies
		\begin{equation*}
			\partial_t u - \Delta u = \partial_t (2nt) - \Delta(|x|^2) = 0.
		\end{equation*}
		
		By analogy, one would like to say that the function $v(t,x):=c_{n,s}t-|x|^{2s}$, for a suitable constant $c_{n,s}>0$ depending only on $n$ and $s$, satisfies
		\begin{equation*}
			\partial_t v + (-\Delta)^s v = 0.
		\end{equation*}
	\end{rem}
	
	Let $u:(-\infty,0]\times\R^n\to\R$ be a function with polynomial growth at infinity in space. The main idea is to consider a sequence of smooth cut-off functions $\{\eta_R\}_R$ such that the fractional Laplacian of $\eta_R u$ is well-defined. Then, we obtain that the space-fractional heat operator of~$u$,
	that is, 
	\begin{equation*}
		\partial_t u + (-\Delta)^s u,
	\end{equation*}
	diverges as~$R\to+\infty$, and it can be  expressed as a given function plus a (non-uniquely determined) polynomial remainder. We will make a precise sense of this in Definition~\ref{def::k_viscosity_sol} below.
	\medskip
	
	To conclude this introductory part, we present a possible application of our notion of solution for non-local parabolic equations. 
	
	The problem addressed in this paper stems from a question by Felix Otto: if a smooth function solves a homogeneous parabolic equation modulo a polynomial in space-time and has polynomial growth, must it itself be a space-time polynomial? This question, as well as our notion of solution, was motivated by problems arising in singular 
	stochastic partial differential equations (SPDEs)
	and Hairer's regularity structures\footnote{		
		We recall that regularity structures are defined as triplets $(A,T,G)$ of:
		\begin{enumerate}[i)]
			\item a locally finite index set $A\subset\R$ bounded from below;
			\item a graded model space $T:=\bigoplus_{\alpha\in A} T_\alpha$, where each $T_\alpha$ is a Banach space representing elements of ``$\alpha$-homogeneity'' (generalizing the classical space of monomials);
			\item a structure group $G\subseteq \mbox{Aut}(T)$ of linear operators acting on $T$ such that, for every $\Gamma\in G$, $\alpha\in A$, and $a\in T_\alpha$, we have that
			$$ \Gamma a - a \in \bigoplus_{\beta < \alpha} T_\beta .$$
		\end{enumerate}
	}, where polynomials play a distinguished role. 
	
	Indeed, in the framework of ``subcritical'' semilinear SPDEs, Hairer's theory proposes a generalized notion of Taylor expansion that splits functions into a ``rough'' term and a regular polynomial remainder. 	
	Due to the local subcriticality, non-linear fluctuations remain controllable at a small scale. This implies that solutions of the semilinear equation resemble those of the corresponding linear problems, and can be constructed via a suitable approximation and limit procedure from a regularized problem.
	
	We refer to~\cite{MR3289027,MR3274562, MR3971357, MR3925533, MR4683413} for both introductory overviews and detailed discussions on regularity structures and related topics.
	
	
	Building on the formalism of solutions ``up to a polynomial'', 
	we answer Felix Otto's question with
	our main result, that is a {\bf Liouville-type theorem} which establishes that any function that is $s$-caloric, up to a polynomial, on all of time-space and has polynomial growth must be a polynomial both in space and time (see Theorem~\ref{th::liouville} below).
	
	This Liouville-type result also relies on 
	a sharp version, valid ``up to polynomials'', of the
	{\bf Schauder estimates} for solutions to a space-fractional heat equation (see
	Theorem~\ref{th::schauder_modk} below).
	
	A useful structural result is also
	that {\bf solutions modulo~$k$ are stable} under
	a uniform control of the mass and locally uniform convergence
	(see the forthcoming Theorem~\ref{th::stability}).
	
	The extension of the elliptic case, introduced in~\cite{elliptic_poly}, to the parabolic setting requires significant additional arguments. Indeed, besides the technical transition from elliptic to parabolic regularity theory, a bespoke analysis is needed to account for the time dependence of the spatial polynomial representation of the solutions, in particular to distinguish between general time-dependent coefficients in the spatial polynomials and genuinely time-dependent polynomial coefficients.
	
	
	\section{Notation and main results}
	
	For our purposes, let us recall that if $\alpha:=(\alpha_1,\dots,\alpha_n)\in\N^n$ is a multi-index, for every $x\in\R^n$, we write
	\begin{equation*}
		x^\alpha := x_1^{\alpha_1}\dots x_n^{\alpha_n}.
	\end{equation*}
	Moreover, we say that a multi-index $\alpha\in\N^n$ is of degree at most $k$, and we simply write $|\alpha|\leq k$, if
	\begin{equation*}
		|\alpha| := \alpha_1 + \dots + \alpha_n \leq k.
	\end{equation*}
	
	Then, for a general vector space $W$ (finite-dimensional or not), let us define 
	\begin{align*}
		W_k^{ [x] } &:= \left\{\mbox{$P(x)$ is a polynomial of degree at most $k$ in $x$ with coefficients in $W$}\right\} \\
		&= \left\{\sum_{|\alpha|\leq k} c_\alpha x^\alpha \mbox{ with }c_\alpha\in W\right\}\\
		\mbox{and}\quad& W^{ [x] } := \bigcup_{k\in\N} W_k^{ [x] }.
	\end{align*}
	
	We also adopt the conventions
	\begin{equation} \label{eq::negative_deg}
		\begin{split}
			&W_0^{ [x] }:=W,\\
			\mbox{and}\quad&W_k^{ [x] }:=\{0_W\},\quad\mbox{ for all }k<0.
		\end{split}
	\end{equation}
	
	In what follows, we will focus on
	\begin{equation*}
		\R_k^{ [t,x] } := \left\{P(t,x):=\sum_{j+|\alpha|\leq k} c_{j,\alpha} t^j x^\alpha \mbox{ with }c_{j,\alpha}\in\R\right\}
	\end{equation*}
	and, for any interval $I\subseteq\R$ and for any $\sigma>0$,
	\begin{equation*}
		\left(\cont^{\sigma}(I)\right)_k^{ [x] } := \left\{P(t,x):=\sum_{|\alpha|\leq k} c_\alpha(t) x^\alpha \mbox{ with }c_\alpha\in\cont^{\sigma}(I) \right\}.
	\end{equation*}
	Namely, 
	$\left(\cont^{\sigma}(I)\right)_k^{ [x] }$
	is the family of polynomials
in~$x$ with degree at most~$k$ and whose coefficients are $\cont^{\sigma}(I)$-functions of~$t$,
while~$\R_k^{ [t,x] }$ is the smallest space in which the coefficients are monomials in~$t$, with the additional constraint that the degree of the total polynomial in~$x$ and~$t$
has degree at most~$k$.

	Furthermore, given $\gamma>0$, we conform to the notation $\cont^\gamma=\cont^{\intpart{\gamma},\{\gamma\}}$ for the space of all $\cont^{\intpart{\gamma}}$ functions with $\{\gamma\}$-H\"older continuous highest-order derivatives. Here, as customary, $\intpart{\gamma}$ and $\{\gamma\}$ denote the integer part and the fractional part of $\gamma$, respectively. 
	
	Throughout the rest of this paper, we consider a fixed $\sigma\in(0,1)$. 
	\medskip
	
	Since we work within the framework of viscosity solutions for parabolic equations (see Definition~\ref{def::viscosity_sol} below), we recall the definition of the standard parabolic cylinder of radius $r>0$ and height $r^{2s}$, namely
	\begin{equation*}
		Q_r :=(-r^{2s},0]\times B_r,
	\end{equation*}
	extending also this notation to $Q_\infty:=(-\infty,0]\times \R^n$.
	
Besides, for later convenience, let us introduce the closed, restricted cylinder
	\begin{equation*}
		\overline{Q}^*_r :=\left[-r^{2s},-\frac{r^{2s}}{2}\right]\times \overline{B}_r.
	\end{equation*}
	
	Moreover, let us consider, for any $k\in\N$, the measure $\omega_k$ on $\R^n$, which is absolutely continuous with respect to the Lebesgue measure, and is defined as
	\begin{equation*}
		d\omega_k(y):= \frac{dy}{1+|y|^{n+2s+k}}.
	\end{equation*}
We also set
	\begin{equation*}
		L^1(d\omega_k) := \left\{v:\R^n\to\R \mbox{ s.t. } \norm{v}_{L^1(d\omega_k)}:=\int_{\R^n} |v(y)|\,d\omega_k(y)<+\infty\right\}.
	\end{equation*}
	
	Then, for any $\xi\geq 0$ and for any interval $I \subseteq \R$, we define the space $V_k^\xi(I)$ 
	of all functions $u: I\to L^1(d\omega_k)$ such that, for every $\overline{t}\in I$,
	\begin{equation*}
		\lim_{t\to\overline{t}} \norm{u(t,\cdot)-u(\overline{t},\cdot)}_{L^1(d\omega_k)} = 0 ,
	\end{equation*}
	and 
	\begin{equation} \label{eq::t_poly_growth}
		\int_{\R^n} \frac{|D^{\intpart{\xi}}_t u(t_1,y)-D^{\intpart{\xi}}_t u(t_2,y)|}{|t_1-t_2|^{\{\xi\}}}\, d\omega_k(y) < +\infty.
	\end{equation}

	In other words,
	\begin{equation*}
		V_k^\xi(I) = \{ u\in\cont^0(I;L^1(d\omega_k)) \mbox{ s.t. $u$ satisfies~\eqref{eq::t_poly_growth}}\}.
	\end{equation*}
	For the rest of the present paper, we assume $\xi\geq \sigma$.
	
	In the wake of~\cite{elliptic_poly}, given $T>0$, $r>0$, and a function $u:(-T,0]\to L^1(d\omega_k)$, we define the $k$-tail of~$u$ in the interval $(-T,0]$ as
	\begin{equation} \label{eq::tail_def}
		J^k_{T,r}(u) := \sup_{t\in(-T,0]} \int_{B_{r}^c} \frac{|u(t,y)|}{|y|^{n+2s+k}}\,dy.
	\end{equation}
	For simplicity, we denote $J^k_R(u):=J^k_{R^{2s},R/2}(u)$.
	\medskip

	Now, we recall the definition of viscosity solutions for the space-fractional heat equation together with the definition of a suitable space of test functions.
	\begin{definition}[{\cite[Definition~2.5]{parabolic_viscosity}}]
		We define the space of test functions $\TF=\TF(Q_1)$ as the set of all pairs $(\varphi,(t-\tau,t]\times B_r(x))$ such that
		\begin{align*}
			&\varphi\in \cont^0((-1,0];L^1(d\omega_0))\cap \cont^2(Q_1)\\
			\mbox{and}\quad&(t-\tau,t]\times B_r(x)\subseteq Q_1.
		\end{align*}
	\end{definition}
	
	\begin{definition}[{\cite[Definition~2.8]{parabolic_viscosity}}] \label{def::viscosity_sol}
		An upper semi-continuous (lower semi-continuous) function $u:Q_1\to\R$ is said to be a viscosity sub-solution (super-solution) of
		\begin{equation} \label{eq::viscosity_heat}
			\partial_t u + (-\Delta)^s u = f \quad\mbox{in }Q_1, 
		\end{equation}
		and we write
		\begin{equation*}
			\begin{split}
				&\partial_t u + (-\Delta)^s u \leq f \quad\mbox{in }Q_1\\
				&\quad\big(\mbox{resp. }\partial_t u + (-\Delta)^s u \geq f\quad\mbox{in }Q_1\big),
			\end{split}
		\end{equation*}
		if, for every $(\varphi,(t-\tau,t]\times B_r(x))\in\TF$ touching $u$ by above (by below) at $(t,x)$, i.e. 
		\begin{enumerate}[i)]
			\item $\varphi(t,x)=u(t,x)$;
			\item $\varphi>u$ in $\left((t-\tau,t]\times B_r(x)\right)\setminus\{(t,x)\}$ (resp. $\varphi<u$), 
		\end{enumerate}
		we have 
		\begin{equation*}
			\begin{split}
				&\partial_t \varphi(t,x) + (-\Delta)^s \varphi(t,x) \leq f(t,x)\\
				&\quad\big(\mbox{resp. }\partial_t \varphi(t,x) + (-\Delta)^s \varphi(t,x) \geq f(t,x)\big).
			\end{split}
		\end{equation*}
		\smallskip
		
		We say that a function $u\in\cont^0(Q_1)$ is a viscosity solution for
		\begin{equation*} 
			\partial_t u + (-\Delta)^s u = f \quad\mbox{in }Q_1, 
		\end{equation*}
		if it is both a viscosity sub- and a viscosity super-solution.
	\end{definition}
	
	Let $\{\eta_R:\R^n\to\R\}_{R>0}$ be a family of 
	cut-off functions such that
	\begin{equation*}
		\eta_R(x):=
		\begin{cases}
			1\quad&\mbox{in }B_R,\\
			0\quad&\mbox{in }B_{R+1}^c.
		\end{cases}
	\end{equation*}
	
	In the wake of~\cite{elliptic_poly}, we give the following definition (compare with~\cite[Definition~1.1]{elliptic_poly}):
	\begin{definition}[Viscosity solution modulo~$k$] \label{def::k_viscosity_sol}
		Let $s\in(0,1)$, $k\in\N$, and $\xi\geq\sigma$. Let also $u:Q_\infty\to\R$, and~$f:Q_1\to\R$. Suppose that $u\in\cont^0\left(Q_1\right)\cap V_k^\xi((-1,0])$.
		
		We say that $u$ is a viscosity solution modulo a polynomial of degree $k$ of~\eqref{eq::viscosity_heat} (or simply a solution modulo~$k$), and we write
		\begin{equation} \label{eq::viscosity_heat_modk}
			\partial_t u + (-\Delta)^s u  \eqp{k} f  \quad\mbox{in }Q_1
		,\end{equation}
		if there exist a family of functions  $\{f_R:Q_1\to\R\}_{R>0}$ and a family of polynomials $\{P_R\}_{R>0}$ such that 
		\begin{align}
			\label{eq::k_viscosity_f_ptw} &\lim_{R\to+\infty} f_R(t,x) = f(t,x),\quad\mbox{for all }(t,x)\in Q_1,\\
			\label{eq::k_viscosity_poly} &P_R \in\left(\cont^{\sigma}((-1,0])\right)_{k-1}^{ [x] }, \quad\mbox{for all }R>0,\\
			\label{eq::k_viscosity_cutoff_sol} \mbox{and}\quad&\partial_t (\eta_Ru) + (-\Delta)^s (\eta_Ru) = f_R+P_R \quad\mbox{in $Q_1$, in the viscosity sense.}
		\end{align}
		
		If~$f=0$, namely
		\begin{equation*}
			\partial_t u + (-\Delta)^s u  \eqp{k} 0  \quad\mbox{in }Q_1,
		\end{equation*}
		then we say that $u$ is $s$-caloric modulo~$k$ in $Q_1$.
		
	\end{definition}

	\begin{rem} \label{rem::mod_poly}
		Notice that $f$ in Definition~\ref{def::k_viscosity_sol} is defined up to a polynomial of degree at most $k-1$, namely if~$\partial_t u + (-\Delta)^s u  \eqp{k} f$ then~$\partial_t u + (-\Delta)^s u  \eqp{k} f+P$, for any~$P \in \left(\cont^{\sigma}((-1,0])\right)_{k-1}^{ [x] }$.
		
Indeed, if $P\in \left(\cont^{\sigma}((-1,0])\right)_{k-1}^{ [x] }$, then we can write
		\begin{equation*}
			f_R+P_R = f_R + P_R \pm P = \widetilde{f}_R + \widetilde{P}_R,
		\end{equation*} 
		where 
		\begin{equation*}
			\widetilde{f}_R := f_R + P \to f+P,
		\end{equation*}
		as $R\to+\infty$, and
		\begin{equation*}
			\widetilde{P}_R := P_R-P \in \left(\cont^{\sigma}((-1,0])\right)_{k-1}^{ [x] }.
		\end{equation*}
		
		We refer to~\cite[Lemma~1.2]{elliptic_poly} for a proof of this for elliptic equations.
	\end{rem}
	
	In light of Remark~\ref{rem::mod_poly}, we introduce the following notation for the H\"older seminorm modulo a polynomial (see~\cite{elliptic_poly} for the elliptic case). Let $R>0$, $k\in\N$, and $\gamma>0$, then, for any $f\in\cont^0((-R^{2s},0];\cont^\gamma(B_R))$, we use the short notation
	\begin{equation} \label{eq::def_seminorm_mod_poly}
		\sup_{t\in(-R^{2s},0]}[f]_{\cont^\gamma(B_R;k)} = \sup_{t\in(-R^{2s},0]} \inf_{P\in\left(\cont^{\sigma}((-R^{2s},0])\right)_{k-1}^{ [x] }} [f-P]_{\cont^\gamma(B_R)}.
	\end{equation}
In line with~\cite[formula~(1.2)]{elliptic_poly}, when~$\gamma>k-1$ we have that
\begin{equation}\label{eq::def_seminorm_mod_poly22}
 \sup_{t\in(-R^{2s},0]}[f]_{\cont^\gamma(B_R;k)}=\sup_{t\in(-R^{2s},0]}[f]_{\cont^\gamma(B_R)}.\end{equation}
	\medskip
	
	As main result, we obtain sharp Schauder estimates in space for viscosity solutions to a space-fractional heat equation modulo~$k$.
	
	\begin{theorem}[Schauder estimates for solutions to a space-fractional heat equation modulo~$k$] \label{th::schauder_modk}
		Let $s\in(0,1)$, $k\in\N$,~$\xi\geq\sigma$, and~$u:Q_\infty\to\R$. Suppose that $u\in\cont^0\left(Q_1\right)\cap V_k^\xi((-1,0])$ is a viscosity solution modulo~$k$ of~\eqref{eq::viscosity_heat} with initial datum $\overline{u}\in L^\infty(B_1)\cap L^1(d\omega_k)$, i.e. 
		
		\begin{equation*}
			\begin{cases}
				\partial_t u + (-\Delta)^s u \eqp{k} f \quad&\mbox{in }Q_1,\\
				u(-1,\cdot) = \overline{u}\quad&\mbox{in }\R^n.
			\end{cases}
		\end{equation*}
		
		Then, for any $\gamma>0$ such that $\gamma\not\in \N$ and $\gamma+2s\not\in \N$, there exists a constant $C>0$ depending only on $n$, $s$, $k$, and $\gamma$ such that 
		\begin{equation} \label{eq::schauder}
			\sup_{t\in(-1/2,0]}\left(\norm{u}_{\cont^{\gamma+2s}(B_{1/2})}\right)
			\leq C\left(J^k_1(u) + \sup_{t\in(-1,0]}[f]_{\cont^{\gamma}(B_1;k)}+\norm{u}_{L^\infty(\overline{Q}^*_1)}\right).
		\end{equation} 
	\end{theorem}          
	
	We point out that, in the absence of any additional assumption on the external datum (for $t\in(-1,0]$), we cannot obtain more than H\"older-regularity in time for solutions of~\eqref{eq::viscosity_heat}, even for $k=0$,~$f=0$, and~$u(-1,\cdot)\equiv0$. For an explicit counterexample, we refer to~\ref{sec::counterex} (see also \cite[Subsection~2.4.1]{MR3115838}).
	
	On the other hand, assuming H\"older-regularity in time of the external datum and the forcing term, the solution enjoys higher interior time-regularity and, in turn, this will be useful when proving Theorem~\ref{th::schauder_modk}:
	
	\begin{prop} \label{prop::1_zeta_time_reg}
		Let $s\in(0,1)$, $k\in\N$,~$\xi\geq\sigma$, and~$\alpha\in(0,1)$. Let also~$u:Q_\infty\to\R$ be such that $u\in\cont^0\left(Q_1\right)\cap V_k^\xi((-1,0])$, and $u$ a viscosity solution modulo~$k$ of~\eqref{eq::viscosity_heat} with forcing term $f\in \cont^0(B_1;\cont^\alpha((-1,0]))$, and initial datum $\overline{u}\in L^\infty(B_1)\cap L^1(d\omega_k)$, i.e. 
		
		\begin{equation} \label{eq::init_value_problem_timereg}
			\begin{cases}
				\partial_t u + (-\Delta)^s u \eqp{k} f \quad&\mbox{in }Q_1,\\
				u(-1,\cdot) = \overline{u}\quad&\mbox{in }\R^n.
			\end{cases}
		\end{equation}
		
		Then, $u\in \cont^0(B_1;\cont^{1+\min\{\alpha,\sigma\}}((-1,0]))$, where $\sigma$ is as in Definition~\ref{def::k_viscosity_sol}.
	\end{prop}            
	
	We emphasize that the assumption in~\eqref{eq::t_poly_growth} preserves the time-regularity of the right-hand side in~\eqref{eq::init_value_problem_timereg}. Indeed, suppose for instance that $u$ satisfies
	$$ \partial_t u + (-\Delta)^s u \eqp{0} f \quad\mbox{in }Q_1, $$
	for some $f\in \cont^0(B_1;\cont^\alpha((-1,0]))$. Then, for all $R>1$, we have
	$$ \partial_t u + (-\Delta)^s (u \eta_R) = f_R \quad\mbox{in }Q_1, $$
	where
	$$ f_R(t,x) := f(t,x) + \int_{B_R^c} \frac{u(t,x)-u(t,y)}{|x-y|^{n+2s}}\, dy.$$
	Consequently,~\eqref{eq::t_poly_growth} yields $f_R\in \cont^0(B_1;\cont^{\min\{\alpha,\sigma\}}((-1,0]))$.
	
	We prove show Theorem~\ref{th::schauder_modk} and Proposition~\ref{prop::1_zeta_time_reg} in Section~\ref{sec::proof_schauder}.
	\medskip
	
	Another interesting result is that solutions of~\eqref{eq::viscosity_heat} in the sense of Definition~\ref{def::k_viscosity_sol} are stable under local uniform convergence.
	
	\begin{theorem}[Stability for solutions to a space-fractional heat  equation modulo~$k$] \label{th::stability}
		Let $s\in(0,1)$, $k\in\N$, and~$\xi\geq\sigma$. Let us consider sequences~$\{u_j:Q_\infty\to\R\}_j$ and $\{f_j:Q_\infty\to\R\}_j$ such that, for all $j$, the function~$u_j\in\cont^0\left(Q_1\right)\cap V_k^\xi((-1,0])$ is a viscosity solution modulo~$k$ of~\eqref{eq::viscosity_heat} with forcing term $f_j$ in $Q_1$. Moreover, suppose that the sequence $\{u_j\}_j$ has uniformly bounded $L^1(d\omega_k)$-energy. Namely,  
		
		\begin{align}
			\label{eq::heat_eq_modk_j} &\partial_t u_j + (-\Delta)^s u_j \eqp{k} f_j  \quad\mbox{in $Q_1$, for all $j$},\\
			\label{eq::unif_energy_bound} \mbox{and}\quad&\sup_{j\in\N}\, \sup_{t\in (-1,0]} \int_{\R^n} |u_j(t,y)|\, d\omega_k (y)<+\infty .
		\end{align}
		
		Suppose that there exist functions $u\in\cont^0((-1,0];L^1(\R^n))$ and $f:Q_1\to\R$ such that, as $j\to+\infty$,
		\begin{align}
			\label{eq::local_unif_conv_u} &u_j\to u\quad\mbox{in }L^\infty_{\loc}(Q_1),\\
			\label{eq::local_unif_conv_f} &f_j\to f\quad\mbox{in }L^\infty_{\loc}(Q_1),\\
			\label{eq::weak_conv_u} \mbox{and}\quad&\int_{B_1^c} \frac{u_j(t,y)\,\varphi(y)}{|y|^{n+2s+k}}\,dy \to \int_{\R^n}\frac{ u(t,y)\,\varphi(y)}{|y|^{n+2s+k}}\,dy\mbox{, for all $t\in(-1,0]$ and $\varphi \in L^\infty(B_1)$}.
		\end{align}
		
		Then, 
		\begin{equation} \label{eq::limit_eq}
			\partial_t u + (-\Delta)^s u \eqp{k} f  \quad\mbox{in }Q_1. 
		\end{equation} 
	\end{theorem}
	
	A proof of Theorem~\ref{th::stability} can be found in Section~\ref{sec::proof_stability}.
	\medskip
	
	We also consider solutions of the space-fractional heat equation in the whole half-space $Q_\infty$.
	\begin{definition} \label{def::viscosity_sol_mod_k_whole_space}
 		Let $s\in(0,1)$, $k\in\N$, $u:Q_\infty\to\R$, and~$f:Q_\infty\to\R$.  Let also
 		\begin{equation} \label{eq::def_d_ks}
 			d(k,s):=
 			\begin{cases}
 				k,\quad&\mbox{if }s\in\left(0,\frac{1}{2}\right],\\
 				k+1,\quad&\mbox{if }s\in\left(\frac{1}{2},1\right).
 			\end{cases}
 		\end{equation}
 		
 		Suppose that $u\in\cont^0\left(Q_\infty\right) \cap V_k^\xi((-\infty,0])$ satisfies the growth condition
 		\begin{equation} \label{eq::ancient_sol_0}
 			\norm{u}_{L^\infty({Q}_R)} \leq CR^{d(k,s)},\quad\mbox{for all }R>0.
 		\end{equation}
 		
 		We say that $u$ is a viscosity solution modulo~$k$ of~\eqref{eq::viscosity_heat} in $Q_\infty$ (or that $u$ is a global solution modulo~$k$), and we write
 		\begin{equation*} 
 			\partial_t u + (-\Delta)^s u  \eqp{k} f  \quad\mbox{in }Q_\infty,
 		\end{equation*}
 		if 
 		\begin{equation*} 
 			\partial_t u + (-\Delta)^s u  \eqp{k} f  \quad\mbox{in }Q_R,
 		\end{equation*}
 		for all $R>0$.
	\end{definition} 
	
	\begin{rem}
		The growth condition in~\eqref{eq::ancient_sol_0} is the usual control as $t\to-\infty$ imposed for the so-called ancient solutions. Compare it with~\cite[formula~(1.3)]{MR3987724} for the classical setting and~\cite[formula~3.2]{MR3385173} for the non-local case with sub-critical growth.
	\end{rem}
	
	Our third main result establishes a rigidity property. In particular, we prove a Liouville-type theorem, which guarantees that any globally $s$-caloric function must be a polynomial in space and time.
	
	\begin{theorem}[Liouville theorem for $s$-caloric functions modulo~$k$] \label{th::liouville}
		Let $s\in(0,1)$ and~$k\in\N$.		
		Let also~$\xi\in[\sigma,k]$ and $u\in\cont^0\left(Q_\infty\right) \cap V_k^\xi((-\infty,0])$.
		Suppose that $u$ satisfies the growth condition
		\begin{equation*} 
			\norm{u}_{L^\infty(Q_R)} \leq CR^{d(k,s)},\quad\mbox{for all }R>0.
		\end{equation*}
		
		Then, $u$ is globally $s$-caloric modulo~$k$, i.e.
		\begin{equation*}
			\partial_t u + (-\Delta)^s u \eqp{k} 0 \quad\mbox{in }Q_\infty,
		\end{equation*}
		if and only if $u=P$, with $P\in\left(\cont^{1+\sigma}((-\infty,0])\right)^{[x]}_{d(k,s)}$.
	\end{theorem}
	
	A proof of Theorem~\ref{th::liouville} will be presented in Section~\ref{sec::proof_liouville}.
	\medskip
	
	To conclude our analysis of viscosity solutions modulo~$k$, we show that such solutions exist. In particular, in Section~\ref{sec::existence}, we establish the existence to a parabolic initial-value problem modulo~$k$ with Dirichlet external condition. The precise statement goes as follows. 
	
	\begin{prop}[Existence of viscosity solutions modulo~$k$] \label{prop::existence}		
		Let $s\in(0,1)$,~$k\in\N$, and~$\xi\geq\sigma$. Let also~$\overline{u}:Q_\infty\to\R$ such that $\overline{u}\in V_k^\xi((-1,0])$,~$f\in \cont^0(B_1;\cont^\alpha((-1,0]))$, and~$v\in L^\infty(B_1)\cap L^1(d\omega_k)$.
		
		Then, there exists a function $u:Q_\infty\to\R$ such that
		\begin{equation} \label{eq::dirichlet_modk}
			\begin{cases}
				\partial_t u + (-\Delta)^s u  \eqp{k} f \quad &\mbox{in }Q_1,\\
				u=\overline{u}\quad &\mbox{in }B_1^c\times(-1,0],\\
				u(\cdot,0) = v\quad &\mbox{in }\R^n.
			\end{cases}
		\end{equation}
	\end{prop}

	\subsection{Viscosity solution strongly modulo~$k$}
	We point out that the condition in~\eqref{eq::k_viscosity_poly} of Definition~\ref{def::k_viscosity_sol} is quite general, allowing for any $\cont^\sigma$ dependence in time. On the other hand, restricting $P_R$ in Definition~\ref{def::k_viscosity_sol} to be a polynomial in both space and time allows us to strengthen the rigidity result in Theorem~\ref{th::liouville}. With this aim, we introduce the notion of viscosity solutions \emph{strongly} modulo~$k$ and we establish that any $s$-caloric function \emph{strongly} modulo~$k$ is a space-time polynomial itself.
	
	\begin{definition}[Viscosity solution strongly modulo~$k$] \label{def::k_viscosity_sol_strong}
		Let $s\in(0,1)$, $k\in\N$, $u:Q_\infty\to\R$, and~$f:Q_1\to\R$. Suppose that $u\in\cont^0\left(Q_1\right)\cap V_k^\xi((-1,0])$.
		
		We say that $u$ is a viscosity solution strongly modulo polynomial of degree $k$ of~\eqref{eq::viscosity_heat}, and we write
		\begin{equation} \label{eq::viscosity_heat_modk_strong}
			\partial_t u + (-\Delta)^s u  \eqps{k} f  \quad\mbox{in }Q_1,
		\end{equation}
		if there exist a family of functions  $\{f_R:Q_1\to\R\}_{R>0}$ and a family of polynomials $\{P_R\}_{R>0}$ such that 
		\begin{align}
			\label{eq::k_viscosity_f_ptw2} &\lim_{R\to+\infty} f_R(t,x) = f(t,x),\quad\mbox{for all }(t,x)\in Q_1,\\
			\label{eq::k_viscosity_poly_strong} &P_R \in \R_{k-1}^{ [t,x] }, \quad\mbox{for all }R>0,\\
			\label{eq::k_viscosity_cutoff_sol2} \mbox{and}\quad&\partial_t (\eta_Ru) + (-\Delta)^s (\eta_Ru) = f_R+P_R \quad\mbox{in $Q_1$, in the viscosity sense.}
		\end{align}
		
		Moreover, we say that $u$ is $s$-caloric strongly modulo~$k$ in $Q_1$, if
		\begin{equation*}
			\partial_t u + (-\Delta)^s u  \eqps{k} 0  \quad\mbox{in }Q_1.
		\end{equation*}
		\smallskip
		
		Analogously to Definition~\ref{def::viscosity_sol_mod_k_whole_space}, we say that $u$ satisfies
		\begin{equation*} 
			\partial_t u + (-\Delta)^s u  \eqps{k} f  \quad\mbox{in }Q_\infty
		\end{equation*}
		if 
		\begin{equation*} 
			\partial_t u + (-\Delta)^s u  \eqps{k} f  \quad\mbox{in }Q_R,
		\end{equation*}
		for all $R>0$.
	\end{definition}
	
	\begin{rem} \label{rem::equiv_mod_0}
		One can compare Definitions~\ref{def::k_viscosity_sol} and~\ref{def::k_viscosity_sol_strong}, noticing that
		assumption~\eqref{eq::k_viscosity_f_ptw2}
		coincides with~\eqref{eq::k_viscosity_f_ptw}, as well as assumption~\eqref{eq::k_viscosity_cutoff_sol2} is the same as~\eqref{eq::k_viscosity_cutoff_sol}. In fact, the difference of the two settings is that the polynomials in~\eqref{eq::k_viscosity_poly_strong}
		are assumed to be in a smaller family than those in~\eqref{eq::k_viscosity_poly}, see the framework in~\eqref{eq::negative_deg} and the comments below it.
		
		Accordingly, Definition~\ref{def::k_viscosity_sol_strong} is stronger than Definition~\ref{def::k_viscosity_sol} since
		\begin{equation*}
			\R_{k-1}^{ [t,x] }  \subset \left(\cont^{\sigma}((-1,0])\right)_{k-1}^{ [x] }.
		\end{equation*}
		Therefore,
		\begin{equation*}
			\partial_t u + (-\Delta)^s u \eqps{k} f \quad\implies\quad \partial_t u + (-\Delta)^s u \eqp{k} f.
		\end{equation*}
		
		Moreover, we point out that the two definitions are equivalent for $k=0$. Indeed,
		\begin{equation*}
			\left(\cont^{\sigma}((-1,0])\right)_{-1}^{ [x] } = \R_{-1}^{ [t,x] } = \{0\}
		\end{equation*}
		entails that
		\begin{equation*}
			\partial_t u + (-\Delta)^s u \eqp{0} f \quad\implies\quad \partial_t u + (-\Delta)^s u \eqps{0} f.
		\end{equation*}
	\end{rem}
	
	Our third main result establishes a rigidity property. In particular, we prove a Liouville-type theorem, which guarantees that any globally $s$-caloric function (strongly modulo~$k$) must be a polynomial in space and time.
	
	\begin{theorem}[Liouville theorem for $s$-caloric functions strongly modulo~$k$] \label{th::liouville_strong}
		Let $s\in(0,1)$ and~$k\in\N$.		
		Let also~$\xi\in[\sigma,k]$ and $u\in\cont^0\left(Q_\infty\right) \cap V_k^\xi((-\infty,0])$.
		Suppose that $u$ satisfies the growth condition
		\begin{equation*} 
			\norm{u}_{L^\infty(Q_R)} \leq CR^{d(k,s)},\quad\mbox{for all }R>0.
		\end{equation*}
		
		Then, $u$ is globally $s$-caloric strongly modulo~$k$, i.e.
		\begin{equation*}
			\partial_t u + (-\Delta)^s u \eqps{k} 0 \quad\mbox{in }Q_\infty,
		\end{equation*}
		if and only if $u=P+H$, with $P\in\R_k^{ [t,x] }$ and $H\in\R_{d(k,s)}^{ [x] }$. 
		
		In particular, $u$ is of degree at most $\lfloor\xi\rfloor$ in $t$.
	\end{theorem}
	A proof of Theorem~\ref{th::liouville_strong} will be presented in Section~\ref{sec::proof_liouville_strong}.
	\smallskip 
	
	We point out that the statement of Theorem~\ref{th::liouville_strong} can be reformulated introducing the notion of \emph{parabolic degree}.
	\begin{definition}\label{def::par_deg}
		Let $a\in\N$ and $b\in\N^n$. We define the parabolic degree of the monomial $t^a x^b$ as
		\begin{equation*}
		\degpar(t^ax^b) := 2sa+b.
		\end{equation*}
	\end{definition}
	In light of Definition~\ref{def::par_deg}, we rewrite Theorem~\ref{th::liouville_strong} as follows:
	\begin{cor}
		Let $s\in(0,1)$ and~$k\in\N$.		
		Let also~$\xi\in[\sigma,k]$ and $u\in\cont^0\left(Q_\infty\right) \cap V_k^\xi((-\infty,0])$.
		Suppose that $u$ satisfies the growth condition
		\begin{equation*} 
			\norm{u}_{L^\infty(Q_R)} \leq CR^{d(k,s)},\quad\mbox{for all }R>0.
		\end{equation*}
		
		Then, $u$ is globally $s$-caloric strongly modulo~$k$
		if and only if $u$ is a polynomial of parabolic degree at most $d(k,s)$ in both space and time, i.e., $u=P$, for some $P\in\R^{ [t,x] }$ with $\degpar(P)\leq d(k,s)$. 
	\end{cor}

	
	\section{Preliminaries and known results}
	\label{sec::preliminary}
	This section is devoted to outlining some preliminary results and well-known facts from the current literature concerning viscosity solutions and the regularity of solutions to a space-fractional heat equation.
	
	To begin with, extending~\cite[Lemma~2.3]{elliptic_poly} to the parabolic case, we show that if $u$ satisfies
	\begin{equation}
		\partial_t u + (-\Delta)^s u \eqp{k} 0,\quad\mbox{in }Q_1,
	\end{equation}
	for some $k\in\N$, and $u=0$ in $\left((-1,0]\times B_{1/2}^c\right)\cup\left(\{-1\}\times\R^n \right)$, then $u$ must vanish identically.
	\begin{lemma} \label{lemma::max_principle}
		Let $R>r>0$ and $T>0$. Let also $P\in\left(\cont^{\sigma}((-T,0])\right)^{ [x] }$ and $u$ be a viscosity solution of 
		\begin{equation*}
			\partial_t u + (-\Delta)^s u = P\quad\mbox{in }(-T,0]\times B_{R}.
		\end{equation*}
		
		If $u=0$ in $\left((-T,0]\times B_{r}^c\right)\cup\left(\{-T\}\times\R^n \right)$, then $u\equiv0$.
	\end{lemma}
	\begin{proof}
		We argue inductively on the degree of $P$.
		
		Suppose first that $P\in\left(\cont^{\sigma}((-T,0])\right)_{-1}^{ [x] }$, i.e. $P\equiv 0$. Then, by the maximum principle for viscosity solutions of the space-fractional heat equation, we infer that $u\equiv0$.
		\smallskip
		
		Assume now that the conclusion of the lemma holds true for any polynomial of degree at most $k-1$, and consider $P\in\left(\cont^{\sigma}((-T,0])\right)_{k}^{ [x] }$. Let us take $r'$ and $R'$ such that $r<r'<R'<R$. Let  us also consider $\theta\in\R^n$ such that $|\theta|<\min\{r'-r,R-R'\}$. Then, define
		\begin{equation} \label{eq::finite_diff_theta_sigma}
			\begin{split}
				&u_{\theta}(t,x) := u(t,x+\theta)-u(t,x)\\
				\mbox{and}\quad&P_{\theta}(t,x) := P(t,x+\theta)-P(t,x).
			\end{split}
		\end{equation}
		
		Thus, we have
		\begin{equation*}
			\partial_t u_{\theta} + (-\Delta)^su_{\theta} = P_{\theta} ,\quad\mbox{in }(-T,0]\times B_{R'},
		\end{equation*}
		in the viscosity sense, and
		\begin{equation*}
			u_{\theta}=0 \quad\mbox{in }\left((-T,0]\times B_{r'}^c\right)\cup\left(\{-T\}\times\R^n \right).
		\end{equation*}
		Moreover, $P_{\theta}\in\left(\cont^{\sigma}((-T,0])\right)_{k-1}^{ [x] }$. 
		
		Therefore, by the inductive step, we infer that $u_{\theta}\equiv 0$, and hence 
		\begin{equation*}
			u(t,x) = g(t)\quad\mbox{in }(-T,0]\times B_{r'},
		\end{equation*}
		for some function $g$. Since $u=0$ in $\left((-T,0]\times B_{r}^c\right)\cup\left(\{-T\}\times\R^n \right)$ by assumption, we have that 
		\begin{equation*}
			g(t) = u(t,x) = 0,\quad\mbox{for all }t\in(-T,0].
		\end{equation*}
		Thus, we conclude that $u\equiv0$.
	\end{proof}
	\medskip
	
	Moving the term involving the time derivative to the right-hand side in \eqref{eq::viscosity_heat}, we obtain generalized versions of Theorem~3.1 and Corollary~3.2 in~\cite{elliptic_poly} to the parabolic case. 
	
	\begin{prop} \label{prop::parabolic_3_1}
		Let $s\in(0,1)$, $k\in\N$,~$\xi\geq\sigma$, and~$u:Q_\infty\to\R$. Suppose that $u\in\cont^1\left((-1,0];\cont^{\alpha}_{\loc}(B_1)\right)\cap V_k^\xi((-1,0])$, for some $\alpha>2s$.
		
		Let  $\eta:\R^n\to\R$ be a function with compact support and such that $\eta\equiv 1$ in $B_2$.
		
		Then, there exist a function $\psi: B_1\times B_2^c\to\R$ satisfying
		\begin{equation} \label{eq::estimate_psi}
			\sup_{x\in B_1,\, y\in B_2^c} |\partial_x^\zeta \psi(x,y)| < +\infty,
		\end{equation}
		for any $\zeta\in\N^n$, and a polynomial $P_{u,\eta}\in(\cont^{\sigma}((-1,0]))_{k-1}^{ [x] }$ such that
		\begin{equation*}
			\partial_t (\eta u) + (-\Delta)^s (\eta u) = P_{u,\eta} + f_{u,\eta},
		\end{equation*}
		where
		\begin{equation*}
			f_{u,\eta} := \partial_t (\eta u) + f_{1,u}+f_{2,u}+ f_{u,\eta}^\star
		\end{equation*}
		with
		\begin{equation} \label{eq::def_f1_f2_fstar}
			\begin{split}
				&f_{1,u}:= \int_{B_2} \frac{u(t,x)-u(t,y)}{|x-y|^{n+2s}}\,dy,\\
				&f_{2,u}:= \int_{B_2^c} \frac{u(t,x)}{|x-y|^{n+2s}}\,dy,\\
				\mbox{and}\quad&f_{u,\eta}^\star := \int_{B_2^c} \frac{\eta(y)u(t,y)\psi(x,y)}{|y|^{n+2s+k}}\,dy .
			\end{split}
		\end{equation}
	\end{prop}
	
	\begin{cor} \label{cor::parabolic_3_2}
		Let $s\in(0,1)$, $k\in\N$,~$\xi\geq\sigma$, and~$u:Q_\infty\to\R$. Suppose that $u\in\cont^1\left((-1,0];\cont^{\alpha}_{\loc}(B_1)\right)\cap V_k^\xi((-1,0])$, for some $\alpha>2s$.
		
		Let $\{\eta_R:\R^n\to\R\}_{R>0}$ be a sequence of functions such that
		\begin{align*}
			&\mbox{supp}(\eta_R)\subseteq B _R,\quad\mbox{for all }R>0, \\
			\mbox{and}\quad &\lim_{R\to+\infty}\eta_R=1\quad\mbox{a.e. in }\R^n. 
		\end{align*}
		
		Then, there exist a function $f_u$ and a family of polynomials $\{P_{u,\eta_R}\}_{R>0}\subset \left(\cont^{\sigma}((-1,0])\right)_{k-1}^{ [x] }$ such that
		\begin{equation*}
			\lim_{R\to+\infty} \partial_t(\eta_R u)+ (-\Delta)^s(\eta_R u)+P_{u,\eta_R} = f_u\quad\mbox{pointwise in }Q_1,
		\end{equation*}
		where
		\begin{equation*}
			f_{u} := \partial_t u + f_{1,u}+f_{2,u}+ f_{3,u},
		\end{equation*}
		with
		\begin{equation} \label{eq::def_f3} 
			f_{3,u}(t,x):=\int_{B_2^c} \frac{u(t,y)\psi(x,y)}{|y|^{n+2s+k}}\,dy
		\end{equation} and~$\psi$ satisfying~\eqref{eq::estimate_psi}.
	\end{cor}
	
	\begin{rem} \label{rem::parabolic_3_5}
		We point out that, by Definition~\ref{def::k_viscosity_sol} and Corollary~\ref{cor::parabolic_3_2}, we have that
		\begin{equation*}
			\partial_t u + (-\Delta)^s u \eqp{k} f_u\quad\mbox{in }Q_1,
		\end{equation*}
		whenever $u\in\cont^1\left((-1,0];\cont^{\alpha}_{\loc}(B_1)\right)\cap V_k^\xi((-1,0])$.
	\end{rem}
	
	As a direct consequence of Proposition~\ref{prop::parabolic_3_1}, arguing as in the proof of~\cite[Corollary~3.7]{elliptic_poly}, we obtain the following convergence result for $f_R$ in Definition~\ref{def::k_viscosity_sol}.
	\begin{cor} \label{cor::parabolic_3_7}
		Let $k\in\N$ and let $u$ and $f_R$ be as in Definition~\ref{def::k_viscosity_sol}.
		
		Then, for every $m\in\N$ and $R'>R$, we have that
		\begin{equation} \label{eq::parabolic_3_7}
			\sup_{t\in(-1,0]} \inf_{P\in\left(\cont^{\sigma}((-1,0])\right)_{k-1}^{ [x] }}\norm{D^m(f_{R'}-f_R-P)}_{L^\infty(B_1)} \leq C\sup_{t\in(-1,0]}  \int_{B_R^c} \frac{|u(t,y)|}{|y|^{n+2s+k}}\,dy, 
		\end{equation}
		for some constant $C>0$ depending only on $n$, $s$, and~$m$.
	\end{cor}
	
	For the sake of completeness, the proofs of these results can be found respectively in Appendices~\ref{sec::parabolic_3_1}, \ref{sec::parabolic_3_2}, and~\ref{sec::parabolic_3_7}.
	\medskip
	
	As a byproduct of Corollary~\ref{cor::parabolic_3_7}, we infer the following:
	\begin{cor} \label{cor::parabolic_3_8}
		Let $u\in\cont^1\left((-1,0];\cont^0(B_1)\right) \cap V_0^\xi((-1,0])$ . Then,
		\begin{equation} \label{eq::strong_mod_0}
			\partial_t u + (-\Delta)^s u \eqps{0} f \quad\mbox{in }Q_1
		\end{equation}
		if and only if
		\begin{equation} \label{eq::weak_mod_0}
			\partial_t u + (-\Delta)^s u \eqp{0} f \quad\mbox{in }Q_1
		\end{equation}
		if and only if 
		\begin{equation} \label{eq::viscosity_mod_0}
			\partial_t u + (-\Delta)^s u = f \quad\mbox{in }Q_1,
		\end{equation}
		in the viscosity sense.
	\end{cor}

	\begin{proof}
		Observe that the equivalence between~\eqref{eq::strong_mod_0} and~\eqref{eq::weak_mod_0} follows from Remark~\ref{rem::equiv_mod_0}. Thus, we focus on showing that also~\eqref{eq::weak_mod_0} and~\eqref{eq::viscosity_mod_0} are equivalent.
		\smallskip
		
		Let us suppose first that 
		\begin{equation*}
			\partial_t u + (-\Delta)^s u = f \quad\mbox{in }Q_1,
		\end{equation*}
		in the viscosity sense. Consider the sequence $\{\eta_R:\R^n\to\R\}_{R}$ defined as
		\begin{equation*}
			\eta_R :=
			\begin{cases}
				1\quad&\mbox{in }B_R,\\
				0\quad&\mbox{in }B_R^c,
			\end{cases}
		\end{equation*}
		for any $R>0$. Then, for every $(t,x)\in Q_1$ and every $R$ large enough (e.g $R>10$), we also have that
		\begin{equation} \label{eq::crop_viscosity}
			\partial_t (\eta_{R/2}u) + (-\Delta)^s (\eta_{R/2}u) = f +\int_{\R^n} \frac{(1-\eta_{R/2})u(t,y)}{|x-y|^{n+2s}}\,dy\quad\mbox{in }Q_1,
		\end{equation}
		in the viscosity sense.
		
		Let us define 
		\begin{equation*}
			v := (\eta_R-\eta_{R/2})u.
		\end{equation*}
		Notice that $\eta_R v= v$ and $v\equiv 0$ in $(-1,0]\times B_{R/2}$.
		
		Therefore, from Proposition~\ref{prop::parabolic_3_1} with $k=0$, we deduce that, in $Q_1$,
		\begin{equation} \label{eq::v_eq}
			\begin{split}
				&\partial_t v+ (-\Delta)^sv = (-\Delta)^s (\eta_R v) = f_{1,v}+f_{2,v}+f_{v,\eta_R}^* \\
				&\qquad= \int_{B_R\setminus B_2} \frac{v(t,y)\psi(x,y)}{|y|^{n+2s}}\,dy = \int_{B_R\setminus B_{R/2}} \frac{u(t,y)\psi(x,y)}{|y|^{n+2s}}\,dy .
			\end{split}
		\end{equation}
		
		From~\eqref{eq::crop_viscosity} and~\eqref{eq::v_eq}, we obtain that, in~$Q_1$,
		\begin{equation*}
			\begin{split}
				&\partial_t(\eta_R u) + (-\Delta)^s(\eta_R u)\\
				&\qquad= \partial_t [(\eta_R-\eta_{R/2})u] + (-\Delta)^s [(\eta_R-\eta_{R/2})u] + \partial_t (\eta_{R/2}u) + (-\Delta)^s (\eta_{R/2}u)\\
				&\qquad= \partial_t v+ (-\Delta)^sv + \partial_t (\eta_{R/2}u) + (-\Delta)^s (\eta_{R/2}u)\\
				&\qquad= \int_{B_R\setminus B_{R/2}} \frac{u(t,y)\psi(x,y)}{|y|^{n+2s}}\,dy + f +\int_{\R^n} \frac{(1-\eta_{R/2})u(t,y)}{|x-y|^{n+2s}}\,dy=:f_R.
			\end{split}
		\end{equation*}
		
		Moreover, recalling~\eqref{eq::estimate_psi}, we see that, in~$Q_1$,
		\begin{eqnarray*}
		&&	\lim_{R\to+\infty}f_R-f=
		\lim_{R\to+\infty} \left(\int_{B_R\setminus B_{R/2}} \frac{u(t,y)\psi(x,y)}{|y|^{n+2s}}\,dy + \int_{\R^n} \frac{(1-\eta_{R/2})u(t,y)}{|x-y|^{n+2s}}\,dy\right) = 0.\end{eqnarray*}
				Namely,
		\begin{equation*}
			\partial_t u + (-\Delta)^s u \eqp{0} f \quad\mbox{in }Q_1
		\end{equation*}
		in the sense of Definition~\ref{def::k_viscosity_sol} (with~$k=0$, and so~$P_R=0$).
		
		Suppose now that
		\begin{equation*}
			\partial_t u + (-\Delta)^s u \eqp{0} f \quad\mbox{in }Q_1.
		\end{equation*}
		Since $\left(\cont^{\sigma}((-1,0])\right)_{-1}^{ [x] }=\{0\}$, according to Definition~\ref{def::k_viscosity_sol}, there exists a sequence of functions $\{f_R\}_{R>0}$ such that
		\begin{align*}
			&f_R\to f,\quad\mbox{as }R\to+\infty\\
			\mbox{and}\quad& \partial_t(\eta_R u) + (-\Delta)^s(\eta_R u) = f_R \quad\mbox{in }Q_1,
		\end{align*}
		in the viscosity sense.
		
		Now, observe that $\eta_R u$ converges to $u$ locally uniformly, as $R\to+\infty$.
		
		Moreover, by Corollary~\ref{cor::parabolic_3_7} (used here with $m=0$ and~$k=0$), we have that, for every $(t,x)\in Q_1$,
		\begin{equation} \label{eq::f_conv_unif}
			\begin{split}
				& |f(t,x)-f_R(t,x)| = \lim_{R'\to+\infty} |f_{R'}(t,x)-f_R(t,x)| \\
				&\qquad\qquad \leq  \lim_{R'\to+\infty} \sup_{t\in(-1,0]} \norm{f_{R'}-f_R}_{L^\infty(B_1)} \leq C \sup_{t\in(-1,0]} \int_{B_R^c} \frac{|u(t,y)|}{|y|^{n+2s}}\,dy.
			\end{split}
		\end{equation}
		In particular, taking the limit as $R\to+\infty$, we infer from~\eqref{eq::f_conv_unif} that $f_R$ converges uniformly to $f$ in $Q_1$.
		
		As a consequence, we use~\cite[Theorem~5.6]{MR3115838} (see also~\cite[Lemma~5]{MR2781586} for the elliptic case) to deduce that 
		\begin{equation*}
			\partial_t u + (-\Delta)^s u = f \quad\mbox{in }Q_1,
		\end{equation*}
		in the viscosity sense, concluding the proof.
	\end{proof}
	\medskip
	
	Another important feature of Definition~\ref{def::k_viscosity_sol} is the consistency with respect to lower-order growth. 
	\begin{lemma} \label{lemma::parabolic_3_9}
		Let $s\in(0,1)$, $\xi\geq\sigma$, and~$k\in\N$. Let $u\in\cont^0\left(Q_1\right)\cap V_k^\xi((-1,0])$ be a viscosity solution modulo~$k$ of~\eqref{eq::viscosity_heat} with forcing term $f\in \cont^0(B_1;\cont^\xi((-1,0]))$, i.e. 
		
		\begin{equation} \label{eq::eq_u_k_f}
			\partial_t u + (-\Delta)^s u \eqp{k} f \quad\mbox{in }Q_1.
		\end{equation}
		
		Suppose that
		\begin{equation} \label{eq::tail_lower_order}
			\sup_{t\in(-1,0]} \int_{\R^n} |u(t,y)| \,d\omega_j(y) <+\infty,
		\end{equation}
		for some $j\leq k$.
		
		Then, there exists a polynomial $P\in\left(\cont^{\sigma}((-1,0])\right)_{k-1}^{ [x] }$ such that
		\begin{equation*}
			\partial_t u + (-\Delta)^s u \eqp{j} f+P \quad\mbox{in }Q_1.
		\end{equation*}
	\end{lemma}
		
	\begin{proof}
		Let $\{\eta_R\}_R$ be a sequence of  smooth cut-off functions such that $\eta_R\in\cont^\infty_c(B_{R+1})$ and~$\eta_R=1$ in~$B_R$, for all~$R\geq4$. Let us also consider functions $v:=(1-\eta_4)u$ and $w:=\eta_4 u$. In particular, $v$ is smooth in~$Q_1$ (being identically zero) and, by~\eqref{eq::tail_lower_order},
		\begin{equation*}
			\sup_{t\in(-1,0]} \int_{\R^n} |v(t,y)| \,d\omega_j(y) \leq \sup_{t\in(-1,0]} \int_{\R^n} |u(t,y)| \,d\omega_j(y) <+\infty.
		\end{equation*}
		Thus, Corollary~\ref{cor::parabolic_3_2} yields
		\begin{equation} \label{eq::heat_v_j}
			\begin{split}
				&\partial_t v(t,x) + (-\Delta)^s v(t,x) \eqp{j} f_v(t,x) \\
				&\qquad  = -\int_{B_2} \frac{v(t,y)}{|x-y|^{n+2s}}\,dy + \int_{B_2^c} \frac{v(t,y)\psi(x,y)}{|y|^{n+2s+j}}\,dy\\
				&\qquad= \int_{B_4^c} \frac{u(t,y)\psi(x,y)}{|y|^{n+2s+j}}\,dy,
			\end{split}
		\end{equation}
		where in the last step we have used that $v=0$ in $B_4$ and $v=u$ in $B_4^c$.
		
		Recalling Definition~\ref{def::k_viscosity_sol}, we obtain from~\eqref{eq::heat_v_j} that
		\begin{equation} \label{eq::heat_v_j_2}
			\partial_t (\eta_Rv)(t,x) + (-\Delta)^s (\eta_Rv)(t,x) = \int_{B_4^c} \frac{v(t,y)\psi(x,y)}{|y|^{n+2s+j}}\,dy + \varphi_R(t,x) + P_R(t,x)\quad\mbox{in }Q_1,
		\end{equation}
		in the viscosity sense, for some function $\varphi_R\to0$, as $R\to+\infty$, and for some polynomial $P_R \in\left(\cont^{\sigma}((-1,0])\right)_{j-1}^{ [x] }\subset \left(\cont^{\sigma}((-1,0])\right)_{k-1}^{ [x] }$.
		
		On the other hand, by~\eqref{eq::eq_u_k_f} and Definition~\ref{def::k_viscosity_sol}, we deduce that
		\begin{equation} \label{eq::eq_u_k_f_2}
			\partial_t (\eta_Ru)(t,x) + (-\Delta)^s (\eta_Ru)(t,x) = f(t,x) + \theta_R(t,x) + H_R(t,x)\quad\mbox{in }Q_1,
		\end{equation}
		in the viscosity sense, for some function\footnote{
			From Definition~\ref{def::k_viscosity_sol}, we know that there exists a sequence of functions $\{f_R\}_R$ such that $f_R\to f$, as $R\to+\infty$. Therefore, it suffices to set $\theta_R:= f_R-f$.
		} $\theta_R\to0$, as $R\to+\infty$, and for some polynomial $H_R \in \left(\cont^{\sigma}((-1,0])\right)_{k-1}^{ [x] }$.
		
		Therefore,~\eqref{eq::heat_v_j_2} and~\eqref{eq::eq_u_k_f_2} together entail that, in~$Q_1$,
		\begin{equation} \label{eq::heat_w_j}
			\begin{split}
				&\partial_t (\eta_4w)(t,x) + (-\Delta)^s (\eta_4w)(t,x) = \partial_t( \eta_Rw)(t,x) + (-\Delta)^s (\eta_Rw)(t,x) \\
				&\qquad= \partial_t\big( \eta_R(u-v)\big)(t,x) + (-\Delta)^s\big( \eta_R(u-v)\big)(t,x)\\
				&\qquad = f(t,x)-\int_{B_4^c} \frac{u(t,y)\psi(x,y)}{|y|^{n+2s+j}}\,dy + \theta_R(t,x)-\varphi_R(t,x) + H_R(t,x)-P_R(t,x).
			\end{split}
		\end{equation}
		Moreover, we have that
		\begin{equation} \label{eq::limit_poly}
			\lim_{R\to+\infty} \varphi_R-\theta_R+H_R-P_R = \lim_{R\to+\infty} H_R-P_R =: P 
		\end{equation}
		is a polynomial of degree at most $k-1$.
		
		Now, let $g:=f+P$ and $\phi_R:=\theta_R+H_R-P_R-P$. Then, from~\eqref{eq::eq_u_k_f_2} and~\eqref{eq::limit_poly}, we obtain
		\begin{equation*}
			\partial_t( \eta_Ru) + (-\Delta)^s( \eta_Ru) = g+\phi_R+P_R,
		\end{equation*}
		in the viscosity sense, with $\phi_R\to0$, as $R\to+\infty$, and $P_R \in \left(\cont^{\sigma}((-1,0])\right)_{j-1}^{ [x] }$. Hence, it follows that
		\begin{equation*}
			\partial_t u + (-\Delta)^s u \eqp{j} g \quad\mbox{in }Q_1.
		\end{equation*}
		as desired.
	\end{proof}
	
	The results obtained so far allow us to reduce to the case of a ``classical'' space-fractional heat equation, satisfied by~$u$, which is the solution modulo~$k$, multiplied by a given cut-off function. The details are as follows:
		
	\begin{cor} \label{cor::parabolic_3_10}
		Let $s\in(0,1)$, $\xi\geq\sigma$, and~$k\in\N$. Let $u\in\cont^0\left(Q_1\right)\cap V_k^\xi((-1,0])$ be a viscosity solution modulo~$k$ of~\eqref{eq::viscosity_heat} with forcing term $f\in \cont^0(B_1;\cont^\xi((-1,0]))$, i.e. 
		
		\begin{equation*}
			\partial_t u + (-\Delta)^s u \eqp{k} f \quad\mbox{in }Q_1.
		\end{equation*}
		
		Let also $\rho\geq1$ and $$\chi_\rho:=\begin{cases}1 &\quad{\mbox{ in }} B_\rho,\\
		 0 &\quad{\mbox{ in }} B_\rho^c.\end{cases}$$
		
		Then, there exists a polynomial $P\in \left(\cont^{\sigma}((-1,0])\right)_{k-1}^{ [x] }$ such that
		\begin{equation} \label{eq::eq_viscosity_cutoff}
			\begin{split}
				&\partial_t (\chi_\rho u)(t,x) + (-\Delta)^s (\chi_\rho u)(t,x) \\
				&\qquad= f (t,x)+ P(t,x) + \int_{B_2\setminus B_\rho} \frac{u(t,y)}{|x-y|^{n+2s}}\,dy - \int_{B_2^c \cap B_\rho^c} \frac{u(t,y)\psi(x,y)}{|y|^{n+2s+k}} \,dy\quad\mbox{in }Q_1
			\end{split}
		\end{equation}
		in the viscosity sense, where $\psi$ is as in Proposition~\ref{prop::parabolic_3_1}.
	\end{cor}

	\begin{proof}
		Let $\{\eta_R\}_R$ be a sequence of  smooth cut-off function such that $\eta_R\in\cont^\infty_c(B_{R+1})$ and $\eta_R=1$ in $B_R$, for all $R\geq1$. By Definition~\ref{def::k_viscosity_sol}, we have that
		\begin{equation} \label{eq::eq_u_R}
			\partial_t (\eta_Ru)(t,x) + (-\Delta)^s (\eta_Ru)(t,x) = f(t,x) + \varphi_R(t,x) + P_R(t,x)\quad\mbox{in }Q_1,
		\end{equation}
		in the viscosity sense, for some function $\varphi_R\to0$, as $R\to+\infty$, and for some polynomial $P_R \in \left(\cont^{\sigma}((-1,0])\right)_{k-1}^{ [x] }$. 
		
		Let $w:=\chi_\rho u$ and $v:=(1-\chi_\rho)u$, for some $\rho\geq1$. Then, $v$ is smooth in $Q_1$ (being identically zero), and accordingly, thanks to Corollary~\ref{cor::parabolic_3_2}, we also have
		\begin{equation*}
			\begin{split}
				&\partial_t v(t,x) + (-\Delta)^s v(t,x) \eqp{k} f_v(t,x) \\
				&\qquad  = -\int_{B_2\setminus B_\rho} \frac{u(t,y)}{|x-y|^{n+2s}}\,dy + \int_{B_2^c\cap B_\rho^c} \frac{u(t,y)\psi(x,y)}{|y|^{n+2s+j}}\,dy\quad\mbox{in }Q_1.
			\end{split}
		\end{equation*} Notice that we are using the notation~$B_2\setminus B_\rho=\varnothing$ if~$\rho\ge2$.
		
		Hence, recalling Definition~\ref{def::k_viscosity_sol} once again, we infer that
		\begin{equation}  \label{eq::eq_v_R}
			\begin{split}
				&\partial_t( \eta_Rv)(t,x) + (-\Delta)^s (\eta_Rv)(t,x)\\
				&\qquad= -\int_{B_2\setminus B_\rho} \frac{u(t,y)}{|x-y|^{n+2s}}\,dy + \int_{B_2^c\cap B_\rho^c} \frac{u(t,y)\psi(x,y)}{|y|^{n+2s+j}}\,dy +\phi_R(t,x)+H_R(t,x)\quad\mbox{in }Q_1,
			\end{split}
		\end{equation}
		in the viscosity sense, with for some function $\phi_R\to0$, as $R\to+\infty$, and for some polynomial $H_R \in \left(\cont^{\sigma}((-1,0])\right)_{k-1}^{ [x] }$.
		
		Therefore, putting~\eqref{eq::eq_u_R} and~\eqref{eq::eq_v_R} together, we deduce that
		\begin{equation*}
			\begin{split}
				&\partial_t (\eta_Rw)(t,x) + (-\Delta)^s (\eta_Rw)(t,x)\\
				=& \partial_t\big( \eta_R(u-v)\big)(t,x) + (-\Delta)^s \big(\eta_R(u-v)\big)(t,x)\\
				=& f(t,x) + \int_{B_2\setminus B_\rho} \frac{u(t,y)}{|x-y|^{n+2s}}\,dy - \int_{B_2^c\cap B_\rho^c} \frac{u(t,y)\psi(x,y)}{|y|^{n+2s+j}}\,dy +\varphi_R(t,x)-\phi_R(t,x)+P_R(t,x)-H_R(t,x)
			\end{split}
		\end{equation*}
		in $Q_1$, in the viscosity sense, which reads
		\begin{equation} \label{eq::eqk_w}
			\partial_t w(t,x) + (-\Delta)^s w(t,x) \eqp{k}  f(t,x) + \int_{B_2\setminus B_\rho} \frac{u(t,y)}{|x-y|^{n+2s}}\,dy - \int_{B_2^c\cap B_\rho^c} \frac{u(t,y)\psi(x,y)}{|y|^{n+2s+j}}\,dy\quad\mbox{in }Q_1.
		\end{equation}
		
		Moreover, since $w$ is supported in $B_{\rho}$, we have that
		\begin{equation*}
			\sup_{t\in(-1,0]} \int_{\R^n} |w(t,y)|\, d\omega_0(y) \leq C \sup_{t\in(-1,0]} \int_{\R^n} |u(t,y)|\,d\omega_k(y) <+\infty.
		\end{equation*}
		
		From this,~\eqref{eq::eqk_w} and Lemma~\ref{lemma::parabolic_3_9}, we obtain that
		\begin{equation*}
			\begin{split}
				&\partial_t w(t,x) + (-\Delta)^s w(t,x) \\
				&\qquad\eqp{0}  f(t,x)+P(t,x) + \int_{B_2\setminus B_\rho} \frac{u(t,y)}{|x-y|^{n+2s}}\,dy - \int_{B_2^c\cap B_\rho^c} \frac{u(t,y)\psi(x,y)}{|y|^{n+2s+j}}\,dy\quad\mbox{in }Q_1,
			\end{split}
		\end{equation*}
		for some polynomial $P\in \left(\cont^{\sigma}((-1,0])\right)_{k-1}^{ [x] }$.
		
		Therefore, thanks to Corollary~\ref{cor::parabolic_3_8}, we conclude that
		\begin{equation*}
			\begin{split}
				&\partial_t w(t,x) + (-\Delta)^s w(t,x) \\&\qquad=  f(t,x)+P(t,x) + \int_{B_2\setminus B_\rho} \frac{u(t,y)}{|x-y|^{n+2s}}\,dy - \int_{B_2^c\cap B_\rho^c} \frac{u(t,y)\psi(x,y)}{|y|^{n+2s+j}}\,dy\quad\mbox{in }Q_1,
			\end{split}
		\end{equation*}
		in the viscosity sense, showing~\eqref{eq::eq_viscosity_cutoff}.
	\end{proof}

	
	\section{Schauder estimates for solutions to a space-fractional heat equation} 
	\label{sec::nonlocal_schauder}
	Here, we establish some Schauder estimates for solutions to the space-fractional heat equation, proving higher interior regularity under H\"older-continuous forcing terms. The setting considered, for the moment,
	is not ``up to a polynomial''.
	
	With this aim, we will employ the following version of interior Schauder estimates (compare with~\cite[Theorem~5.2]{parabolic_viscosity} and~\cite[Theorem~2.2]{MR3385173}):
	
	\begin{theorem}\label{th::rescaled_known_schauder}
		Let $s\in(0,1)$, $\xi\geq\sigma$, and $f\in\cont^0(Q_1)$. Let $u\in \cont^0(Q_1)\cap V_0^\xi((-1,0])$ be a viscosity solution of
		\begin{equation*}
			\partial_t u +(-\Delta)^s u = f \quad\mbox{in }Q_1.
		\end{equation*}
		
		Then, $u\in\cont^{\alpha}\left(Q_{1/2}\right)$, for some $\alpha\in(0,1)$, and
		
		\begin{equation} \label{eq::rescaled_known_schauder}
				[u]_{\cont^{\alpha}\left(Q_{1/2}\right)} \leq C\left(J^0_1(u)+\norm{f}_{L^\infty(Q_1)} +\norm{u}_{L^\infty(\overline{Q}^*_1)}\right).
		\end{equation}
	\end{theorem}
	
	We observe that, since $u$ and $f$ are uniformly bounded in time in $Q_1$, we argue as in~\cite[Lemma~5.2]{MR3661864} (see also~\cite[Lemma~4.1]{elliptic_poly}) to obtain the following auxiliary result\footnote{Compare~\eqref{eq::rescaled_known_schauder} with the conclusion of~\cite[Theorem~5.2]{parabolic_viscosity}.}.
	
	\begin{lemma} \label{lemma::partial_Linfty_norm_estimate}
		Let $s\in(0,1)$ and~$\xi\geq\sigma$. Let $u\in \cont^0(Q_1)\cap V_0^\xi((-1,0])$ be such that 
		\begin{align*}
			&J^0_1(u) \leq 1\\ 
			\mbox{and}\quad&\partial_t u +(-\Delta)^s u \leq 1 \quad\mbox{in }Q_1.
		\end{align*}
		
		Then, 
		\begin{equation*} 
				\norm{u}_{L^\infty(Q_1)} \leq C+\norm{u}_{L^\infty(\overline{Q}^*_1)}.
		\end{equation*}
	\end{lemma}
	
	\begin{proof}
		Since $u$ is uniformly bounded in time, we argue as in~\cite[Lemma~5.2]{MR3661864}. Therefore, up to rescalings, we assume that
		\begin{align*}
			&J^0_1(u) \leq \delta_0\\ 
			\mbox{and}\quad&\partial_t u +(-\Delta)^s u \leq \delta_0 \quad\mbox{in }Q_1,
		\end{align*}
		for some small $\delta_0>0$.
		Moreover, by~\eqref{eq::t_poly_growth}, there exist $a\geq0$ such that 
		\begin{equation*}
			|u(t,x)| \leq a(1+t)^\xi+\sup_{x\in B_1} |u(-1,x)|,
		\end{equation*}
		for every $(t,x)\in Q_1$.
		
		Thus, comparing $u$ with the barrier 
		\begin{equation*} \label{eq::Linfty_norm_barrier}
			\psi(t,x):=\frac{a(1+t)^\xi}{(1-|x|^2)^{n}}+\norm{u}_{L^\infty(\overline{Q}^*_1)},
		\end{equation*}
		we deduce the desired estimate. Indeed, up to vertical translations, we assume that the function $\psi_h:=\psi+h$ lies above $u$ and touches $u$ at some point $(t_0,x_0)\in Q_1$, namely
		\begin{align*}
			&u \leq \psi_h\quad\mbox{in }Q_1\\
			\mbox{and}\quad& u(t_0,x_0) = \psi_h(t_0,x_0)
		\end{align*} 
		
		Then, by the computations in~\cite[Lemma~5.2]{MR3661864} and recalling Definition~\ref{def::viscosity_sol}, we infer that
		\begin{equation*}
			\partial_t \psi_h(t_0,x_0) +(-\Delta)^s \psi_h(t_0,x_0) \geq c,
		\end{equation*}
		for some constant $c>0$. However, this leads to a contradiction whenever $\delta_0$ is small enough.
	\end{proof}
	
	A direct implication of Lemma~\ref{lemma::partial_Linfty_norm_estimate} is the following:
	
	\begin{cor} \label{cor::Linfty_norm_estimate}
		Let $s\in(0,1)$ and $\xi\geq\sigma$. Let $u\in \cont^0(Q_1)\cap V_0^\xi((-1,0])$ be such that 
		\begin{equation*}
			\partial_t u +(-\Delta)^s u = f \quad\mbox{in }Q_1.
		\end{equation*}
		
		Then,
		\begin{equation} \label{eq::Linfty_norm_estimate}
				\norm{u}_{L^\infty(Q_1)} \leq C\left( J^0_1(u) + \norm{f}_{L^\infty(Q_1)}  \right) +\norm{u}_{L^\infty(\overline{Q}^*_1)}.
		\end{equation}
	\end{cor}
	
	\begin{proof}[Proof of Theorem~\ref{th::rescaled_known_schauder}]
		Thanks to~\cite[Theorem~5.2]{parabolic_viscosity} and Corollary~\ref{cor::Linfty_norm_estimate}, we have that
		\begin{equation*}
			\begin{split}
				&[u]_{\cont^{\alpha}(Q_{1/2})} = \sup_{(t_1,x_1),(t_2,x_2)\in Q_1}\frac{|u(t_1,x_1)-u(t_2,x_2)|}{(|t_1-t_2|^{1/2s}+|x_1-x_2|)^\alpha}\\
				&\qquad\leq C\left(\norm{u}_{L^\infty(\overline{Q}_1)}+J^0_1(u)+\norm{f}_{L^\infty(Q_1)}\right) \leq C\left(J^0_1(u)+\norm{f}_{L^\infty(Q_1)}+\norm{u}_{L^\infty(\overline{Q}^*_1)}\right),
			\end{split}
		\end{equation*}
		showing~\eqref{eq::rescaled_known_schauder}.
	\end{proof}
	
Another consequence of Theorem~\ref{th::rescaled_known_schauder} is that we are able to estimate the full $\cont^m$-norm of the forcing term by the H\"older semi-norm of its highest-order derivative. The precise result reads as follows.

	\begin{cor} \label{cor::schauder_seminorm}
		Let $s\in(0,1)$,~$\xi\geq\sigma$, and~$r\in (0,9/10)$. Let also~$m\in\N$ and $\theta\in(0,1)$. Let~$u:Q_\infty\to\R$ and~$f\in\cont^{m+\theta}(Q_1)$.
		
		Suppose that $u\in\cont^0(Q_1)\cap V_0^\xi((-1,0])$ is a solution of
		\begin{equation*}
			\partial_t u +(-\Delta)^s u= f\quad\mbox{in }Q_1
		\end{equation*}
		in the viscosity sense.
		
		Then, there exists a constant $C>0$, depending only on $n$, $s$, $m$, and $\xi$, such that
		\begin{equation} \label{eq::schauder_seminorm}
			\begin{split}
				&\norm{u}_{L^\infty(Q_{9/10})}+\sup_{t\in(-99/100,0]} \norm{f}_{\cont^m(B_{99/100})}\\
				&\qquad\qquad \leq C\left(J^0_{1,r}(u)+\sup_{t\in(-1,0]}[f]_{\cont^{m+\theta}(B_1)}+\norm{u}_{L^\infty(\overline{Q}^*_1)}\right).
			\end{split}
		\end{equation}
	\end{cor} 
	The proof of Corollary~\ref{cor::schauder_seminorm} adapts to our setting the strategy implemented to establish~\cite[Lemma~4.1]{elliptic_poly}. For the convenience of the reader, it can be found in Appendix~\ref{sec::parabolic_4_1}.
	\medskip
	
	Next, we state a version of the Schauder estimates obtained by an inductive argument from~\cite[Theorem~5.1]{parabolic_viscosity}. Since the proof follows the ideas of~\cite[Lemma~4.2]{elliptic_poly}, we postpone the details to Appendix~\ref{sec::known_suboptimal_schauder}.
	
	 \begin{prop} \label{prop::known_suboptimal_schauder}
	 	Let $s\in(0,1)$, $\xi\geq\sigma$, and $f:Q_1\to\R$. Let $u\in \cont^0(Q_1)\cap V^\xi_0((-1,0])$ be a viscosity solution of
	 	\begin{equation*}
	 		\partial_t u +(-\Delta)^s u = f \quad\mbox{in }Q_1.
	 	\end{equation*}
	 	
	 	Then, for any $\gamma>0$ such that $\gamma\not\in\N$ and $\gamma+2s\not\in\N$,
	 	\begin{equation} \label{eq::known_suboptimal_schauder}
	 	\sup_{t\in(-1/2,0]}\norm{u(t,\cdot)}_{\cont^{\gamma+2s}(B_{1/2})}
	 			\leq C\left(J^0_1(u)+\sup_{t\in(-1,0]}\norm{f}_{\cont^\gamma(B_1)} +\norm{u}_{L^\infty(Q_1)}\right),
	 	\end{equation}
	 	for some constant $C>0$ depending only on $n$,~$s$, and~$\gamma$.
	 \end{prop}
	 
	From Corollary~\ref{cor::schauder_seminorm} and Proposition~\ref{prop::known_suboptimal_schauder}, we infer sharper Schauder estimates for solutions to the space-fractional heat equation.
	\begin{theorem} \label{th::ultimate_schauder}
		Let $s\in(0,1)$ and $\xi\geq\sigma$. Let $u:Q_\infty\to\R$ and $f:Q_1\to \R$.
		
		Suppose that $u\in\cont^0(Q_1)\cap V_0^\xi((-1,0])$ is a solution of 
		\begin{equation*}
				\partial_t u +(-\Delta)^s u= f\quad\mbox{in }Q_1
		\end{equation*}
		in the viscosity sense.
		
		Then, for any $\gamma>0$ such that $\gamma\not\in\N$ and $\gamma+2s\not\in\N$,
		\begin{equation} \label{eq::ultimate_schauder}
			\sup_{t\in(-1/2,0]}\norm{u(t,\cdot)}_{\cont^{\gamma+2s}(B_{1/2})} \leq C\left(J^0_1(u)+\sup_{t\in(-1,0]}[f]_{\cont^{\gamma}(B_1)}+\norm{u}_{L^\infty(\overline{Q}^*_1)}\right),
		\end{equation}
		for some constant $C>0$ depending only on $n$, $s$, and~$\gamma$.
	\end{theorem} 
	
	\begin{proof} 
		
		Up to resizing the domain, we assume that
		\begin{equation*}
			\partial_t u +(-\Delta)^s u= f\quad\mbox{in }Q_4.
		\end{equation*}
		
		Let $\eta\in\cont^\infty_c(B_3;[0,1])$ be a smooth cut-off function such that $\eta=1$ in $B_2$, and let $v(t,x):=\eta(x) u(t,x)$. Then, we have
		\begin{equation*}
			\partial_t v + (-\Delta)^s v =g\quad\mbox{in }Q_{3/2},
		\end{equation*}
		where $g:= f+h$ and
		\begin{equation*}
			h(t,x):=\int_{B_2^c} \frac{(1-\eta(y))u(t,y)}{|x-y|^{n+2s}}\,dy.
		\end{equation*}
		
		Thus, Proposition~\ref{prop::known_suboptimal_schauder} entails that
		\begin{equation} \label{eq::v_partial_schauder}
			\begin{split}
				&\sup_{t\in(-1,0]}\norm{v(t,\cdot)}_{\cont^{\gamma+2s}(B_{1})} 
				\leq C\left(J^0_{3/2}(v)+\sup_{t\in(-(3/2)^{2s},0]}\norm{g}_{\cont^\gamma(B_{3/2})} +\norm{v}_{L^\infty(Q_{3/2})}\right)\\
				&\qquad\qquad\leq C\left(J^0_{3}(v)+\sup_{t\in(-(3/2)^{2s},0]}\norm{f}_{\cont^\gamma(B_{3/2})}+\sup_{t\in(-(3/2)^{2s},0]}\norm{h}_{\cont^\gamma(B_{3/2})} +\norm{v}_{L^\infty(Q_{3/2})}\right),
			\end{split}
		\end{equation}
		up to relabeling $C$. Here, in the last step, we have also used that
		\begin{equation*}
			\begin{split}
				J^0_{3/2}(v)& = \sup_{t\in(-(3/2)^{2s},0]} \left(\int_{B_{3/2}^c} \frac{|v(t,y)|}{|y|^{n+2s}}\,dy + \int_{B_{3/2}\setminus B_{3/4}} \frac{|v(t,y)|}{|y|^{n+2s}}\,dy\right) \\
				&\leq J^0_{3}(v) + C\norm{v}_{L^\infty(Q_{3/2})}.
			\end{split}
		\end{equation*}
		
		Now, we estimate the four terms on the right hand side of~\eqref{eq::v_partial_schauder}.
		First, notice that, by Corollary~\ref{cor::schauder_seminorm}, we have that
		\begin{equation} \label{eq::v_Linfty_norm_estiamte}
			\begin{split}
				&\norm{v}_{L^\infty(Q_{3/2})}\leq \norm{u}_{L^\infty(Q_{3/2})} 
				 \leq C\left(J^0_3(u)+\sup_{t\in(-3^{2s},0]}[f]_{\cont^\gamma(B_3)}+\norm{u}_{L^\infty(\overline{Q}^*_3)}\right),
			\end{split}
		\end{equation}
		up to relabeling the constant $C>0$.
		
		Moreover, observe that, for all $m\in\N$,
		\begin{align*}
			&|D^m_xh(t,x)| \leq C\int_{B_2^c}\frac{|u(t,y)|}{|x-y|^{n+2s+m}}\,dy\\
			\mbox{and}\quad&|D^m_x h(t,x_1)-D^m_x h(t,x_2)| \leq C|x_1-x_2|\int_{B_2^c}\frac{|u(t,y)|}{|x_1-y|^{n+2s+m+1}}\,dy.
		\end{align*}
		It thereby follows that
		\begin{equation} \label{eq::h_holder_estimate}
			\begin{split}
				&\sup_{t\in(-(3/2)^{2s},0]}\norm{h}_{\cont^\gamma(B_{3/2})} \leq C\sup_{(t,x)\in Q_{3/2}}\left(\int_{B_2^c}\frac{|u(t,y)|}{|x-y|^{n+2s+m}}\,dy\right)\\
				&\qquad  \leq C\sup_{t\in (-(3/2)^{2s},0]}\left(\int_{B_{2}^c}\frac{|u(t,y)|}{|y|^{n+2s}}\,dy\right) \leq CJ^0_3(u),
			\end{split}
		\end{equation}
		for some positive constant $C$ depending only on $n$,~$s$, and~$m$, and possibly changing from line to line.
		
		Plugging~\eqref{eq::v_Linfty_norm_estiamte} and~\eqref{eq::h_holder_estimate} into~\eqref{eq::v_partial_schauder}, and recalling that $v=u$ in $Q_1$, we find
		\begin{equation*}
			\begin{split}
				&\sup_{t\in(-1,0]}\norm{u(t,\cdot)}_{\cont^{\gamma+2s}(B_{1})} = \sup_{t\in(-1,0]}\norm{v(t,\cdot)}_{\cont^{\gamma+2s}(B_{1})} \\
				&\qquad\qquad\leq C\left(J^0_3(u)+\sup_{t\in(-(3/2)^{2s},0]}\norm{f}_{\cont^\gamma(B_{3/2})}+\sup_{t\in(-3^{2s},0]}[f]_{\cont^\gamma(B_3)}+\norm{u}_{L^\infty(\overline{Q}^*_3)}\right),
			\end{split}
		\end{equation*}
		Employing again Corollary~\ref{cor::schauder_seminorm} to estimate $\displaystyle \sup_{t\in (-(3/2)^{2s},0]}\norm{f}_{\cont^\gamma(B_{3/2})}$, up to relabeling $C$, we conclude that
		\begin{equation*}
			\sup_{t\in(-1,0]}\norm{u(t,\cdot)}_{\cont^{\gamma+2s}(B_{1})} \leq C\left(J^0_3(u)+\sup_{t\in (-3^{2s},0]}[f]_{\cont^\gamma(B_3)}+\norm{u}_{L^\infty(\overline{Q}^*_3)}\right).
		\qedhere\end{equation*}
	\end{proof}
	
	\section{Proof of Theorem~\ref{th::schauder_modk} \\ Schauder estimates for solutions to a space-fractional heat equation modulo~$k$} 
	\label{sec::proof_schauder}
	In this section, we establish the Schauder estimates and higher time-regularity for viscosity solutions with polynomial growth stated in Theorem~\ref{th::schauder_modk}. 
	
	To this end, we first prove Proposition~\ref{prop::1_zeta_time_reg}, which relies on the following general time-regularity result.
	
	\begin{prop} \label{prop::schauder_time_reg}
		Let $s\in(0,1)$, $k\in\N$, $\alpha\in(0,1)$, and~$u:Q_\infty\to\R$. Suppose that $u\in\cont^0\left(Q_1\right)\cap \cont^0((-1,0]; L^1(d\omega_k))$ is a viscosity solution modulo~$k$ of~\eqref{eq::viscosity_heat} with forcing term $f\in \cont^0(B_1;\cont^\alpha((-1,0]))$, external datum $g\in \cont^0(B_1^c;\cont^\alpha((-1,0]))$, and initial datum $\overline{u}\in L^\infty(B_1)\cap L^1(d\omega_k)$, i.e. 
		
		\begin{equation*}
			\begin{cases}
				\partial_t u + (-\Delta)^s u \eqp{k} f \quad&\mbox{in }Q_1\\
				u=g\quad&\mbox{in }(-1,0]\times B_1^c,\\
				u(-1,\cdot) = \overline{u}\quad&\mbox{in }\R^n.
			\end{cases}
		\end{equation*}
		
		Then, $u\in \cont^0(B_1;\cont^{1+\min\{\alpha,\sigma\}}((-1,0]))$. 
	\end{prop}
	\begin{proof}
		Without loss of generality, we assume $\alpha\leq\sigma$, otherwise we just interchange roles.
		
		Recalling Definition~\ref{def::k_viscosity_sol}, we have that there exist a function $f_4$ and a polynomial $P_4\in\left(\cont^{\sigma}((-1,0])\right)_{k-1}^{ [x] }$ of degree at most $k-1$ such that
		\begin{equation*}
			\partial_t (\eta_4 u)+(-\Delta)^s(\eta_4 u) = f_4+P_4\quad\mbox{in }Q_1,
		\end{equation*}
		where $\eta_4\in\cont^\infty_c(B_5)$ is a smooth cut-off function such that $\eta_4=1$ in $B_4$.
		
		Then, since $\eta_4 u = u$ in $Q_1$, \cite[Theorem~2.3]{MR3628879} yields
		\begin{equation*}
			\begin{split}
				&\norm{\partial_t u}_{\cont^\alpha(Q_{1/2})} = \norm{\partial_t (\eta_4u)}_{\cont^\alpha(Q_{1/2})} \\
				\leq& C \left( J^0_1(\eta_4 u)+ \sup_{t_1,t_2\in (-1,0]} \int_{B_1^c} \frac{|(\eta_4g)(t_1,y)-(\eta_4g)(t_2,y)|}{|t_1-t_2|^\alpha}\,d\omega_0(y) + \sup_{x\in B_1} \norm{f_4+P_4}_{\cont^\alpha((-1,0])} +\norm{\eta_4 u}_{L^\infty(\overline{Q}_1^*)} \right)\\
				\leq& C \left( \sup_{t\in(-1,0]} \int_{B_5\setminus B_{1/2}} \frac{|u(t,y)|}{|y|^{n+2s}}\,dy+  \int_{B_5\setminus B_1}[g(\cdot,y)]_{\cont^\alpha((-1,0])} \,d\omega_0(y) + \sup_{x\in B_1} \norm{f_4+P_4}_{\cont^\alpha((-1,0])} +\norm{u}_{L^\infty(\overline{Q}_1^*)} \right)\\
				\leq& C \left( J^k_1(u) + \sup_{t_1,t_2\in (-1,0]} \int_{B_1^c}[g(\cdot,y)]_{\cont^\alpha((-1,0])} \,d\omega_k(y) + \sup_{x\in B_1} \norm{f_4+P_4}_{\cont^\alpha((-1,0])} + \norm{u}_{L^\infty(\overline{Q}_1^*)} \right),
			\end{split}
		\end{equation*}
		for some constant $C>0$ depending only on $n$, $s$, and $k$, and possibly changing from line to line. 
		
		In particular, $u\in \cont^0(B_1;\cont^{1+\alpha}((-1,0]))$, as desired.
	\end{proof}
	
	\begin{proof}[Proof of Proposition~\ref{prop::1_zeta_time_reg}]
One observes that the assumptions of Proposition~\ref{prop::schauder_time_reg} are fulfilled and therefore
	Proposition~\ref{prop::1_zeta_time_reg} follows from Proposition~\ref{prop::schauder_time_reg}.
\end{proof}
	
	We now state and prove two partial results required for the proof of Theorem~\ref{th::schauder_modk}. These results provide suitable upper bounds on the spatial $\cont^\gamma$-norm for solutions to a space-fractional heat equation modulo~$k$ with forcing term.
	
	\begin{prop}\label{prop::schauder_cutoff}
		Let $s\in(0,1)$,~$k\in\N$, and~$\xi\geq\sigma$. Let also $u\in\cont^1\left((-1,0];\cont^{\alpha}_{\loc}(B_1)\right)\cap V_k^\xi((-1,0])$, for some $\alpha>2s$, and $f_u$ be as in Corollary~\ref{cor::parabolic_3_2}.
		
		Then, for all $\gamma>k-1$ such that $\gamma\not\in\N$ and $\gamma+2s\not\in\N$, we have that
		\begin{equation} \label{eq::schauder_cutoff}
			\sup_{t\in (-1/2,0]}\norm{u}_{\cont^{\gamma+2s}(B_{1/2})} \leq C\left(J^k_1(u)+\sup_{t\in (-1,0]}[f_u]_{\cont^{\gamma}(B_1)}+\norm{u}_{L^\infty(\overline{Q}^*_1)}\right),
		\end{equation}
		for some constant $C>0$ depending only on $n$, $s$, $k$, and~$\gamma$.
	\end{prop}
	
	\begin{proof}
		Let $\gamma=m+\theta$, with $m\in\N$ and $\theta\in(0,1)$, and observe that $\gamma>k-1$ reads $m\geq k-1$.
		
		Let also $\{\eta_R\}_R$ be a sequence of cut-off functions as in Definition~\ref{def::k_viscosity_sol}. Thus, employing Proposition~\ref{prop::parabolic_3_1} with $\eta=\eta_4$ and setting $v:=\eta_4 u$, we obtain that
		\begin{equation*}
			\partial_t v + (-\Delta)^s v = P_{u,4}+f_{u,4}\quad\mbox{in }Q_1,
		\end{equation*}
		where $P_{u,4}\in \left(\cont^{\sigma}((-1,0])\right)_{k-1}^{ [x] }$.
		
We also recall the notation in~\eqref{eq::def_seminorm_mod_poly} and the observation in~\eqref{eq::def_seminorm_mod_poly22} to infer that
\begin{eqnarray*}&& \sup_{t\in(-1,0]}[P_{u,4}+f_{u,4}]_{\cont^{\gamma}(B_1)}=
\sup_{t\in(-1,0]}[P_{u,4}+f_{u,4}]_{\cont^{\gamma}(B_1;k)}\\&&\quad
= \sup_{t\in(-1,0]} \inf_{P\in\left(\cont^{\sigma}((-1,0])\right)_{k-1}^{ [x] }} [P_{u,4}+f_{u,4}-P]_{\cont^\gamma(B_1)}
\le \sup_{t\in(-1,0]} [f_{u,4}]_{\cont^\gamma(B_1)}.
\end{eqnarray*}
		Therefore, Theorem~\ref{th::ultimate_schauder} entails that
		\begin{equation} \label{eq::u_estimate_schauder_cutoff}
			\begin{split}
				&\sup_{t\in(-1/2,0]} \norm{u(t,\cdot)}_{\cont^{\gamma+2s}(B_{1/2})} = \sup_{t\in(-1/2,0]} \norm{v(t,\cdot)}_{\cont^{\gamma+2s}(B_{1/2})}\\
				&\qquad \leq C\left(J^0_1(v)+\sup_{t\in(-1,0]}[P_{u,4}+f_{u,4}]_{\cont^{\gamma}(B_1)}+\norm{v}_{L^\infty(\overline{Q}^*_1)}\right) \\
				&\qquad \leq C\left(J^0_1(v)+\sup_{t\in(-1,0]}[f_{u,4}]_{\cont^{\gamma}(B_1)}+\norm{u}_{L^\infty(\overline{Q}^*_1)}\right) .
			\end{split}
		\end{equation}
		
		Moreover, recalling the definitions of $f_{u,4}^\star$ and $f_{3,u}$ in~\eqref{eq::def_f1_f2_fstar} and~\eqref{eq::def_f3}, respectively, we find
		\begin{equation}  \label{eq::f_estimate_schauder_cutoff}
			\begin{split}
				&\sup_{t\in(-1,0]} [f_{u,4}]_{\cont^{\gamma}(B_1)} = \sup_{t\in(-1,0]} [f_{u}+f_{u,4}^\star-f_{3,u}]_{\cont^{\gamma}(B_1)}\\
				&\qquad \leq \sup_{t\in(-1,0]} [f_{u}]_{\cont^{\gamma}(B_1)} + \sup_{t\in(-1,0]} [f_{u,4}^\star - f_{3,u}]_{\cont^{\gamma}(B_1)}\\
				&\qquad \leq \sup_{t\in(-1,0]} [f_{u}]_{\cont^{\gamma}(B_1)} + C\int_{B_2^c} \frac{|u(y)|}{|y|^{n+2s+k}}\,dy\\
				&\qquad \leq \sup_{t\in(-1,0]} [f_{u}]_{\cont^{\gamma}(B_1)} + CJ^k_1(u).
			\end{split}
		\end{equation}
		
		Also, we have
		\begin{equation}  \label{eq::tail_estimate_schauder_cutoff}
			J^0_1(v) = \int_{B_5\setminus B_{1/2}} \frac{|\eta_4 u(y)|}{|y|^{n+2s}}\,dy \leq C\int_{B_5\setminus B_{1/2}} \frac{|\eta_4 u(y)|}{|y|^{n+2s+k}}\,dy \leq CJ^k_1(u).
		\end{equation}
		
		Therefore, plugging~\eqref{eq::f_estimate_schauder_cutoff} and~\eqref{eq::tail_estimate_schauder_cutoff} in \eqref{eq::u_estimate_schauder_cutoff}, up to renaming $C$, we obtain
		\begin{equation*}
			\sup_{t\in(-1/2,0]} \norm{u(t,\cdot)}_{\cont^{\gamma+2s}(B_{1/2})} \leq C\left(J^k_1(u)+\sup_{t\in(-1,0]}[f_{u}]_{\cont^{\gamma}(B_1)}+\norm{u}_{L^\infty(\overline{Q}^*_1)}\right) ,
		\end{equation*}
		as desired.
	\end{proof}
	
	\begin{cor} \label{cor::high_schauder_modk}
		Let $s\in(0,1)$,~$k\in\N$, and~$\xi\geq\sigma$. Let also $u:Q_\infty\to\R$ and $f:Q_1\to\R$.
		
		Assume that $u\in\cont^0(Q_1)\cap V_k^\xi((-1,0])$ satisfies
		\begin{equation*}
			\partial_t u + (-\Delta)^s u \eqp{k} f\quad\mbox{in }Q_1.
		\end{equation*}
		
		Then, for all $\gamma>k-1$ such that $\gamma\not\in\N$ and $\gamma+2s\not\in\N$, we have that
		\begin{equation} \label{eq::high_schauder_modk}
			\sup_{t\in (-1/2,0]}\norm{u}_{\cont^{\gamma+2s}(B_{1/2})} \leq C\left(J^k_1(u)+\sup_{t\in (-1,0]}[f]_{\cont^{\gamma}(B_1)}+\norm{u}_{L^\infty(\overline{Q}^*_1)}\right),
		\end{equation}
		for some constant $C>0$ depending only on $n$, $s$, $k$, and~$\gamma$.
	\end{cor}
	
	\begin{proof}
		Let us assume first that 
		\begin{equation} \label{eq::u_extra_reg}
			u\in \cont^1\left((-1,0];\cont^{\alpha}_{\loc}(B_1)\right).
		\end{equation}
			Then, by Corollary~\ref{cor::parabolic_3_2}, we have that
		\begin{equation*}
			\partial_t u + (-\Delta)^s u \eqp{k} f_u\quad\mbox{in }Q_1.
		\end{equation*}
			In particular, $f-f_u\in \left(\cont^{\sigma}((-1,0])\right)_{k-1}^{ [x] }$. Therefore, Proposition~\ref{prop::schauder_cutoff} yields
		\begin{equation} \label{eq::high_schauder_modk_weaker}
			\begin{split}
				\sup_{t\in (-1/2,0]}\norm{u}_{\cont^{\gamma+2s}(B_{1/2})}& \leq C\left(J^k_1(u)+\sup_{t\in (-1,0]}[f_u]_{\cont^{\gamma}(B_1)}+\norm{u}_{L^\infty(\overline{Q}^*_1)}\right) \\
				& \leq C\left(J^k_1(u)+\sup_{t\in (-1,0]}[f]_{\cont^{\gamma}(B_1)}+\norm{u}_{L^\infty(\overline{Q}^*_1)}\right),
			\end{split}
		\end{equation}
		which is precisely~\eqref{eq::high_schauder_modk} under the extra assumption~\eqref{eq::u_extra_reg}.
		
		For the general case, let us consider a sequence of mollifiers $\{\rho_\epsilon\}_\epsilon$, where $\rho_\epsilon(x):=\epsilon^{-n}\rho(x/\epsilon)$, for some~$\rho\in \cont^\infty_c(B_1)$. We define
		\begin{equation*}
			u_\epsilon := u*\rho_\epsilon\quad\mbox{and}\quad f_\epsilon := f*\rho_\epsilon.
		\end{equation*}
		Thus, we have
		\begin{equation*}
			\partial_t u_\epsilon + (-\Delta)^s u_\epsilon = f_\epsilon\quad\mbox{in }Q_{99/100},
		\end{equation*}
		whenever $\epsilon$ is small enough. 
		
		Moreover, recalling also Proposition~\ref{prop::1_zeta_time_reg}, we have that $u_\epsilon\in \cont^1\left((-1,0];\cont^{\alpha}_{\loc}(B_1)\right)$. Then, up to resizing the domain, we infer from Proposition~\ref{prop::schauder_cutoff} that\footnote{
		The estimate in~\eqref{eq::reg_u_epsilon} follows from Proposition~\ref{prop::schauder_cutoff} and a suitable choice of $r>0$ for the term $J^0_{1,r}((u))$ in~\eqref{eq::schauder_seminorm} (see Corollary~\ref{cor::schauder_seminorm}).
		}
		\begin{equation} \label{eq::reg_u_epsilon}
			\sup_{t\in \left(-\left(\frac{9}{10}\right)^{2s},0\right]}\norm{u_\epsilon}_{\cont^{\gamma+2s}(B_{9/10})} \leq C\left(J^k_{\frac{99}{100},\frac{3}{4}}(u_\epsilon) +\sup_{t\in \left(-\left(\frac{99}{100}\right)^{2s},0\right]}[f_\epsilon]_{\cont^{\gamma}(B_{99/100})}+\norm{u_\epsilon}_{L^\infty(\overline{Q}^*_{99/100})}\right).
		\end{equation}
		
		Taking limits as $\epsilon\to0$, we find
		\begin{equation} \label{eq::some_limits}
			\begin{split}
				&\lim_{\epsilon\to0} \sup_{t\in \left(-\left(\frac{9}{10}\right)^{2s},0\right]}\norm{u_\epsilon}_{\cont^{\gamma+2s}(B_{9/10})}  \geq \sup_{t\in (-1/2,0]}\norm{u}_{\cont^{\gamma+2s}(B_{1/2})} ,\\
				&\lim_{\epsilon\to0} \sup_{t\in \left(-\left(\frac{99}{100}\right)^{2s},0\right]}[f_\epsilon]_{\cont^{\gamma}(B_{99/100})} \leq \sup_{t\in (-1,0]}[f]_{\cont^{\gamma}(B_{1})},\\
				\mbox{and}\quad&\lim_{\epsilon\to0} \norm{u_\epsilon}_{L^\infty(\overline{Q}^*_{1})} \leq \norm{u}_{L^\infty(\overline{Q}^*_{1})}.
			\end{split}
		\end{equation}
		
		Moreover, for every $y\in B_{99/200}^c$ and $z\in B_\epsilon(y)$, we have that
		\begin{equation*}
			|z| \leq |y-z|+|y| \leq \epsilon+|y| \leq 2|y|,
		\end{equation*}
		and hence
		\begin{equation} \label{eq::limit_tail}
			\begin{split}
				&J^k_{\frac{99}{100},\frac{3}{4}}(u_\epsilon)  = \sup_{t\in \left(-\left(\frac{99}{100}\right)^{2s},0\right]}\int_{B_{3/4}^c} \frac{|u_\epsilon(t,y)|}{|y|^{n+2s+k}}\,dy \\
				&\qquad = \sup_{t\in\left(-\left(\frac{99}{100}\right)^{2s},0\right]}\int_{B_{3/4}^c} \int_{B_\epsilon(y)} \frac{|\rho_\epsilon(y-z)|\,|u(t,z)|}{|y|^{n+2s+k}}\,dz\,dy\\
				&\qquad \leq C \sup_{t\in\left(-\left(\frac{99}{100}\right)^{2s},0\right]}\int_{B_{3/4}^c} \int_{B_\epsilon(y)} \frac{|\rho_\epsilon(y-z)|\,|u(t,z)|}{|z|^{n+2s+k}}\,dz\,dy\\
				&\qquad \leq C \sup_{t\in\left(-\left(\frac{99}{100}\right)^{2s},0\right]}\int_{B_{3/4}^c} \frac{|u(t,z)|}{|z|^{n+2s+k}}\left(\int_{\R^n}|\rho_\epsilon(y-z)|\,dy\right)\,dz\\
				&\qquad \leq C \sup_{t\in\left(-\left(\frac{99}{100}\right)^{2s},0\right]}\int_{B_{3/4}^c} \frac{|u(t,z)|}{|z|^{n+2s+k}}\,dy\\
				&\qquad \leq J^k_1(u),
			\end{split}
		\end{equation}
		up to relabeling $C$.
		
		Therefore, putting~\eqref{eq::some_limits} and~\eqref{eq::limit_tail} together with~\eqref{eq::reg_u_epsilon}, we deduce~\eqref{eq::high_schauder_modk}, as desired.
	\end{proof}
	
	We are now ready to prove Theorem~\ref{th::schauder_modk}.
	
	\begin{proof}[Proof of Theorem~\ref{th::schauder_modk}]
		Let us first observe that, when $\gamma>k-1$, the claim in~\eqref{eq::schauder} follows immediately from Corollary~\ref{cor::high_schauder_modk},
		the definition
		in~\eqref{eq::def_seminorm_mod_poly} and the observation in~\eqref{eq::def_seminorm_mod_poly22}. Thus, from now on
		we focus on the case~$\gamma\leq k-1$.
		
		To this purpose, we consider a solution $v$ of the space-fractional heat equation with zero-external datum and zero-initial condition, i.e.,
		\begin{equation*}
			\begin{cases}
				\partial_t v + (-\Delta)^s v =f\quad&\mbox{in }Q_1,\\
				v=0\quad&\mbox{in }(-1,0]\times B_1^c,\\
				v(-1,\cdot)=0\quad&\mbox{in }\R^n.
			\end{cases}
		\end{equation*}
		Then, thanks to~\cite[Theorem~5.2]{parabolic_viscosity}, we have
		\begin{equation*}
			\sup_{t\in(-1,0]} \norm{v}_{\cont^\zeta(\R^n)} \leq C\norm{f}_{L^\infty(Q_1)},
		\end{equation*}
		for some $\zeta\in(0,2)$ and some $C>0$. It thereby follows that
		\begin{equation} \label{eq::estimate_v_0_data}
			J^0_1(v) \leq C \sup_{t\in(-1,0]}  \norm{v}_{\cont^\zeta(\R^n)} \leq C\norm{f}_{L^\infty(Q_1)},
		\end{equation}
		and that
	\begin{equation} \label{eq::estimate_v_0_data22} \norm{v}_{L^\infty(\overline{Q}^*_1)}\le \sup_{t\in(-1,0]} \norm{v}_{\cont^\zeta(\R^n)} \le C\norm{f}_{L^\infty(Q_1)},\end{equation}
		up to renaming $C$. 
		
		From these observations and Theorem~\ref{th::ultimate_schauder}, we infer that
		\begin{equation} \label{eq::estimate_v_0_data2}
			\begin{split}
				\sup_{t\in(-1/2,0]} \norm{v}_{\cont^{\gamma+2s}(B_{1/2})} &\leq C\left(J^0_1(v)+\sup_{t\in(-1,0]}[f]_{\cont^{\gamma}(B_1)}+\norm{v}_{L^\infty(\overline{Q}^*_1)}\right) \\
				& \leq C\left(\norm{f}_{L^\infty(Q_1)}+\sup_{t\in(-1,0]}[f]_{\cont^{\gamma}(B_1)}\right).
			\end{split}
		\end{equation}
		
		Moreover, by Corollary~\ref{cor::parabolic_3_8}, we have that
		\begin{equation*}
			\partial_t v + (-\Delta)^s v \eqp{0}f\quad\mbox{in }Q_1,
		\end{equation*}
		and hence
		\begin{equation*}
			\partial_t v + (-\Delta)^s v \eqp{k}f\quad\mbox{in }Q_1.
		\end{equation*}
		
		Thus, the function $w:=u-v$ satisfies
		\begin{equation*}
			\partial_t w + (-\Delta)^s w \eqp{k}0\quad\mbox{in }Q_1.
		\end{equation*}
		
		Let us pick $\epsilon>0$ small, and define $\overline{\gamma}:=k-1+\epsilon>k-1$. Therefore, Corollary~\ref{cor::high_schauder_modk} entails that
		\begin{equation} \label{eq::estimate_w_above}\begin{split}
			\sup_{t\in(-1/2,0]} \norm{w}_{\cont^{\overline{\gamma}+2s}(B_{1/2})} &\leq C\left(J^k_1(w)
			+\norm{w}_{L^\infty(\overline{Q}^*_1)}\right)
			\\& \leq C\left(J^k_1(u)+J^k_1(v)
			 +\norm{u}_{L^\infty(\overline{Q}^*_1)}
			 +\norm{v}_{L^\infty(\overline{Q}^*_1)}
			 \right) \\& \leq C\left(J^k_1(u) +\norm{u}_{L^\infty(\overline{Q}^*_1)}+
			 \norm{f}_{L^\infty(Q_1)}\right),\end{split}
		\end{equation}
		where in the last inequality we have used~\eqref{eq::estimate_v_0_data} and~\eqref{eq::estimate_v_0_data22} .
	
	On the other hand, $\overline{\gamma}>\gamma$ yields
	\begin{equation} \label{eq::estimate_w_below}
		\sup_{t\in(-1/2,0]} \norm{w}_{\cont^{\overline{\gamma}+2s}(B_{1/2})} \geq \sup_{t\in(-1/2,0]} \norm{u}_{\cont^{ {\gamma}+2s}(B_{1/2})} -\sup_{t\in(-1/2,0]} \norm{v}_{\cont^{ {\gamma}+2s}(B_{1/2})}.
	\end{equation}
	
	From~\eqref{eq::estimate_v_0_data2},~\eqref{eq::estimate_w_above}, and~\eqref{eq::estimate_w_below}, we deduce that
	\begin{equation*} 
		\begin{split}
			\sup_{t\in(-1/2,0]} \norm{u}_{\cont^{\gamma+2s}(B_{1/2})} &\leq C\left(\sup_{t\in(-1/2,0]} \norm{v}_{\cont^{\gamma+2s}(B_{1/2})}+
			 \norm{f}_{L^\infty(Q_1)}+J^k_1(u) +\norm{u}_{L^\infty(\overline{Q}^*_1)}\right) \\
			& \leq C\left(\norm{f}_{L^\infty(Q_1)}+\sup_{t\in(-1,0]}[f]_{\cont^{\gamma}(B_1)}+J^k_1(u)+\norm{u}_{L^\infty(\overline{Q}^*_1)}\right),
		\end{split}
	\end{equation*}
	up to relabeling $C$.
	
	Now, notice that, in light of Remark~\ref{rem::mod_poly}, we also have
	\begin{equation} \label{eq::u_semifinal_estimate}
		\sup_{t\in(-1/2,0]} \norm{u}_{\cont^{\gamma+2s}(B_{1/2})} \leq C\left(\norm{f+P}_{L^\infty(Q_1)}+\sup_{t\in(-1,0]}[f+P]_{\cont^{\gamma}(B_1)}+J^k_1(u)+\norm{u}_{L^\infty(\overline{Q}^*_1)}\right),
	\end{equation}
	for any $P\in\left(\cont^{\sigma}((-1,0])\right)_{k-1}^{ [x] }$.
	
	From~\eqref{eq::u_semifinal_estimate} and Lemma~\ref{lemma::C1}, we conclude that
	\begin{equation*}
		\sup_{t\in(-1/2,0]}\left(\norm{u(t,\cdot)}_{\cont^{\gamma+2s}(B_{1/2})}\right)	\leq C\left(\sup_{t\in(-1,0]}[f]_{\cont^{\gamma}(B_1;k)}+J^k_1(u) + \norm{u}_{L^\infty(\overline{Q}^*_1)}\right),
	\end{equation*}
	also in the case $\gamma\leq k-1$.
	\end{proof}
	

	\section{Proof of Theorem~\ref{th::stability} \\ Stability of solutions to a space-fractional heat equation modulo~$k$}
	\label{sec::proof_stability}
	
	This section is devoted to the proof of the stability result for solutions to a space-fractional heat equation modulo~$k$ in Theorem~\ref{th::stability}. In the wake of~\cite{elliptic_poly}, we first address the case in which a solution $u$ to~\eqref{eq::viscosity_heat_modk} vanishes in~$Q_1$.
	
	\begin{lemma} \label{lemma::stability_vanishing}
		Let $s\in(0,1)$, $k\in\N$, and~$\xi\geq\sigma$. Let us consider sequences~$\{u_j:Q_\infty\to\R\}_j$ and $\{g_j:Q_1\to\R\}_j$ such that, for all $j$, $u_j\in\cont^0\left(Q_1\right)\cap V_k^\xi((-1,0])$ is a viscosity solution modulo~$k$ of~\eqref{eq::viscosity_heat_modk} with forcing term $g_j$ in $Q_1$, and $u_j=0$ in $Q_1$. Moreover, suppose that the sequence $\{u_j\}_j$ has uniformly bounded $L^1(d\omega_k)$-norm. Namely, 
		
		\begin{eqnarray}&&
			\label{eq::heat_eq_modk_system_j}
				u_j=0\qquad{\mbox{and}}\qquad
				\partial_t u_j + (-\Delta)^s u_j \eqp{k} g_j\qquad\mbox{in $Q_1$, for all $j$,}
			\\	
			\label{eq::unif_energy_bound2}\mbox{and}\quad&&\sup_{j\in\N}\, \sup_{t\in (-1,0]} \int_{\R^n} |u_j(t,y)|\, d\omega_k (y)<+\infty.
		\end{eqnarray}
		
		Suppose also that there exist functions $u\in\cont^0((-1,0];L^1(\R^n))$ and $g:Q_1\to\R$ such that $u=0$ in $Q_1$ and, as $j\to+\infty$,
		\begin{align}
			\label{eq::ptw_conv_g} &g_j\to g\quad\mbox{a.e. in }Q_1\\
			\label{eq::weak_conv_u2} \mbox{and}\quad&\int_{B^c_1}\frac{ u_j(t,y)\,\varphi(y)}{|y|^{n+2s+k}}\, dy \to \int_{B^c_1}\frac{ u(t,y)\,\varphi(y)}{|y|^{n+2s+k}}\, dy,\quad\mbox{ for all $t\in(-1,0]$ and $\varphi \in L^\infty(B_1^c)$}.
		\end{align}
		
		Then, 
		\begin{equation*} 
			\partial_t u + (-\Delta)^s u \eqp{k} g \quad\mbox{in }Q_1. 
		\end{equation*} 
	\end{lemma}
	
	\begin{proof}
		Since $u_j=0$ in $Q_1$, we are in the setting of Corollary~\ref{cor::parabolic_3_2}. Therefore, for every $j$ and every $(t,x)\in Q_1$, we obtain
		\begin{equation*}
			\partial_t u_j(t,x) + (-\Delta)^s  u_j(t,x) \eqp{k} f_{u_j}(t,x) = -\int_{B_{2}\setminus B_1} \frac{u_j(t,y)}{|x-y|^{n+2s}}\,dy + \int_{B_{2}^c} \frac{u_j(t,y) \psi(x,y)}{|y|^{n+2s+k}}\,dy.
		\end{equation*}
		Thus, we infer that 
		\begin{equation*}
			P_j := g_j - f_{u_j} \in \left(\cont^{\sigma}((-1,0])\right)_{k-1}^{ [x] }.
		\end{equation*}
		
		Moreover, since the function~$y\mapsto |x-y|^{-n-2s}$ belongs to~$L^\infty(B_1^c)$ for any fixed~$x\in B_1$, by~\eqref{eq::ptw_conv_g} and~\eqref{eq::weak_conv_u2}, we have that
		\begin{equation} \label{eq::limit_P}
			\lim_{j\to+\infty} P_j(t,x) = g(t,x) + \int_{B_{2}\setminus B_1} \frac{u(t,y)}{|x-y|^{n+2s}}\,dy - \int_{B_{2}^c} \frac{u(t,y) \psi(x,y)}{|y|^{n+2s+k}}\,dy =: P(t,x),
		\end{equation}
		for every given $(t,x)\in Q_1$. In particular, $P\in \left(\cont^{\sigma}((-1,0])\right)_{k-1}^{ [x] }$ (see for instance~\cite[Lemma~2.1]{elliptic_poly}).
		
We also claim that		\begin{equation}\label{cjwtr879ywfkjsa98764kjjdsak45678ger}
\sup_{t\in(-1,0]}\int_{B_1^c} \frac{|u(t,y)|}{|y|^{n+2s+k}}\,dy<+\infty.\end{equation}
		Indeed, by~\eqref{eq::weak_conv_u2}, for all~$t\in(-1,0]$,
\begin{eqnarray*}
&&\int_{B_1^c} \frac{|u(t,y)|}{|y|^{n+2s+k}}\,dy=
\int_{B_1^c} \frac{u(t,y)\chi_{\{u>0\}}}{|y|^{n+2s+k}}\,dy-\int_{B_1^c} \frac{u(t,y)\chi_{\{u\leq0\}}}{|y|^{n+2s+k}}\,dy
\\&&\quad=\lim_{j\to+\infty}
\int_{B_1^c} \frac{u_j(t,y)\chi_{\{u>0\}}}{|y|^{n+2s+k}}\,dy-\int_{B_1^c} \frac{u_j(t,y)\chi_{\{u\leq0\}}}{|y|^{n+2s+k}}\,dy\\
&&\quad\le 2\sup_{j\in\N}\int_{B_1^c} \frac{|u_j(t,y)|}{|y|^{n+2s+k}}\,dy\le 2\sup_{t\in(-1,0]}\,\sup_{j\in\N}\int_{B_1^c} \frac{|u_j(t,y)|}{|y|^{n+2s+k}}\,dy\\
&&\quad\le  2\sup_{j\in\N}\,\sup_{t\in(-1,0]}\int_{\R^n} |u_j(t,y)|\,d\omega_k(y).
\end{eqnarray*}		
Thus the claim in~\eqref{cjwtr879ywfkjsa98764kjjdsak45678ger} follows from~\eqref{eq::unif_energy_bound2}.
		
		Therefore, we can apply Corollary~\ref{cor::parabolic_3_2} to $u$ and obtain that
		\begin{equation*}
			\partial_t u(t,x) + (-\Delta)^s  u(t,x) \eqp{k} f_u(t,x)=-\int_{B_{2}\setminus B_1} \frac{u(t,y)}{|x-y|^{n+2s}}\,dy + \int_{B_{2}^c} \frac{u(t,y) \psi(x,y)}{|y|^{n+2s+k}}\,dy.
		\end{equation*}
		This, together with~\eqref{eq::limit_P}, reads
		\begin{equation*}
			\partial_t u + (-\Delta)^s u \eqp{k} g \quad\mbox{in }Q_1.
		\end{equation*}
		as desired.
	\end{proof}
	
	\begin{proof}[Proof of Theorem~\ref{th::stability}]
		Let 
		\begin{equation*}
			\eta := 
			\begin{cases}
				1\quad&\mbox{in }B_1,\\
				0\quad&\mbox{in }B_1^c
			\end{cases}
		\end{equation*}
	and define
		\begin{equation*}
			\begin{split}
				&v_j := (1-\eta)u_j\\
				&w_j :=\eta u_j
			\end{split}
			\qquad\mbox{and}\qquad
			\begin{split}
				&v := (1-\eta)u\\
				&w :=\eta u.
			\end{split}
		\end{equation*}
Recalling~\eqref{eq::local_unif_conv_u}, we have that $v_j\to v$ and $w_j\to w$ locally uniformly in $Q_1$, as $j\to+\infty$.
		
		Now, notice that, by~\eqref{eq::heat_eq_modk_j} and Corollary~\ref{cor::parabolic_3_10} (used here with~$\rho:=1$), there exists a polynomial $P_j\in\left(\cont^{\sigma}((-1,0])\right)_{k-1}^{ [x] }$ such that
		\begin{equation*}
			\begin{split}
				&\partial_t w_j(t,x) + (-\Delta)^s w_j(t,x) \\
				&\qquad\qquad= f_j(t,x)+ P_j(t,x) +\int_{B_{2}\setminus B_1} \frac{u_j(t,y)}{|x-y|^{n+2s}}\,dy - \int_{B_{2}^c} \frac{u_j(t,y) \psi(x,y)}{|y|^{n+2s+k}}\,dy\\
				&\qquad\qquad  =: h_j(t,x).
			\end{split}
		\end{equation*}
Thus, from Corollary~\ref{cor::parabolic_3_8} we find
$$ \partial_t w_j(t,x) + (-\Delta)^s w_j(t,x)\eqp{0}h_j(t,x)$$
		and therefore
	\begin{equation} \label{eq::heat_eq_modk_w_j}\partial_t w_j(t,x) + (-\Delta)^s w_j(t,x)\eqp{k}h_j(t,x).		\end{equation}
		
		It thereby follows that
		\begin{equation}  \label{eq::heat_eq_modk_v_j}
			\begin{split}
				&\partial_t v_j(t,x) + (-\Delta)^s v_j(t,x) = \partial_t u_j(t,x) + (-\Delta)^s u_j(t,x) - \partial_t w_j(t,x) - (-\Delta)^s w_j(t,x)\\
				&\qquad \eqp{k} f_j(t,x) - h_j(t,x)\\
					&\qquad \eqp{k} P_j(t,x) +\int_{B_{2}\setminus B_1} \frac{u_j(t,y)}{|x-y|^{n+2s}}\,dy - \int_{B_{2}^c} \frac{u_j(t,y) \psi(x,y)}{|y|^{n+2s+k}}\,dy
					\\
				&\qquad \eqp{k} \int_{B_{2}\setminus B_1} \frac{u_j(t,y)}{|x-y|^{n+2s}}\,dy - \int_{B_{2}^c} \frac{u_j(t,y) \psi(x,y)}{|y|^{n+2s+k}}\,dy \\&\qquad=:g_j(t,x).
			\end{split}
		\end{equation}
		
		Observe that, as~$j\to+\infty$, 
		\begin{equation} \label{eq::g_conv_unif}
			g_j \to g:= \int_{B_{2}\setminus B_1} \frac{u(t,y)}{|x-y|^{n+2s}}\,dy - \int_{B_{2}^c} \frac{u(t,y) \psi(x,y)}{|y|^{n+2s+k}}\,dy\quad\mbox{locally uniformly in $Q_1$.}
		\end{equation}
		Indeed, thanks to~\eqref{eq::local_unif_conv_u}  and~\eqref{eq::estimate_psi}, $g_j\to g$ pointwise in $Q_1$. In addition, for any $r\in(0,1)$, we have
		\begin{equation*}
			\begin{split}
				&\sup_{(t,x)\in Q_{r}} |D_x g_j| \leq C\sup_{(t,x)\in Q_{r}} \left(\int_{B_{2}\setminus B_1} \frac{|u_j(t,y)|}{|x-y|^{n+2s+1}}\,dy + \int_{B_{2}^c} \frac{|u_j(t,y)| \,|D_x\psi(x,y)|}{|y|^{n+2s+k}}\,dy \right)\\
				&\qquad\leq C \sup_{t\in(-r^{2s},0]} \left(\int_{B_{2}\setminus B_1} \frac{|u_j(t,y)|}{(1-r)^{n+2s+1}}\,dy + \int_{B_2^c} |u_j(t,y)|\,d\omega_k(y)\right) \\
				&\qquad\leq C\sup_{j\in\N}\, \sup_{t\in (-1,0]} \int_{\R^n} |u_j(t,y)|\, d\omega_k (y) \leq C, 
			\end{split}
		\end{equation*}
		for some positive constant $C$, possibly changing from line to line. Here, for the last inequality we have used~\eqref{eq::unif_energy_bound}. 
With this uniform bound one obtains~\eqref{eq::g_conv_unif}.

		
		Thus, we are in the position of exploiting Lemma~\ref{lemma::stability_vanishing} for~$v$ and conclude that
		\begin{equation} \label{eq::limit_eq_v}
			\partial_t v + (-\Delta)^s v \eqp{k} g\quad\mbox{in }Q_1.
		\end{equation}
		
		Now, let us consider $\varphi\in\cont^\infty_c(Q_1)$, and let $U:=\mbox{supp}(\varphi) \Subset Q_1$. Let us also pick $\epsilon>0$ and consider the function $\rho_\epsilon(x):= \epsilon^{-n}\rho(x/\epsilon)$, for some $\rho\in\cont^\infty_c(B_1)$. Then, we set
		\begin{equation*}
			w_{j,\epsilon} := w_j*\rho_\epsilon \qquad\mbox{and}\qquad h_{j,\epsilon} := h_j*\rho_\epsilon.
		\end{equation*}
		
		Thus, the spatial support of~$w_{j,\epsilon}$ is compactly contained in~$ B_{3/2}$, and
		\begin{equation} \label{eq::heat_eq_modk_w_j_eps}
			\partial_t w_{j,\epsilon} + (-\Delta)^s w_{j,\epsilon} = h_{j,\epsilon} \quad\mbox{in }U,
		\end{equation}
		in the classical sense (see e.g.~\cite{MR3161511} for the elliptic case). Therefore, integrating by parts,
		\begin{equation*}
			\begin{split}
				&\int_{U} h_{j,\epsilon}\,\varphi = \int_{U} (\partial_t w_{j,\epsilon}+(-\Delta)^s w_{j,\epsilon})\,\varphi = \int_{(-1,0]\times\R^n} (\partial_t w_{j,\epsilon}+(-\Delta)^s w_{j,\epsilon})\,\varphi\\
				&\qquad=  \int_{(-1,0]\times\R^n} w_{j,\epsilon}\, (\partial_t\varphi+(-\Delta)^s \varphi) =  \int_{(-1,0]\times B_{3/2}} w_{j,\epsilon}\, (\partial_t\varphi+(-\Delta)^s \varphi).
			\end{split}
		\end{equation*}
		From this, for every $j,\ell\in\N$, we obtain
		\begin{equation} \label{eq::h_j_cauchy}
			\begin{split}
				&\left| \int_{U} (h_{j}-h_{\ell}) \,\varphi\right| = \lim_{\epsilon\to0} \left| \int_{U} (h_{j,\epsilon}-h_{\ell,\epsilon}) \,\varphi\right| \leq \lim_{\epsilon\to0}  \int_{(-1,0]\times B_{3/2}} |w_{j,\epsilon}-w_{\ell,\epsilon}| \,|\partial_t\varphi+(-\Delta)^s \varphi| \\
				&\qquad  = \int_{(-1,0]\times B_{3/2}} |w_{j}-w_{\ell}| \,|\partial_t\varphi+(-\Delta)^s \varphi| \leq C\norm{w_{j}-w_{\ell}}_{L^1((-1,0]\times B_{3/2})} \norm{\varphi}_{\cont^2(Q_\infty)}.
			\end{split}
		\end{equation}
		
		Moreover,
		\begin{equation} \label{eq::u_L1_conv}
			\lim_{j\to+\infty} \norm{w_j-\eta u}_{L^1(Q_1)} = \lim_{j\to+\infty} \norm{u_j-u}_{L^1(Q_1)} =0.
		\end{equation}
		Therefore, from~\eqref{eq::h_j_cauchy} and~\eqref{eq::u_L1_conv}, we deduce that $\{h_j\}_j$ is a Cauchy sequence with respect to the norm~$\norm{\cdot}_{\star}$, defined for any function $\phi$ as
		\begin{equation*}
			\norm{\phi}_{\star} := \sup_{\substack{\varphi\in\cont^2_c(Q_1) \\ \norm{\varphi}_{\cont^2(Q_\infty)}\le1} } \int_{Q_1}  \phi\,\varphi\, .
		\end{equation*}
		Thanks to~\eqref{eq::f_conv_unif}, we find that also $\{f_j\}_j$ is a Cauchy sequence in the norm $\norm{\cdot}_\star$, and, by~\eqref{eq::g_conv_unif} so is $\{g_j\}_j$.
		It thereby follows that
		\begin{equation} \label{eq::P_j_cauchy}
			P_j = g_j+h_j-f_j
		\end{equation}
		is a Cauchy sequence in the norm $\norm{\cdot}_\star$, and hence there exists $P\in\left(\cont^{\sigma}((-1,0])\right)_{k-1}^{ [x] }$ such that
		\begin{equation} \label{eq::P_conv_unif}
			\text{$P_j\to P$ uniformly in $U$, for any $U\Subset Q_1$.}
		\end{equation} 
		
		Thus, we infer from~\eqref{eq::f_conv_unif} and~\eqref{eq::P_conv_unif} that also $\{f_j+P_j\}_j$ converges locally uniformly in~$Q_1$, and so does~$\{h_j\}_j$. Let~$h$ be such a limit. Then, by~\cite[Theorem~5.6]{parabolic_viscosity} and taking the limit as~$j\to+\infty$ in~\eqref{eq::P_j_cauchy}, we deduce that
		\begin{equation*}
			\partial_t w + (-\Delta)^s w = h = f-g+P\quad\mbox{in }Q_1,
		\end{equation*}
		in the viscosity sense.
		
		Now, recalling Corollary~\ref{cor::parabolic_3_8}, we obtain
		\begin{equation*} 
			\partial_t w + (-\Delta)^s w \eqp{0} f-g+P\quad\mbox{in }Q_1
		\end{equation*}
		and hence
		\begin{equation} \label{eq::limit_eq_w}
			\partial_t w + (-\Delta)^s w \eqp{k} f-g\quad\mbox{in }Q_1.
		\end{equation}
		
		Therefore, from~\eqref{eq::limit_eq_v} and~\eqref{eq::limit_eq_w}, we infer that, in~$Q_1$,
		\begin{equation*}
				\partial_t u + (-\Delta)^s u = \partial_t v + (-\Delta)^s v +\partial_t w + (-\Delta)^s w 
				\eqp{k} g + f-g = f,
			\end{equation*}
		concluding the proof.
	\end{proof}
	
	
	\section{Proof of Theorem~\ref{th::liouville} \\ Liouville theorem for $s$-caloric functions modulo~$k$}
	\label{sec::proof_liouville}
	Here, we prove that any function satisfying the growth condition~\eqref{eq::ancient_sol_0} is $s$-caloric modulo~$k$ if and only if it
	is a polynomial in space and time, as stated in Theorem~\ref{th::liouville}. We see Theorem~\ref{th::liouville} as the analog of the Liouville theorem for~$s$-caloric functions with polynomial growth.
	
	Theorem~\ref{th::liouville} follows immediately by employing Proposition~\ref{prop::liouville_weak} below.
	
	\begin{prop} \label{prop::liouville_weak}
		Let $s\in(0,1)$ and $k\in\N$.		
		Let also~$\xi\geq\sigma$ and $u\in\cont^0\left(Q_\infty\right) \cap V_k^\xi((-\infty,0])$.
		Suppose that~$u$ satisfies the growth condition
		\begin{equation} \label{eq::ancient_sol}
			\norm{u}_{L^\infty(Q_R)} \leq CR^{d(k,s)},\quad\mbox{for all }R>0.
		\end{equation}
		
Suppose also that~$u$ is $s$-caloric strongly modulo~$k$, i.e.
		\begin{equation*} 
			\partial_t u + (-\Delta)^s u \eqp{k} 0 \quad\mbox{in }Q_\infty.
		\end{equation*}
	
	Then, $u \in\left(\cont^{1+\sigma}((-\infty,0])\right)_{d(k,s)}^{ [x] }$. 
	\end{prop}
	
	\begin{proof} 
		Let $u$ satisfy the assumptions of Proposition~\ref{prop::liouville_weak} and suppose that 
		\begin{equation} \label{eq::liouville_modk}
			\partial_t u +(-\Delta)^s u \eqp{k} 0 \quad\mbox{in }Q_\infty.
		\end{equation}
		First, we show that
		\begin{equation} \label{eq::poly_worse_regularity}
			u \in \left(\cont^{\sigma}((-\infty,0])\right)_{d(k,s)}^{ [x] }.
		\end{equation}
		
		To this aim, let us take $m:=d(k,s)+1$.
		Recall that, from Proposition~\ref{prop::1_zeta_time_reg}, we know that $u$ is continuously differentiable in time. Besides, thanks to Theorem~\ref{th::schauder_modk} and recalling~\eqref{eq::ancient_sol}, for any $j\in\N$ and $t\in(-\infty,0]$, we have that 
		\begin{equation}\label{vhwerui34t972830r8gcuk0}\begin{split}
		\norm{D_x^m u(t,\cdot)}_{L^\infty(B_{j/2})} &\leq  C\left(j^{k+2s-m}J^k_j(u)+j^{-m}\norm{u}_{L^\infty(\overline{Q}^*_j)}\right)
		\\&\leq C\left(j^{k+2s-m}J^k_j(u)+j^{-m}\norm{u}_{L^\infty(\overline{Q}_j)}\right)\\&\leq
 C\left(j^{k+2s-m}J^k_j(u)+j^{d(k,s)-m}\right),
	\end{split}	\end{equation}
	up to renaming~$C>0$.
	
Now, we observe that, using polar coordinates
and recalling~\eqref{eq::ancient_sol}, for all~$t\in(-j^{2s},0]$,
\begin{equation}\label{vhwerui34t972830r8gcuk}\begin{split}&
\int_{B_{j/2}^c} \frac{|u(t,y)|}{|y|^{n+2s+k}}\,dy
=\int_{j/2}^{+\infty}\left(\int_{\partial B_{r}} |u(t,r\omega)|\,d{\mathcal{H}}^{n-1}_\omega\right)\frac{dr}{r^{n+2s+k}}\\
&\quad\le \int_{j/2}^{+\infty}\left(\int_{\partial B_{r}}
\sup_{\substack{\tau\in(-j^{2s},0] \\ z\in  B_r}}
|u(\tau,z)|\,d{\mathcal{H}}^{n-1}_\omega\right)\frac{dr}{r^{n+2s+k}} \\
&\quad\le C
\int_{j/2}^{+\infty}
\sup_{\substack{\tau\in(-(2r)^{2s},0] \\ z\in  B_{2r}}}
|u(\tau,z)| \,\frac{dr}{r^{1+2s+k}} \\
&\quad\le C \int_{j/2}^{+\infty}\frac{dr}{r^{1+2s+k-d(k,s)}}.
\end{split}\end{equation}
We also recall from~\eqref{eq::def_d_ks} that if~$s\in\left(0,\frac{1}{2}\right]$, then~$ 1+2s+k-d(k,s)=1+2s>1$, and if~$s\in\left(\frac{1}{2},1\right)$,
then~$1+2s+k-d(k,s)=2s>1$.
These observations and~\eqref{vhwerui34t972830r8gcuk} give that
$$ 
J^k_j(u)= \sup_{t\in(-j^{2s},0]} \int_{B_{j/2}^c} \frac{|u(t,y)|}{|y|^{n+2s+k}}\,dy\le
C j^{d(k,s)-2s-k}.
$$

Plugging this information into~\eqref{vhwerui34t972830r8gcuk0}
we find that
$$ \norm{D_x^m u(t,\cdot)}_{L^\infty(B_{j/2})} \leq  Cj^{d(k,s)-m}.$$
	Taking the limit as $j\to+\infty$, we find that $D_x^m u(t,\cdot)\equiv0$. Recalling also that $m-1=d(k,s)$, we deduce that $u\in(\cont^{\sigma}((-\infty,0]))_{d(k,s)}^{ [x] }$, showing~\eqref{eq::poly_worse_regularity}.
	
	To conclude, we employ Proposition~\ref{prop::1_zeta_time_reg} to infer that $u$ is of class $\cont^{1+\sigma}$ in time, which yields $u\in(\cont^{1+\sigma}((-\infty,0]))_{d(k,s)}^{ [x] }$, as desired.
	\end{proof}
	
	
	\section{Proof of Proposition~\ref{prop::existence} \\ Existence of viscosity solutions modulo~$k$} \label{sec::existence}
	
	In this section, we establish the existence of viscosity solutions (modulo~$k$) to the initial-value problem associated with~\eqref{eq::viscosity_heat_modk} with Dirichlet external condition. The idea of the proof relies on a cut-off argument (see~\cite[Theorem 1.4]{MR4297378} for the elliptic case).
	
	\begin{proof}[Proof of Proposition~\ref{prop::existence}]
		Let $\overline{u}\in V_k^\xi((-1,0])$,~$f\in \cont^0(B_1;\cont^\alpha((-1,0]))$, and~$v\in L^\infty(B_1)\cap L^1(d\omega_k)$, as in the statement of Proposition~\ref{prop::existence}. Then, we define
		\begin{align*}
			&u_1(x,t) := \chi_{B_2}(x)\overline{u}(x,t),\\
			&u_2(x,t) := \chi_{B_2\setminus B_1}(x)\overline{u}(x,t),\\
			\text{and}\qquad &\widetilde{v}(x) := \chi_{B_2}(x)v(x).
		\end{align*}
		
		In particular, $u_1$ is smooth in $Q_1$, and, by Remark~\ref{rem::parabolic_3_5}, satisfies
		\begin{equation} \label{eq::eq_u_1}
			\partial_t u_1 + (-\Delta)^s u_1  \eqp{k} f_{u_1} \quad \mbox{in }Q_1.
		\end{equation}
		
		Let $\widetilde{f}:=f-f_{u_1}$. Then, we consider the Dirichlet problem
		\begin{equation*}
			\begin{cases}
				\partial_t \widetilde{u} + (-\Delta)^s \widetilde{u}  = \widetilde{f} \quad &\mbox{in }Q_1,\\
				\widetilde{u}=u_2\quad &\mbox{in }B_1^c\times(-1,0],\\
				\widetilde{u}(\cdot,0) = \widetilde{v}\quad &\mbox{in }\R^n.
			\end{cases}
		\end{equation*}
		
		Thanks to Corollary~\ref{cor::parabolic_3_8}, we infer that $\widetilde{u}$ satisfies 
		\begin{equation*}
			\begin{cases}
				\partial_t \widetilde{u} + (-\Delta)^s \widetilde{u}  \eqp{k} \widetilde{f} \quad &\mbox{in }Q_1,\\
				\widetilde{u}=u_2\quad &\mbox{in }B_1^c\times(-1,0],\\
				\widetilde{u}(\cdot,0) = \widetilde{v}\quad &\mbox{in }\R^n.
			\end{cases}
		\end{equation*}
		
		Therefore, the function $u:=u_1+\widetilde{u}$ is a viscosity solution (mod~$k$) of~\eqref{eq::dirichlet_modk}, as desired.
	\end{proof}
	
	
	\section{Proof of Theorem~\ref{th::liouville_strong} \\ Liouville theorem for $s$-caloric functions strongly modulo~$k$}
	\label{sec::proof_liouville_strong}
	
	This section is devoted to proving a Liouville theorem for $s$-caloric functions strongly modulo~$k$. Theorem~\ref{th::liouville_strong} follows from Proposition~\ref{prop::liouville_weak} and the forthcoming partial result.
	
	\begin{lemma} \label{lemma::partial_liouville}
		Let $s\in(0,1)$ and $k\in\N$.		
		Let also~$\xi\in[\sigma,k]$ and $u\in\cont^0\left(Q_\infty\right) \cap V_k^\xi((-\infty,0])$.
		Suppose that $u$ satisfies the growth condition
		\begin{equation*} 
			\norm{u}_{L^\infty(Q_R)} \leq CR^{d(k,s)},\quad\mbox{for all }R>0,
		\end{equation*}
		and $u$ is $s$-caloric strongly modulo~$k$, i.e.
		\begin{equation*}
			\partial_t u + (-\Delta)^s u \eqps{k} 0 \quad\mbox{in }Q_\infty.
		\end{equation*}
		
		Then, $u=P+H$, for some polynomials $P\in\R_k^{ [t,x] }$ and $H\in\R_{d(k,s)}^{ [x] }$. 
	\end{lemma}
		
	\begin{proof}
		Let us suppose that $u$ is $s$-caloric strongly modulo~$k$. In particular, $u$ is $s$-caloric modulo~$k$, i.e.
		\begin{equation*}
			\partial_t u + (-\Delta)^s u \eqp{k} 0\quad\mbox{in }Q_\infty.
		\end{equation*}
		Thus, by Proposition~\ref{prop::liouville_weak}, we know that $u\in\left(\cont^{1+\sigma}((-\infty,0])\right)_{d(k,s)}^{ [x] }$,
		namely
		\begin{equation}\label{bcasr4365rfejkwguhwfo6789l98765}
		u(t,x)=\sum_{\substack{\alpha\in\N^n \\ |\alpha|\le d(k,s)}} c_\alpha(t)\,x^\alpha,\end{equation}
		with~$c_\alpha\in \cont^{1+\sigma}((-\infty,0])$.
		
For any~$R>0$, let also $\eta_R$ be a smooth cut-off function such that
		\begin{equation*}
			\eta_R =
			\begin{cases}
				1\quad&\mbox{in }B_R,\\
				0\quad&\mbox{in }B_{R+1}.
			\end{cases}
		\end{equation*}
		Then, by~\cite[Theorem~1.5]{elliptic_poly},
		 for every $\alpha\in\N^n$, we have that
		\begin{equation*}
			(-\Delta)^s(\eta_Rx^\alpha) = \sum_{\substack{\beta\in\N^n \\ |\beta| \leq |\alpha| -1}} d^\alpha_\beta x^\beta + \rho^\alpha_R(x),
		\end{equation*}
		with $d^\alpha_\beta\in\R$, for all $\beta$, and $\rho^\alpha_R\to0$ as $R\to+\infty$.
		
		As a consequence,
	\begin{equation*}
			\begin{split}
				&\partial_t(\eta_R u)+(-\Delta)^s (\eta_R u) = \left(\partial_t+(-\Delta)^s\right)
				\sum_{\substack{\alpha\in\N^n \\ |\alpha|\le d(k,s)}} \eta_R(x) c_\alpha(t) x^\alpha \\
				&\quad=\sum_{\substack{\alpha\in\N^n \\ |\alpha|\le d(k,s)}} \partial_tc_\alpha(t)\eta_R(x)  x^\alpha 
				+ \sum_{\substack{\alpha\in\N^n \\ |\alpha|\le d(k,s)}} c_\alpha(t) (-\Delta)^s(\eta_R(x) x^\alpha) \\
				&\quad=\sum_{\substack{\alpha\in\N^n \\ |\alpha|\le d(k,s)}} \partial_tc_\alpha(t) \eta_R(x) x^\alpha 
				+ \sum_{\substack{\alpha\in\N^n \\ |\alpha|\le d(k,s)}} c_\alpha(t) \left(
				\sum_{\substack{\beta\in\N^n \\ |\beta| \leq |\alpha| -1}} d^\alpha_\beta x^\beta + \rho^\alpha_R(x)\right)\\
				&\quad=\sum_{\substack{\alpha\in\N^n \\ |\alpha|\le d(k,s)}} \partial_tc_\alpha(t)  \eta_R(x)x^\alpha 
				+ \sum_{\substack{\alpha\in\N^n \\ |\alpha|\leq d(k,s)-1}}\left( \sum_{\substack{\beta\in\N^n \\ |\alpha| +1 \leq |\beta| \leq d(k,s)}} c_\beta(t) 
				d_\alpha^\beta x^\alpha \right)+
				\sum_{\substack{\alpha\in\N^n \\ |\alpha|\le d(k,s)}}  \rho^\alpha_R(x)\\
				&\quad
				=\sum_{\substack{\alpha\in\N^n \\ |\alpha| = d(k,s)}} \partial_t c_\alpha(t) \eta_R(x)x^\alpha + \sum_{\substack{\alpha\in\N^n \\ |\alpha|\leq d(k,s)-1}} \left(\partial_tc_\alpha(t)\eta_R(x)+\sum_{\substack{\beta\in\N^n \\ |\alpha| +1 \leq |\beta| \leq d(k,s)}}c_\beta(t)d^\beta_\alpha\right) x^\alpha + \sum_{\substack{\alpha\in\N^n \\ |\alpha|\leq d(k,s)}} \rho^\alpha_R(x).
			\end{split}
		\end{equation*}	
	
	Let now~$M\in(0,R)$. 	
Since, by assumption, $ \partial_t u+(-\Delta)^s u \eqps{k}0$ in $Q_M$, we necessarily have that
		\begin{equation*}
			\sum_{\substack{\alpha\in\N^n \\ |\alpha| = d(k,s)}} \partial_t c_\alpha(t) x^\alpha + \sum_{\substack{\alpha\in\N^n \\ |\alpha|\leq d(k,s)-1}} \left(\partial_tc_\alpha(t) +\sum_{\substack{\beta\in\N^n \\ |\alpha| +1 \leq |\beta| \leq d(k,s)}}c_\beta(t)d^\beta_\alpha\right) x^\alpha  \in \R_{k-1}^{ [t,x] }\quad\mbox{in }Q_M.
		\end{equation*}
		
		We claim that, for any~$\alpha\in\N^n$ with~$|\alpha|\le d(k,s)$,
		\begin{equation}\label{claim:::poly_t}
			\text{$c_\alpha$ is a polynomial in~$t$ of degree at most~$d(k,s)-|\alpha|$.}
		\end{equation}		
		To check~\eqref{claim:::poly_t}, we argue by induction on~$d(k,s)-|\alpha|$.
		
		For any $\alpha\in\N^n$ such that $|\alpha|=d(k,s)$, we have that 
		\begin{equation*}
			 \partial_t c_\alpha(t) x^\alpha  \in \R_{k-1}^{ [t,x] }\quad\mbox{in }Q_M,
		\end{equation*}
		which, recalling that~$d(k,s)>k-1$, entails that~$\partial_t c_\alpha=0$, i.e., $c_\alpha$ is constant. Therefore, 
		\begin{equation}\label{cneguirtuifewjy698}
		{\mbox{\eqref{claim:::poly_t}
		is satisfied when~$|\alpha|=d(k,s)$.}}\end{equation}
		
		For $\alpha\in\N^n$ such that $|\alpha|=d(k,s)-1$, we distinguish between the cases~$s\in\left(0,\frac12\right]$ and~$s\in\left(\frac12,1\right)$.
		If~$s\in\left(0,\frac12\right]$, then~$d(k,s)=k$, thus $|\alpha|=k-1$ and, exploiting~\eqref{cneguirtuifewjy698},
		\begin{equation*}
			\left(\partial_tc_\alpha(t) +\sum_{\substack{\beta\in\N^n \\ |\alpha| +1 \leq |\beta| \leq d(k,s)}}c_\beta(t)d^\beta_\alpha\right) x^\alpha =\left(\partial_tc_\alpha(t) +\sum_{\substack{\beta\in\N^n\\ |\beta|=k}}c_\beta d^\beta_\alpha\right) x^\alpha  \in \R_{k-1}^{ [t,x] }\quad\mbox{in }Q_M
		\end{equation*}
		if and only if 
		\begin{equation*}
			\partial_tc_\alpha(t) +\sum_{\substack{\beta\in\N^n\\ |\beta|=k}}c_\beta d^\beta_\alpha \text{ is constant.}
		\end{equation*}
		
		If instead~$s\in\left(\frac12,1\right)$, we have that~$d(k,s)=k+1$, hence~$|\alpha|=k$ and, recalling~\eqref{cneguirtuifewjy698},
		\begin{equation*}
				\left(\partial_tc_\alpha(t) +\sum_{\substack{\beta\in\N^n\\ |\beta|=k+1}}c_\beta d^\beta_\alpha\right) x^\alpha  \in \R_{k-1}^{ [t,x] }\quad\mbox{in }Q_M
		\end{equation*}
		if and only if
		\begin{equation*}
			\partial_tc_\alpha(t) +\sum_{\substack{\beta\in\N^n\\ |\beta|=k+1}}c_\beta d^\beta_\alpha = 0.
		\end{equation*}
		In both cases, we infer that $c_\alpha(t)$ is a polynomial of degree at most $1$ in time, thus establishing~\eqref{claim:::poly_t}
		when~$|\alpha|=d(k,s)-1$.
		
		Assume now that 
		\begin{equation} \label{eq::inductive_step_liouville}
			\begin{split}
				&\text{$c_\beta$ is a polynomial in~$t$ of degree at most~$d(k,s)-|\beta|$,}\\
				&\text{for all $\beta\in\N^n$ such that $|\alpha|+1\leq |\beta|\leq d(k,s)$. }
			\end{split}
		\end{equation}
		Then,
		\begin{equation*}
			\left(\partial_tc_\alpha(t) +\sum_{\substack{\beta\in\N^n \\ |\alpha| +1 \leq |\beta| \leq d(k,s)}}c_\beta(t)d^\beta_\alpha\right) x^\alpha  \in \R_{k-1}^{ [t,x] }\quad\mbox{in }Q_M
		\end{equation*}
		yields
		\begin{equation}\label{chsdjbe93ueriheflqwiouwefio}
			 \partial_tc_\alpha(t) +\sum_{\substack{\beta\in\N^n \\ |\alpha| +1 \leq |\beta| \leq d(k,s)}}c_\beta(t)d^\beta_\alpha \text{ is a polynomial of degree at most $k-1-|\alpha|$ in time.}
		\end{equation}
		Moreover, by~\eqref{eq::inductive_step_liouville} 
		we know that
		\begin{equation}\label{chsdjbe93ueriheflqwiouwefio2}
		\sum_{\substack{\beta\in\N^n \\ |\alpha| +1 \leq |\beta| \leq d(k,s)}}c_\beta(t)d^\beta_\alpha \text{ is a polynomial of degree at most $d(k,s)-1-|\alpha|$ in time.}
		\end{equation}

Now, if~$s\in\left(0,\frac12\right]$, the statements in~\eqref{chsdjbe93ueriheflqwiouwefio}
and~\eqref{chsdjbe93ueriheflqwiouwefio2} give that	
		\begin{equation*}
			\partial_tc_\alpha(t) \text{ is a polynomial of degree at most $d(k,s)-|\alpha|-1$ in time,}
		\end{equation*}
		which entails~\eqref{claim:::poly_t}.
		
		If instead~$s\in\left(\frac12,1\right)$,
		we deduce from~\eqref{chsdjbe93ueriheflqwiouwefio}
and~\eqref{chsdjbe93ueriheflqwiouwefio2} that
$$ \partial_tc_\alpha(t) +\sum_{\substack{\beta\in\N^n \\ |\beta|=|\alpha|+1}}c_\beta(t)d^\beta_\alpha=0.
$$
Consequently, by~\eqref{eq::inductive_step_liouville},
\begin{equation*}
			\partial_tc_\alpha(t) \text{ is a polynomial of degree at most $d(k,s)-|\alpha|-1$ in time,}
		\end{equation*}
		which completes the proof of~\eqref{claim:::poly_t}.
		
		Now, notice that the conclusion of Lemma~\ref{lemma::partial_liouville} follows immediately from~\eqref{bcasr4365rfejkwguhwfo6789l98765} and~\eqref{claim:::poly_t} when ~$s\in\left(0,\frac12\right]$. If instead~$s\in\left(\frac12,1\right)$, then $d(k,s)=k+1$, and we write $u$ as
		\begin{equation*}
			\begin{split}
				u(t,x)&=\sum_{\substack{\alpha\in\N^n,\beta\in\N \\ |\alpha|+\beta \le k+1}} c_{\alpha,\beta}\,t^\beta\,x^\alpha \\
				&= \sum_{\substack{\alpha\in\N^n \\ |\alpha|=k+1}} c_{\alpha,0}\,x^\alpha + \sum_{\substack{\alpha\in\N^n,\beta\in\N\setminus\{0\} \\ |\alpha|+\beta=k+1}} c_{\alpha,\beta}\,t^\beta\,x^\alpha + P(t,x),
			\end{split}
		\end{equation*}
		where $c_{\alpha,\beta}\in\R$ and $P\in\R_{k}^{[t,x]}$.
		
		Observe that, for any $\alpha\in\N^n$ and $\beta\in\N\setminus\{0\}$ such that $|\alpha|+\beta =k+1$, 
		\begin{equation*}
			\left(\partial_t+(-\Delta)^s\right) (c_{\alpha,\beta}\,t^\beta\,x^\alpha) = \beta c_{\alpha,\beta}\,t^{\beta-1}\,x^\alpha + c_{\alpha,\beta}\,t^\beta\,(-\Delta)^sx^\alpha \in\R_{k}^{[t,x]}.
		\end{equation*}
		Therefore, since $u$ satisfies
		$$ \partial_t u + (-\Delta)^s u \eqps{k} 0\quad\mbox{in }Q_M,$$
		it follows that $c_{\alpha,\beta}=0$ for all $\alpha\in\N^n$ and $\beta\in\N\setminus\{0\}$ such that $|\alpha|+\beta =k+1$.
 		
		Consequently, we find that $u=H+P$, where
		\begin{equation*}
			H(x) := \sum_{\substack{\alpha\in\N^n \\ |\alpha|=k+1}} c_{\alpha,0}\,x^\alpha \in\R_{k+1}^{[x]}=\R_{d(k,s)}^{[x]}
		\end{equation*}
		and $P\in\R_{k}^{[t,x]}$, as desired.
		
		Moreover, since $u$ satisfies~\eqref{eq::t_poly_growth} by assumption, we infer that $u$ is a polynomial of degree at most $\lfloor\xi\rfloor$ in time.
	\end{proof}

	
	\appendix
	
	
	\section{Some technical auxiliary results} \label{sec::tech_lemmata}
	
Here we gather, for the convenience of the reader, some useful technical results from the literature, which establish general properties regarding incremental quotients and the H\"older seminorm modulo polynomials defined in~\eqref{eq::def_seminorm_mod_poly}. We refer to Appendices~A and~C of~\cite{elliptic_poly} for the proof of such results.
	
	Let us recall that for every given $\nu\in\S^{n-1}$ and $h\in\R$, the spatial incremental quotient of a function $\psi$ in direction $\nu$  is defined, at any point $x\in\R^n$, as 
	\begin{equation}
		D^\nu_h\psi(x) := \psi(x+h\nu)-\psi(x),
	\end{equation}
	and, for every $\zeta\in\S^{n-1}$, also $D^{(\nu,\zeta)}_h := D^\nu_h \circ D^\zeta_h$.
	
	Then, we have the following:
	\begin{lemma} \label{lemma::A1}
		There exists $\phi:\R^n\to[0,1]^d$ such that
		\begin{equation*} \label{eq::A1}
			h^{-d}D^{(\nu_1,\dots,\nu_d)}_h\psi(x) = \sum_{1\leq j_1,\dots,j_d\leq n} \frac{\partial^d \psi}{\partial x_{j_1}\dots\partial x_{j_d}}(x+\phi_1(x)\nu_1+\dots+h\phi_d(x)\nu_d)\nu_{1j_1}\dots\nu_{dj_d}.
		\end{equation*}
	\end{lemma}
	
	\begin{lemma} \label{lemma::A2}
		Let $\psi$, $\varphi:\R^n\to\R$ be such that
$$\int_{\R^n} \varphi D^{(\nu_1,\dots,\nu_d)}_h\psi \,dx+\int_{\R^n} \psi D^{(-\nu_1,\dots,-\nu_d)}_h\varphi \,dx <+\infty.
	$$
		
		Then,
		\begin{equation*}
			\int_{\R^n} \varphi D^{(\nu_1,\dots,\nu_d)}_h\psi \,dx = \int_{\R^n} \psi D^{(-\nu_1,\dots,-\nu_d)}_h\varphi \,dx.
		\end{equation*}
	\end{lemma}
	\medskip
	
	Now, we recall that the H\"older seminorm modulo polynomials is defined (see~\eqref{eq::def_seminorm_mod_poly} for the parabolic setting), for every function $\psi\in\cont^\gamma(B_1)$, as
	\begin{equation*}
		[\psi]_{\cont^\gamma(B_1;k)} :=  \inf_{P\in\R_{k-1}^{ [x] }} [\psi-P]_{\cont^\gamma(B_1)}.
	\end{equation*}
	Then, the following estimate holds true.
	
	\begin{lemma} \label{lemma::C1}
		Let $k\in\N$ and $\gamma\in(0,+\infty)\setminus\N$, with $\intpart{\gamma}\leq k-1$. Let also $\psi\in\cont^\gamma(B_1)$.
		
		Then, there exists a constant $C>0$, depending only on $n$, $k$, and $\gamma$, such that
		\begin{equation*}
			\inf_{P\in \R_{k-1}^{ [x] }} \left(\norm{\psi-P}_{L^\infty(B_1)} + [\psi-P]_{\cont^\gamma(B_1)}\right) \leq C [\psi]_{\cont^\gamma(B_1;k)}.
		\end{equation*}
	\end{lemma}


	\section{Proof of Proposition~\ref{prop::parabolic_3_1}} \label{sec::parabolic_3_1}
	
	In this section, we present an adaptation of the proof of~\cite[Theorem~3.1]{elliptic_poly} to the parabolic case.
	
	Let $\eta$ be as in the statement of Proposition~\ref{prop::parabolic_3_1}. Thus, for every $(t,x)\in Q_1$, we have
	\begin{equation} \label{eq::calc_laplace}
		\begin{split}
			&(-\Delta)^s(\eta u)(t,x)  = \int_{\R^n} \frac{\eta(x) u(t,x)-\eta(y) u(t,y)}{|x-y|^{n+2s}}\,dy \\
			&\qquad= \int_{B_2} \frac{ u(t,x)-u(t,y)}{|x-y|^{n+2s}}\,dy 
			+ \int_{B_2^c} \frac{ u(t,x)}{|x-y|^{n+2s}}\,dy 
			- \int_{B_2^c} \frac{ \eta(y)u(t,y)}{|x-y|^{n+2s}}\,dy\\
			&\qquad= f_{1,u}(t,x)+f_{2,u}(t,x) - \int_{B_2^c} \frac{ \eta(y) u(t,y)}{|x-y|^{n+2s}}\,dy.
		\end{split}
	\end{equation}
	
	Now, we rewrite
	\begin{equation*}
		\int_{B_2^c} \frac{ \eta(y) u(t,y)}{|x-y|^{n+2s}}\,dy 
		= \int_{B_2^c} \frac{ \eta(y) u(t,y)}{|y|^{n+2s}|x_y-y_y|^{n+2s}}\,dy
	\end{equation*}
	where $x_y:=x/|y|$ and $y_y:=y/|y|$. Moreover, considering a Taylor expansion (up to the $k$th-order) of $|x_y-y_y|^{-n-2s}$ for $x$ near the origin, we obtain
	\begin{equation*}
		\begin{split}
			&\frac{ \eta(y) u(t,y)}{|y|^{n+2s}|x_y-y_y|^{n+2s}}  = \frac{ \eta(y) u(t,y)}{|y|^{n+2s}}\left(\sum_{|\alpha|\leq k-1} c_{\alpha,y_y}x_y^\alpha + \sum_{|\alpha|=k} \varphi_\alpha(x_y,y_y)x_y^\alpha\right)\\
			&\qquad = \frac{ \eta(y) u(t,y)}{|y|^{n+2s}}\left(\sum_{|\alpha|\leq k-1} \frac{c_{\alpha,y_y}}{|y|^{|\alpha|}}x^\alpha + \sum_{|\alpha|=k} \frac{\varphi_\alpha(x_y,y_y)}{|y|^k}x^\alpha\right),
		\end{split}
	\end{equation*}
	with 
	\begin{equation} \label{eq::bound_coeff}
		\sup_{\substack{|\alpha|\leq k-1 \\ y\in\R^{n}\setminus\{0\}}} c_{\alpha,y_y} + \sup_{\substack{|\alpha|=k \\ y\in\R^{n}\setminus\{0\} \\ x\in B_{|y|/2}}} \partial_x^\zeta \varphi_\alpha(x_y,y_y) < +\infty.
	\end{equation}
	
	We define
	\begin{equation*}
		P_{u,\eta}(t,x) := -\sum_{|\alpha|\leq k-1} x^\alpha\int_{B_2^c} \frac{c_{\alpha,y_y}\eta(y)u(t,y)}{|y|^{n+2s+|\alpha|}}\,dy \in \left(\cont^{\sigma}((-1,0])\right)_{k-1}^{ [x] }
	\end{equation*}
	and 
	\begin{equation*}
		\psi(x,y) := - \sum_{|\alpha|=k} \varphi_\alpha(x_y,y_y)x^\alpha,
	\end{equation*}
	so that
	\begin{equation} \label{eq::P_and_f_star}
		-\int_{B_2^c} \frac{ \eta(y) u(t,y)}{|x-y|^{n+2s}}\,dy = P_{u,\eta}(t,x) + \int_{B_2^c} \frac{ \eta(y) u(t,y)\psi(x,y)}{|y|^{n+2s+k}}\,dy = P_{u,\eta}(t,x) + f_{u,\eta}^*(t,x).
	\end{equation}
	Moreover,~\eqref{eq::bound_coeff} yields~\eqref{eq::estimate_psi}.
	
	Finally, putting~\eqref{eq::calc_laplace} and~\eqref{eq::P_and_f_star} together, we infer
	\begin{equation*}
		\partial_t(\eta u) + (-\Delta)^s(\eta u) = \partial_t(\eta u) + f_{1,u}+f_{2,u} + P_{u,\eta}(t,x) + f_{u,\eta}^*(t,x) =: f_{u,\eta}+ P_{u,\eta}(t,x),
	\end{equation*}
	as desired.
	
	
	\section{Proof of Corollary~\ref{cor::parabolic_3_2}} \label{sec::parabolic_3_2}
	
	Here, we show Corollary~\ref{cor::parabolic_3_2}, by following the ideas of~\cite[Corollary~3.2]{elliptic_poly}.
	
	Let us pick $R>0$. Thus, thanks to Proposition~\ref{prop::parabolic_3_1} applied with $\eta:=\eta_R$, we infer that
	\begin{equation*}
		(-\Delta)^s(\eta u) = f_{1,u}+f_{2,u} + P_{u,\eta}(t,x) + f_{u,\eta}^*(t,x)
	\end{equation*}	
	where $P_{u,\eta_R}\in \left(\cont^{\sigma}((-1,0])\right)_{k-1}^{ [x] }$ and
	\begin{equation*}
		f_{u,\eta_R}^*(t,x) := \int_{B_2^c} \frac{ \eta_R(y) u(t,y)\psi(x,y)}{|y|^{n+2s+k}}\,dy.
	\end{equation*}
	
	Moreover, recalling~\eqref{eq::estimate_psi}, we obtain
	\begin{equation*}
		\frac{ |\eta_R(y) u(t,y)\psi(x,y)|}{|y|^{n+2s+k}} \leq \frac{C|u(t,y)|}{|y|^{n+2s+k}} \in L^1(B_2^c),
	\end{equation*}
	for any $t\in(-1,0]$. Therefore, by the Dominated Convergence Theorem, we deduce that, for every $(t,x)\in Q_1$,
	\begin{equation} \label{eq::lim_f3}
		\lim_{R\to+\infty} f_{u,\eta_R}^*(t,x) = \int_{B_2^c} \frac{u(t,y)\psi(x,y)}{|y|^{n+2s+k}}\,dy =: f_{3,u}(t,x).
	\end{equation}
	
	Besides, since $\eta_R\to1$ a.e. in $\R^n$ as $R\to+\infty$, for every $(t,x)\in Q_1$, we have that
	\begin{equation} \label{eq::conv_Dt}
		\lim_{R\to+\infty} \partial_t (\eta_R u) = \lim_{R\to+\infty} \eta_R \partial_t u = \partial_t u.
	\end{equation}
	
	Thus, putting~\eqref{eq::lim_f3} and~\eqref{eq::conv_Dt} together,we conclude that 
	\begin{equation*}
		\begin{split}
			&\lim_{R\to+\infty} \left(\partial_t (\eta_R u) + (-\Delta)^s(\eta u) - P_{u,\eta_R}\right) = \lim_{R\to+\infty} \left(\partial_t (\eta_R u) + f_{1,u}+f_{2,u} + f_{u,\eta_R}^*\right)\\
			&\qquad =  \partial_t u + f_{1,u}+f_{2,u} + f_{3,u} =: f_u.
		\end{split}
	\end{equation*}
	
	
	\section{Proof of Corollary~\ref{cor::parabolic_3_7}} \label{sec::parabolic_3_7}
	
	In this section, we provide a proof of Corollary~\ref{cor::parabolic_3_7}. To this aim, for every $R>1$, let us consider the
	smooth cut-off function
	\begin{equation*}
		\eta_R :=
		\begin{cases}
			1\quad&\mbox{in }B_R,\\
			0\quad&\mbox{in }B_{R+1}^c.
		\end{cases}
	\end{equation*}
Also, we define $v:=(1-\eta_R)u$, which satisfies
	\begin{align*}
		&v=0\quad\mbox{in }Q_1,\\
		&|v|\leq |u|,\\
		&v\in\cont^1((-1,0];\cont^\alpha(B_1)),\quad\mbox{for some }\alpha>2s,\\
		\mbox{and}\quad&J^k_1(v)\leq J^k_1(u).
	\end{align*}
	
	Moreover, for any $R>2$, we have 
	\begin{equation*}
		w:=(\eta_R-\eta_2)u=(\eta_R-\eta_2)v.
	\end{equation*}
	Therefore, recalling Definition~\ref{def::k_viscosity_sol}, we find
	\begin{equation} \label{eq::cor_35_eq_w_u}
\partial_t w + (-\Delta)^s w = \partial_t \big((\eta_R-\eta_2)u \big)+ (-\Delta)^s\big( (\eta_R-\eta_2)u\big) = f_R-f_2+\widetilde{P}_R,
	\end{equation}
	for some $\widetilde{P}_R\in\left(\cont^{\sigma}((-1,0])\right)_{k-1}^{ [x] }$.
	
	On the other hand, by Proposition~\ref{prop::parabolic_3_1}, we have
	\begin{equation} \label{eq::cor_35_eq_w_v}
	\partial_t w + (-\Delta)^s w = \partial_t (\eta_R-\eta_2)v + (-\Delta)^s (\eta_R-\eta_2)v
			 = P_{v,R} + \int_{B_R\setminus B_2} \frac{u(t,y)\psi(x,y)}{|y|^{n+2s+k}}\,dy,
	\end{equation}
	with $P_{v,R}\in\left(\cont^{\sigma}((-1,0])\right)_{k-1}^{ [x] }$.
	
	Thus, comparing~\eqref{eq::cor_35_eq_w_u} and~\eqref{eq::cor_35_eq_w_v}, we deduce that
	\begin{equation} \label{eq::comparison_w_u_v}
		f_R-f_2+ \widehat{P}_R = \int_{B_R\setminus B_2} \frac{u(t,y)\psi(x,y)}{|y|^{n+2s+k}}\,dy,
	\end{equation}
	where $ \widehat{P}_R :=  \widetilde{P}_R - P_{v,R}$.
	
	Let us now define, for every $r$, $R>0$, the function
	\begin{equation*}
		\Psi_{R,r} :=  \int_{B_R\setminus B_r} \frac{u(t,y)\psi(x,y)}{|y|^{n+2s+k}}\,dy.
	\end{equation*}
	 It follows from~\eqref{eq::comparison_w_u_v} that
	 \begin{equation*}
	 	\begin{split}
	 		&\sup_{t\in(-1,0]}  \norm{D^m_x(f_R+ \widehat{P}_R -f_r- \widehat{P}_r)}_{L^\infty(B_1)} = \sup_{t\in(-1,0]} \norm{D^m_x\Psi_{R,r}}_{L^\infty(B_1)} \\
	 		&\qquad \leq \sup_{t\in(-1,0]}  \norm{\Psi_{R,r}}_{\cont^m(B_1)} \leq \sup_{t\in(-1,0]} \int_{B_r^c} \frac{|u(t,y)|}{|y|^{n+2s+k}}\,dy,
	 	\end{split}
	 \end{equation*}
	 for every $m\geq0$, which gives~\eqref{eq::parabolic_3_7}, as desired.

	
	\section{Proof of Corollary~\ref{cor::schauder_seminorm}} \label{sec::parabolic_4_1}
	
	This section is devoted to proving Corollary~\ref{cor::schauder_seminorm}. 
	
	Arguing by contradiction, we suppose that there exist sequences $\{u_j\}_j$ and $\{f_j\}_j$ such that
	\begin{align*}
		&\partial_t u_j + (-\Delta)^s u_j = f_j \quad\mbox{in }Q_1\\
		\mbox{and}\quad&\Theta_j:=\norm{u_j}_{L^\infty(Q_{9/10})}+\sup_{t\in\left(-\left(\frac{99}{100}\right)^{2s},0\right]} \norm{f_j}_{\cont^m(B_{99/100})} \\
		&\qquad>j\left(J^0_{1,r}(u_j)+\sup_{t\in(-1,0]}[f_j]_{\cont^{m+\theta}(Q_1)}+\norm{u_j}_{L^\infty(\overline{Q}^*_1)}\right).
	\end{align*}
	
	Let also
	\begin{equation*}
		\widetilde{u}_j:=\frac{u_j}{\Theta_j}\quad\mbox{and}\quad\widetilde{f}	_j:=\frac{f_j}{\Theta_j}.
	\end{equation*}
	
	Therefore, we have that 
	\begin{align}
		\label{eq::normalized_eq}&\partial_t \widetilde{u}_j + (-\Delta)^s \widetilde{u}_j = \widetilde{f}_j \quad\mbox{in }Q_1,\\
		\label{eq::small_tail}&J^0_{1,r}(\widetilde{u}_j) = \frac{J^0_{1,r}(u_j)}{\Theta_j} \leq \frac{1}{j},\\
		\label{eq::small_boundary}&\norm{u_j}_{L^\infty(\overline{Q}^*_1)}\leq\frac{1}{j},\\
		\label{eq::tilde_norm}& \sup_{t\in\left(-\left(\frac{99}{100}\right)^{2s},0\right]} \norm{\widetilde{f}_j}_{\cont^m(B_{99/100})}\leq 1,\\
		\label{eq::tilde_seminorm}\mbox{and}\quad&\sup_{t\in(-1,0]} [\widetilde{f}_j]_{\cont^{m+\theta}(B_1)}\leq \frac{1}{j}.
	\end{align}
	
	In particular, $\sup_{t\in\left(-\left(\frac{99}{100}\right)^{2s},0\right]}\norm{\widetilde{f}_j}_{\cont^{m+\theta}(B_{99/100})}$ is bounded uniformly in $j$, and hence there exists
	\begin{equation*}
		\widetilde{f}\in\cont^0\left(\left(-\left(\frac{99}{100}\right)^{2s},0\right];\cont^{m+\theta}(B_{99/100})\right)
	\end{equation*} 
	such that $\widetilde{f}_j\to\widetilde{f}$ uniformly in $Q_{99/100}$, as $j\to+\infty$.
	
	Moreover, 
	\begin{equation*}
		\begin{split}
			&\sup_{t\in\left(-\left(\frac{99}{100}\right)^{2s},0\right]} [\widetilde{f}]_{\cont^{m+\theta}(B_{99/100})}\\ 
			&\qquad\leq \sup_{t\in\left(-\left(\frac{99}{100}\right)^{2s},0\right]}[\widetilde{f}-\widetilde{f}_j]_{\cont^{m+\theta}(B_{99/100})} +\sup_{t\in\left(-\left(\frac{99}{100}\right)^{2s},0\right]} [\widetilde{f}_j]_{\cont^{m+\theta}(B_{99/100})} \\
			&\qquad\leq \sup_{t\in\left(-\left(\frac{99}{100}\right)^{2s},0\right]} [\widetilde{f}-\widetilde{f}_j]_{\cont^{m+\theta}(B_{99/100})} +\frac{1}{j} \xrightarrow[j\to+\infty]{}0,
		\end{split}
	\end{equation*}
	which entails that~$\widetilde{f}\in \left(\cont^0\left(\left(-\left(\frac{99}{100}\right)^{2s},0\right]\right)\right)_m^{ [x] }$ in $Q_{99/100}$.
	
	In addition,~\eqref{eq::small_tail} and~\eqref{eq::small_boundary} together imply that $\widetilde{u}_j\to0$ as $j\to+\infty$ in $(-1,0]\times B_{1/2}^c\cup \overline{Q}^*_1$.
	
	Therefore, from Theorem~\ref{th::rescaled_known_schauder}, we deduce that $\widetilde{u}_j\in\cont^\alpha(Q_{9/10})$, for some $\alpha\in(0,1)$, and $\norm{\widetilde{u}_j}_{\cont^{\alpha}(Q_{9/10})}$ is uniformly bounded in $j$. Thus, we find a function $\widetilde{u}\in\cont^\alpha(Q_{9/10})$ such that
	\begin{align*}
		&\lim_{j\to+\infty}\norm{\widetilde{u}-\widetilde{u_j}}_{L^\infty(Q_{9/10})}+ J^0_{9/10,r}(\widetilde{u}-\widetilde{u}_j) = 0,\\
		&\partial_t \widetilde{u}+ (-\Delta)^s \widetilde{u} = \widetilde{f}\quad\mbox{in }Q_{9/10},\\
		\mbox{and}\quad&\widetilde{u}=0\quad\mbox{in }\left(\left(-\left(\frac{9}{10}\right)^{2s},0\right]\times B_{1/2}^c\right)
		\cup \left(\left[-\left(\frac{9}{10}\right)^{2s},-\frac{1}{2}\right]\times B_1\right).
	\end{align*}
	Recalling Lemma~\ref{lemma::max_principle}, we infer that $\widetilde{u}\equiv0$ in $Q_{9/10}$. Consequently $\widetilde{f}\equiv0$ in $Q_{9/10}$, as well. Since $\widetilde{f}$ is a polynomial, we conclude that~$\widetilde{f}\equiv0$ also in the larger cylinder $Q_{99/100}$.
	
	Thus, we obtain that
	\begin{equation*}
		\begin{split}
			&1 = \lim_{j\to+\infty} \frac{\norm{u_j}_{L^\infty(Q_{9/10})}+\sup_{t\in\left(-\left(\frac{99}{100}\right)^{2s},0\right]} \norm{f_j}_{\cont^{m}(B_{99/100})}}{\Theta_j} \\
			&\qquad = \lim_{j\to+\infty}\norm{\widetilde{u_j}}_{L^\infty(Q_{9/10})}+\sup_{t\in\left(-\left(\frac{99}{100}\right)^{2s},0\right]} \norm{\widetilde{f}_j}_{\cont^{m}(B_{99/100})}\\
			&\qquad = \lim_{j\to+\infty}\norm{\widetilde{u_j}-\widetilde{u}}_{L^\infty(Q_{9/10})}+\sup_{t\in\left(-\left(\frac{99}{100}\right)^{2s},0\right]} \norm{\widetilde{f}_j-\widetilde{f}}_{\cont^{m}(B_{99/100})} = 0,
		\end{split}
	\end{equation*}
	which leads to a contradiction, concluding the proof of Corollary~\ref{cor::schauder_seminorm}.
	\medskip
	
	As we think it may be of independent interest, we stress that an analogous argument, up to replacing the term $\sup_{t\in\left(-\left(\frac{99}{100}\right)^{2s},0\right]} \norm{f}_{\cont^m(B_{99/100})}$ with $\norm{f}_{\cont^m(Q_{99/100})}$ in~\eqref{eq::schauder_seminorm}, proves the following:
	
	\begin{cor}
		Let $s\in(0,1)$,~$\xi\geq\sigma$, and~$r\in (0,9/10)$. Let also~$m\in\N$ and $\theta\in(0,1)$. Let~$u:Q_\infty\to\R$ and~$f\in\cont^{m+\theta}(Q_\infty)$.
		
		Suppose that $u\in\cont^0(Q_1)\cap V_0^\xi((-1,0])$ is a solution of
		\begin{equation*}
			\partial_t u +(-\Delta)^s u= f\quad\mbox{in }Q_1
		\end{equation*}
		in the viscosity sense.
		
		Then, there exists a constant $C>0$ depending only on $n$, $s$, $m$, and $\xi$, such that
		\begin{equation*}
\norm{u}_{L^\infty(Q_{9/10})}+\norm{f}_{\cont^m(Q_{99/100})}
			 \leq C\left(J^0_{1,r}(u)+[f]_{\cont^{m+\theta}(Q_1)}+\norm{u}_{L^\infty(\overline{Q}^*_1)}\right).
		\end{equation*}
	\end{cor} 
	
	
	\section{Proof of Proposition~\ref{prop::known_suboptimal_schauder}} \label{sec::known_suboptimal_schauder}
	
	For the sake of completeness, we sketch here a proof of Proposition~\ref{prop::known_suboptimal_schauder}. Let $\gamma=m+\theta$, with $m\in\N$ and $\theta\in(0,1)$. In the wake of~\cite[Lemma~4.2]{elliptic_poly}, the argument goes by induction over $m$. Also, up to rescalings, we assume $R=1$. For simplicity\footnote{
	In the case $\theta+2s\leq1$, the argument still works, but it has to be carried over step-by-step. For more details about this procedure, see~\cite[Lemma~4.2]{elliptic_poly}.
	}, we deal with the case $\theta+2s>1$.
	
	For $m=0$, the conclusion follows from~\cite[Theorem~5.1]{parabolic_viscosity}. Then, let us suppose that
	\begin{equation} \label{eq::known_suboptimal_schauder2}
		\sup_{t\in(-1/2,0]}\norm{u(t,\cdot)}_{\cont^{m+\theta+2s}(B_{1/2})} \leq C\left(J^0_1(u)+\sup_{t\in(-1,0]}\norm{f}_{\cont^{m+\theta}(B_1)} +\norm{u}_{L^\infty(Q_1)}\right),
	\end{equation}
	holds true for some $m\in\N$, and we show the claim for $m+1$.
	
	To this purpose, we recall that, for every given $\nu\in\S^{n-1}$ and $h\in\R$, the spatial incremental quotient of~$u$ in direction $\nu$  is defined, at any point $(t,x)\in Q_\infty$, as 
	\begin{equation}
		D^\nu_hu(t,x) := u(t,x+h\nu)-u(t,x),
	\end{equation}
	and, for every $\zeta\in\S^{n-1}$, also $D^{(\nu,\zeta)}_h := D^\nu_h \circ D^\zeta_h$.
	
	Then, we consider $\nu_1,\dots,\nu_{m+1}\in\S^{n-1}$, and functions $v:=D^{(\nu_1,\dots,\nu_{m+1})}_h u$ and $g:=D^{(\nu_1,\dots,\nu_{m+1})}_h f$. Thus, for $h$ small enough, we have
	\begin{equation*}
		\partial_t v + (-\Delta)^s v = g\quad\mbox{in }Q_{9/10},
	\end{equation*}
	in the viscosity sense.
	
	Let now $\varphi\in\cont^\infty_c(B_{1/4})$ be a smooth cut-off function such that $\varphi=1$ in $B_{1/8}$, and let $w(t,x):=\varphi(x)v(t,x)$. Therefore, for any $(t,x)\in Q_{1/16}$, we obtain
	\begin{equation} \label{eq::eq_for_w}
		\partial_t w(t,x) + (-\Delta)^s w(t,x) = g(t,x) + \int_{\R^n} \frac{(1-\varphi(y))v(t,y)}{|x-y|^{n+2s}}\,dy
	\end{equation}
	Moreover, for all $x\in B_{1/16}$, we have that $\displaystyle \frac{1-\varphi(y)}{|x-y|^{n+2s}}\in L^1(\R^n)$. Hence, by~Lemma~\ref{lemma::A2},
	\begin{equation} \label{eq::integration_by_parts_quot}
		\begin{split}
			&\int_{\R^n} \frac{(1-\varphi(y))v(t,y)}{|x-y|^{n+2s}}\,dy = \int_{\R^n} \frac{(1-\varphi(y))D^{(\nu_1,\dots,\nu_{m+1})}_h u(t,x)}{|x-y|^{n+2s}}\,dy\\
			&\qquad= \int_{\R^n}D^{(-\nu_1,\dots,-\nu_{m+1})}_h \left( \frac{1-\varphi(y)}{|x-y|^{n+2s}}\right) u(t,x)\,dy =: q(t,x).
		\end{split}
	\end{equation}
	From~\eqref{eq::eq_for_w} and~\eqref{eq::integration_by_parts_quot}, we infer that
	\begin{equation*}
		\partial_t w + (-\Delta)^s w = g + q =: F \quad\mbox{in }Q_{1/16},
	\end{equation*}
	in the viscosity sense.
	
	By virtue of~\eqref{eq::known_suboptimal_schauder2} with $m=0$ and~Lemma~\ref{lemma::A1}, we find
	\begin{equation} \label{eq::inductive_step_knwon_schauder}
		\begin{split}
			&\sup_{t\in(-1/32,0]}\norm{D^{(\nu_1,\dots,\nu_{m+1})}_h u}_{\cont^{\theta+2s}(B_{1/32})} = \sup_{t\in(-1/32,0]}\norm{w}_{\cont^{\theta+2s}(B_{1/32})}\\
			&\qquad\qquad\leq C\left(J^0_{1/16}(w)+\sup_{t\in(-1/16,0]}\norm{F}_{\cont^\theta(B_{1/16})} + \norm{w}_{L^\infty(Q_{1/16})}\right). 
		\end{split}
	\end{equation}
	
	Notice also that, thanks to~Lemma~\ref{lemma::A1}, we have that
	\begin{equation*}
		\sup_{t\in(-1/16,0]}\norm{q}_{\cont^1(B_{1/16})} \leq Ch^{m+1}\left(J^0_{1/4}(u)+\norm{u}_{L^\infty(Q_{1/8})}\right)
	\end{equation*}
	and
	\begin{equation} \label{eq::estimate_F}
		\begin{split}
			&\sup_{t\in(-1/16,0]}\norm{F}_{\cont^\theta(B_{1/16})} \leq \sup_{t\in(-1/16,0]}\norm{g}_{\cont^\theta(B_{1/16})} +\sup_{t\in(-1/16,0]}\norm{q}_{\cont^\theta(B_{1/16})} \\
			&\qquad\qquad \leq  Ch^{m+1}\left(J^0_{1/4}(u)+\sup_{t\in(-1,0]}\norm{f}_{\cont^{m+1+\theta}(B_{1/8})}+\norm{u}_{L^\infty(Q_{1/8})}\right),
		\end{split}
	\end{equation}
	for some positive constant $C$, possibly changing from line to line.
	
	In addition, using again~Lemma~\ref{lemma::A1}, we deduce that
	\begin{equation} \label{eq::estimate_sup_w}
		J^0_{1/16}(w) \leq C\norm{w}_{L^\infty(Q_{1/16})} \leq C\sup_{t\in(-1/4,0]}\norm{u}_{\cont^{m+1}(B_{1/4})}.
	\end{equation}
	
	Therefore,~\eqref{eq::inductive_step_knwon_schauder} together with~\eqref{eq::estimate_F} and~\eqref{eq::estimate_sup_w} yields
	\begin{equation*}
		\begin{split}
			&\sup_{t\in(-1/32,0]}\norm{D^{(\nu_1,\dots,\nu_{m+1})}_h u}_{\cont^{\theta+2s}(B_{1/32})} \\
			&\quad\leq Ch^{m+1} \left(J^0_{1/4}(u)+\sup_{t\in(-1,0]}\norm{f}_{\cont^{m+1+\theta}(B_{1/8})} +\norm{u}_{L^\infty(Q_{1/8})}+\sup_{t\in(-1/4,0]}\norm{u}_{\cont^{m+1}(B_{1/4})}\right).
		\end{split}
	\end{equation*}
	Dividing by $h^{m+1}$ and taking the limit as~$h\to0$, it thereby follows that
	\begin{equation} \label{eq::estimate_last_derivative}
		\begin{split}
			&\sup_{t\in(-1/32,0]}\norm{D^{m+1}_x u}_{\cont^{\theta+2s}(B_{1/32})} \\
			&\quad\leq C \left( J^0_{1/4}(u)+\sup_{t\in(-1,0]}\norm{f}_{\cont^{m+1+\theta}(B_{1/8})} +\norm{u}_{L^\infty(Q_{1/8})}+\sup_{t\in(-1/4,0]}\norm{u}_{\cont^{m+1}(B_{1/4})}\right),
		\end{split}
	\end{equation}
	
	Now, observe that
	\begin{equation*}
		J^0_{1/4}(u) \leq J^0_{1}(u) + C\norm{u}_{L^\infty(Q_1)}.
	\end{equation*}
	up to renaming $C$.
	
	From this,~\eqref{eq::known_suboptimal_schauder2},
	\eqref{eq::estimate_last_derivative}, and the fact that $\theta+2s>1$, we conclude that
	\begin{equation*}
		\sup_{t\in(-1/32,0]}\norm{u}_{\cont^{m+1+\theta+2s}(B_{1/32})} \leq C \left( J^0_{1}(u)+\sup_{t\in(-1,0]}\norm{f}_{\cont^{m+1+\theta}(B_{1})}  +\norm{u}_{L^\infty(Q_{1})}\right),
	\end{equation*}
	as desired.
	
	
	\section{A counterexample to time-regularity} \label{sec::counterex}
	
	In this section, we exhibit a counterexample to time-regularity of solutions to a Dirichlet problem in which the exterior datum is not regular in time in the first place. In~\cite{MR3115838}, the authors showed that, in general, a solution to
	\begin{equation} \label{eq::dirichlet_problem}
		\begin{cases}
			(-\Delta)^s u = 0\quad&\text{in }Q_1,\\
			u(0,\cdot) = 0\quad&\text{in }\{0\}\times\R^n,\\
			u = g\quad&\text{in }(-1,0]\times B_1^c,
		\end{cases}
	\end{equation}
	is no better than Lipschitz if $g$ is discontinuous in time. Here, we show that, in general, a viscosity solution to~\eqref{eq::dirichlet_problem} is no better than $\cont^0$ for a suitable choice of $g$.
	
	To this aim, let us define, for any $k\in\N$, sets
	\begin{equation*}
		I_k:=(-2^{-k},-2^{-k-1})\quad\text{and}\quad C_k:=B_{2^{k+2}}\setminus B_{2^{k+2}-1}.
	\end{equation*}
	
	Let $u:Q_\infty\to\R$ be a (viscosity) solution to~\eqref{eq::dirichlet_problem}, for a suitable function $g$. 
	
	Now, let $\gamma\in(0,1)$ and consider sequences of real numbers $\{h_k\}_k$ and $\{d_k\}_k$ to be made precise later. Let also $\varphi\in\cont^\infty_c(B_2)$ be such that $\varphi\equiv 1$ in $B_1$ and $\norm{\nabla \varphi}_{\infty}\leq c_\varphi$, for some $c_\varphi>0$. Then, we define functions $v:Q_\infty\to\R$ and $w:Q_\infty\to\R$ as
	\begin{equation*}
		\begin{split}
		v(t,x) &:= h_0\left(t+\frac{1}{2}\right)+2^{(2s+\gamma)}\chi_{C_0}(x)\\
		&\qquad+ \sum_{k=1}^{+\infty} (-1)^{k}\left(h_k\left(t+\frac{3}{2^{k+2}}\right)\varphi(x) + 2^{k(2s+\gamma)}\chi_{C_k}(x)\right)\chi_{I_k}(t)
		\end{split}
	\end{equation*}
	and
	\begin{equation*}
		\begin{split}
			w(t,x) &:=d_0\left(t+\frac{1}{2}\right)+2^{(2s+\gamma)}\chi_{C_0}(x)\\
			&\qquad+\sum_{k=1}^{+\infty} (-1)^{k} \left(d_k\left(t+\frac{3}{2^{k+2}}\right)\varphi(x) + 2^{k(2s+\gamma)}\chi_{C_k}(x)\right)\chi_{I_k}(t).
		\end{split}
	\end{equation*}
	
	Moreover, we choose
	\begin{equation*}
		g(t,x) = v(t,x) = w(t,x) = \sum_{k=1}^{+\infty} (-1)^{k+1} 2^{k(2s+\gamma)}\chi_{C_k}(x)\chi_{I_k}(t)\quad\text{in }(-1,0)\times B_1^c.
	\end{equation*}
	
	Computing the fractional Laplacian of $v$ in $B_1$, for any fixed time $t\in I_k$, with odd $k$, we obtain
	\begin{equation} \label{eq::v_laplace_odd}
		\begin{split}
			&(-\Delta)^sv(x,t) = \int_{\R^n} \frac{v(x,t)-v(y,t)}{|x-y|^{n+2s}}\, dy \\
			&\quad = h_k\left(t+\frac{3}{2^{k+2}}\right)\int_{\R^n} \frac{\varphi(x)-\varphi(y)}{|x-y|^{n+2s}}\, dy -2^{k(2s+\gamma)}\int_{C_k} \frac{dy}{|x-y|^{n+2s}} \\
			&\quad = h_k\left(t+\frac{3}{2^{k+2}}\right)\left( \int_{B_2\setminus B_1} \frac{1-\varphi(y)}{|x-y|^{n+2s}}\, dy + \int_{B_2^c} \frac{dy}{|x-y|^{n+2s}}\right) -2^{k(2s+\gamma)}\int_{C_k} \frac{dy}{|x-y|^{n+2s}} \\
			&\quad \leq c_1 2^{-k-2}h_k - c_2 2^{k(\gamma-n)}|C_k| \leq c_12^{-k-2} h_k - c_2 2^{k\gamma},
		\end{split}
	\end{equation}
	for some positive constant $c_1=c_1(n,s,c_\varphi)$ and $c_2=c_2(n,s)$, possibly changing from line to line.
	
	From this, choosing  
	$$h_k := \frac{c_2}{1+2^{-k-2}c_1}2^{k\gamma-1}, \quad\text{for $k\in\N$ odd},$$
	we infer that
	\begin{equation} \label{eq::v_subsol_1}
		\partial_t v + (-\Delta)^s v \leq (1+2^{-k-2}c_1)h_k - c_2 2^{k\gamma} \leq 0.
	\end{equation}
	
	Instead, when $k$ is even, a similar computation to~\eqref{eq::v_laplace_odd} yields 
	\begin{equation*}
		\begin{split}
			&(-\Delta)^sv(x,t) = -h_k\left(t+\frac{3}{2^{k+2}}\right)\int_{\R^n} \frac{\varphi(x)-\varphi(y)}{|x-y|^{n+2s}}\, dy +2^{k(2s+\gamma)}\int_{C_k} \frac{dy}{|x-y|^{n+2s}} \\
			&\quad = h_k\left(t+\frac{3}{2^{k+2}}\right)\left( \int_{B_2\setminus B_1} \frac{1-\varphi(y)}{|x-y|^{n+2s}}\, dy + \int_{B_2^c} \frac{dy}{|x-y|^{n+2s}}\right) -2^{k(2s+\gamma)}\int_{C_k} \frac{dy}{|x-y|^{n+2s}} \\
			&\quad \leq -c_1 2^{-k-2}h_k + c_3 2^{k\gamma},
		\end{split}
	\end{equation*}
	for some constant $c_3$ depending only on $n$ and $s$. Therefore, taking 
	$$h_k := \frac{c_3}{1+2^{-k-2}c_1}2^{k\gamma+1}, \quad\text{for $k\in\N$ even},$$
	we obtain that
	\begin{equation} \label{eq::v_subsol_2}
		\partial_t v + (-\Delta)^s v \leq -(1+c_1)h_k + c_3 2^{k\gamma} \leq 0\quad\text{in }Q_1.
	\end{equation}
	From~\eqref{eq::v_subsol_1} and~\eqref{eq::v_subsol_2}, we deduce that $v$ is a (viscosity) subsolution to~\eqref{eq::dirichlet_problem}, and hence $u\geq v$ in $Q_1$.
	
	On a similar note, whenever
	\begin{equation*}
		\begin{cases}
			\displaystyle d_k \geq \frac{c_2}{1-2^{-k-2}c_1} 2^{k\gamma},\quad&\text{when $k$ is odd},\\
			\\
			\displaystyle d_k \leq \frac{c_3}{1-2^{-k-2}c_1} 2^{k\gamma},\quad&\text{when $k$ is even},
		\end{cases}
	\end{equation*}
	we have that
	\begin{equation*}
		\partial_t w + (-\Delta)^s w = (1-c_1)d_k - c_3 2^{k\gamma} \geq 0\quad\text{in }Q_1.
	\end{equation*}
	Hence, $w$ is a supersolution to~\eqref{eq::dirichlet_problem}. 
	
	Therefore, choosing
	\begin{equation*}
		d_k:=
		\begin{cases}
			c_2 2^{k\gamma+1},\quad&\text{when $k$ is odd},\\
			c_3 2^{k\gamma-1},\quad&\text{when $k$ is even},
		\end{cases}
	\end{equation*}
	we obtain that $v\leq u\leq w$ in $Q_1$.
	
	Moreover, observe that $h_k$, $d_k\to+\infty$ as $k\to+\infty$, so that $u$ is at most $\cont^0$ in $Q_1$.
	
	
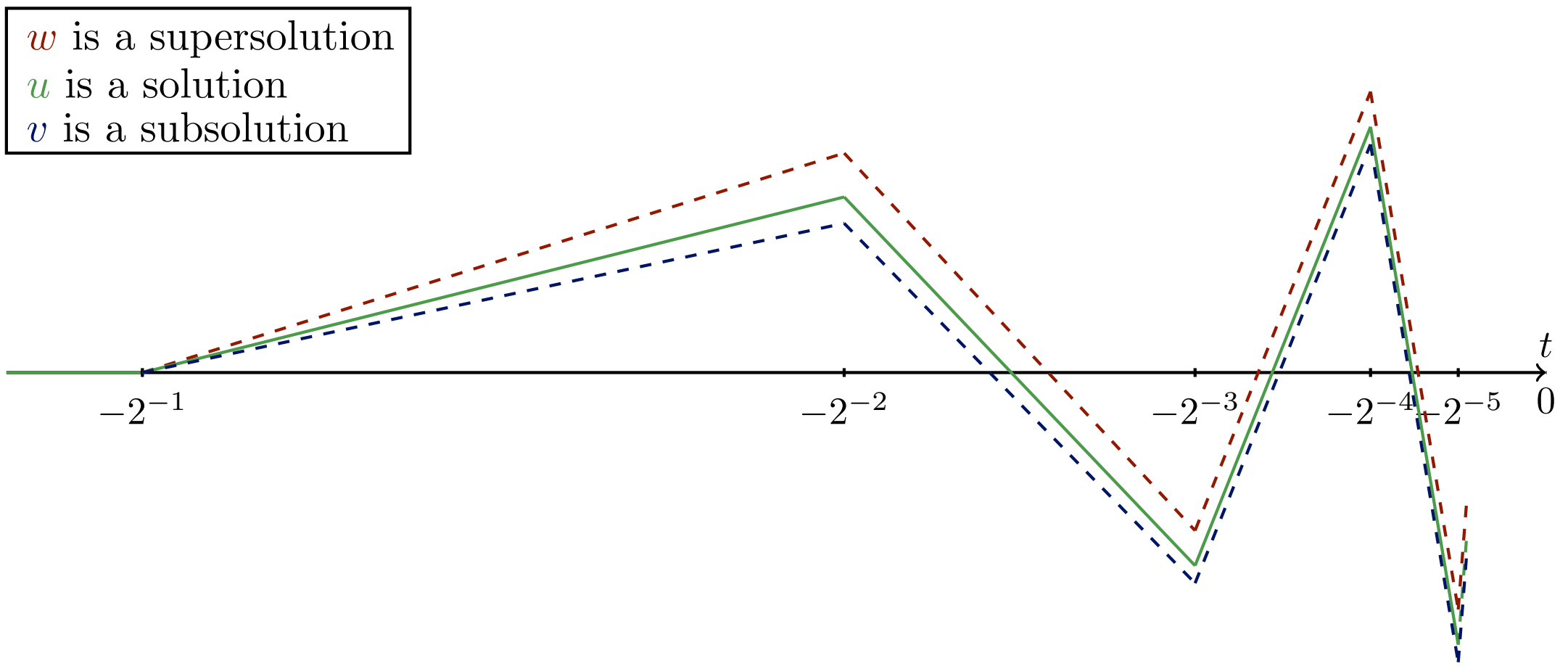
\begin{figure}[htbp]
		\definecolor{subsol}{RGB}{0,20,100}
		\definecolor{sol}{RGB}{75,155,75}
		\definecolor{supersol}{RGB}{145,25,0}
		\centering
		\begin{tikzpicture}[scale=0.8, thick]
			\draw[->] (-1.55,0) -- (16,0) 
			node[right, above, font=\small] {$t$};
			
			\draw (0.0, 0.05) -- (0.0, -0.05) node[below, font=\small] {$-2^{-1}$};
			\draw (8.0, 0.05) -- (8.0, -0.05) node[below, font=\small] {$-2^{-2}$};
			\draw (12.0, 0.05) -- (12.0, -0.05) node[below, font=\small] {$-2^{-3}$};
			\draw (14.0, 0.05) -- (14.0, -0.05) node[below, font=\small] {$-2^{-4}$};
			\draw (15.0, 0.05) -- (15.0, -0.05) node[below, font=\small] {$-2^{-5}$};
			\draw (16, 0) -- (16, 0) node[right, below, font=\small] {$0$};
			
			\draw[sol] (-1.55,0) -- (-0.02,0);
			
			\draw[supersol, dashed] (0,0) -- (8.0, 2.5);
			\draw[sol] (0,0) -- (8.0, 2.0);
			\draw[subsol, dashed] (0,0) -- (8.0, 1.7);
			
			\draw[supersol, dashed] (8.0, 2.5) -- (12.0, -1.8);
			\draw[sol] (8.0, 2.0) -- (12.0, -2.2);
			\draw[subsol, dashed] (8.0, 1.7) -- (12.0, -2.4);
			
			\draw[supersol, dashed] (12.0, -1.8) -- (14.0, 3.2);
			\draw[sol] (12.0, -2.2) -- (14.0, 2.8);
			\draw[subsol, dashed] (12.0, -2.4) -- (14.0, 2.6);
			
			\draw[supersol, dashed] (14.0, 3.2) -- (15, -2.7);
			\draw[sol] (14.0, 2.8) -- (15, -3.1);
			\draw[subsol, dashed] (14.0, 2.6) -- (15, -3.3);
			
			\draw[supersol, dashed] (15, -2.7) -- (15.1, -1.42);
			\draw[sol, dashed] (15, -3.1) -- (15.1, -1.82);
			\draw[subsol, dashed] (15, -3.3) -- (15.1, -2.02);
			
			\draw (-1.55,4.15) rectangle (3.05,2.5);
			\node[supersol,right] at (-1.5, 3.8) {$w$ {\color{black}is a supersolution}};
			\node[sol,right] at (-1.5, 3.3) {$u$ {\color{black}is a solution}};
			\node[subsol,right] at (-1.5, 2.8) {$v$ {\color{black}is a subsolution}};
			
		\end{tikzpicture}
	\end{figure}
	
	\section*{Acknowledgments} \label{sec::acknowledgments}	
	Part of this project was carried out while SD, EV, and RV was visiting Columbia University. For this reason, SD, EV, and RV would like to thank Columbia University for the hospitality.
	
	SD, EV, adn RV are members of the Australian Mathematical Society (AustMS).
	SD and RV are supported by the Australian Research Council Future Fellowship FT230100333 ``New perspectives on nonlocal equations''.
	OS is supported by NSF Grant DMS-2349794.
	EV is supported by the Australian Laureate Fellowship FL190100081 ``Minimal surfaces, free boundaries and partial differential equations''.
	RV is supported by a Scholarship for International Research Fees at the University of Western Australia.
	
	\printbibliography
	
\end{document}